\documentclass{amsart}
\usepackage{amsmath,amssymb,amsthm}
\usepackage[rokicki]{boxedeps}
\HideDisplacementBoxes

\theoremstyle{plain}

\newtheorem{prop}{Proposition}[section]
\newtheorem{thm}[prop]{Theorem}
\newtheorem{lem}[prop]{Lemma}
\newtheorem{cor}[prop]{Corollary}

\theoremstyle{definition}
\newtheorem{defn}[prop]{Definition}

\numberwithin{equation}{section}
 
\theoremstyle{remark}
\newtheorem{example}[prop]{Example}

\def\Z{\Bbb Z}

\def\R{\Bbb R}

\def\C{\Bbb C}

\def\flabel#1{\ifmmode #1\else$ #1$\fi}

\let\angle\undefined
\newsymbol\angle 115C
\def\degrees{\ifmmode^\circ\else$^\circ$\fi}

\def \pict #1 by #2 (#3) {\centerline{
\vbox to #2 {\hrule width #1 height 0pt depth 0pt
\vfill{\special{picture #3 }}}}}

\def \picture #1 by #2 (#3 scaled #4) #5{
\dimen0=#1 \dimen1=#2
\divide\dimen0 by 1000 \multiply\dimen0 by #4
\divide\dimen1 by 1000 \multiply\dimen1 by #4
\vbox{\pict \dimen0 by \dimen1 (#3 scaled #4)
\centerline { #5}}}

\catcode`\@=11

\title[]{Harmonic maps $\bf{M^3 \rightarrow S^1}$ and 2-cycles,
realizing \hfil\break the Thurston norm}

\begin{document}

\author{Gabriel Katz}
\address{Bennington College, Bennington, VT 05201-6001}

\email{gabrielkatz@rcn.com \quad \& \quad gkatz@bennington.edu}

\maketitle

\begin{abstract} 

Let $M^3$ be an oriented 3-manifold. We investigate when one of the fibers or a 
combination of fiber components, 
$F_{best}$, of a \emph{harmonic} map $f: M^3 \rightarrow S^1$ with Morse-type 
singularities delivers the Thurston norm $\chi_-([F_{best}])$ of its homology class 
$[F_{best}] \in H_2(M^3; \Z)$. 

In particular, for a map $f$ with connected fibers and any well-positioned oriented surface 
$\Sigma \subset M$ in the homology class of a fiber, we show that the Thurston number $\chi_-(\Sigma)$
satisfies an inequality
$$\chi_-(\Sigma) \geq \chi_-(F_{best}) - \rho^\circ(\Sigma, f)\cdot Var_{\chi_-}(f).$$
Here the variation $Var_{\chi_-}(f)$ is can be expressed in terms 
of the $\chi_-$-invariants of the fiber components, and the twist $\rho^\circ(\Sigma, f)$ measures 
the complexity of the intersection of $\Sigma$ with a particular set $F_R$ of "bad" fiber components. 
This complexity is tightly linked with the optimal "$\tilde f$-height" of $\Sigma$, 
being lifted to the $f$-induced cyclic cover $\tilde M^3 \rightarrow M^3$. 

Based on these invariants, for any Morse map $f$, we introduce the notion of its 
\emph{twist} $\rho_{\chi_-}(f)$. We prove that, for a harmonic $f$,\, 
$\chi_-([F_{best}]) = \, \chi_-(F_{best})$, if and only if,  $\rho_{\chi_-}(f) = 0$. 

\end{abstract}

%\tableofcontents
%\newpage

\bigskip

%%%%%%%%%%%%%%%%%%%%%%%%%%%%%%%%%%%%%%%

\section{Introduction}

Let $M$ be a compact, oriented 3-manifold, possibly with a boundary. 
With any homology class 
$[\Sigma]\in H_2(M, \partial M; \Z )$, one can associate a number of 
interesting invariants. The first one, $g([\Sigma])$, is the minimum genus 
of an embedded (immersed) oriented surface 
$(\Sigma, \partial\Sigma)  \subset (M, \partial M)$ 
realizing the homology class. 

Let $\Sigma^\oslash$ stand
for the union of all components of $\Sigma$, excluding spheres and disks. 
Put $\chi_-(\Sigma) = |\chi(\Sigma^\oslash)|$ and $\chi_+(\Sigma ) = 
|\chi(\Sigma \setminus \Sigma^\oslash)|$, where $\chi(\sim)$ is the Euler 
number.

The  the Thurston norm $\| [\Sigma] \|_T$ is defined to be the minimum of 
$\{\chi_-(\Sigma)\}$ over all embedded surfaces $\Sigma$ 
representing the homology class $[\Sigma]$. 
\smallskip

In general, the correspondence 
$\{[\Sigma] \Rightarrow \chi_-([\Sigma])\}$ gives rise 
to a \emph{semi-norm}  on the 
vector space $H_2(M, \partial M; \R ) \approx H^1(M; \R )$. In many cases, 
$\|\sim \|_T$ is actually a norm,
with the unit ball in the shape of a convex polyhedron.
\smallskip

Let $\mathcal F$ be a codimension one, oriented foliation with no Reeb components in 
$M$ or along its boundary. Such foliations are characterized by a global 
property: every leaf of $\mathcal F$ is hit by a loop \emph{transversal} to the foliation,
and a similar transversal loop condition is satisfied by $\mathcal F|_{\partial M}$.
In an appropriate metric, the leaves of  $\mathcal F$ 
 are \emph{minimal hypersurfaces} [Su].  Foliations with this property are 
called \emph{taut}.

In [T] Thurston showed that any compact leaf of a taut foliation
$\mathcal F$  attains the \emph{minimal} value of $\chi_-(\sim )$ in its homology class.
On the other hand, Gabai proved that, if a surface $\Sigma \subset M$ is minimizing 
the $\chi_-$-value in its non-trivial homology class and has no toral components, 
then it is a compact leaf of a smooth, taut foliation [G]. Thus, surfaces which realize 
the  norm $\|\sim \|_T$, are compact leaves of taut foliations. 
\smallskip

In general, taut foliations are hard to construct. In contrast, closed, or even
harmonic differential forms are easy to produce. 
If an oriented foliation $\mathcal F$ is generated by the kernels of a closed, 
\emph{non-singular} 1-form 
$\omega$, the foliation is automatically taut. Then, all the leaves of $\mathcal F$
are non-compact, or 
alternatively, they all are compact. In the second case, $M$ fibers over a circle, 
and $\mathcal F$ is comprised of the fibers of the corresponding map 
$f_\omega : M \rightarrow S^1$. In fact, $\omega$ and $f_\omega$ are \emph{harmonic} in 
an appropriate metric. In this setting, the harmonicity of a 1-form is equivalent to the 
tautness of the associated foliation. As a result, compact leaves of a foliation 
generated by a \emph{harmonic non-singular} form, realize the Thurston norm of their homology class.  
\smallskip

Closed 1-forms with singularities produce \emph{singular} foliations, which exhibit rich 
and drastically different behavior from the classical non-singular species [FKL]. 
In the paper, we will be concerned with the foliations generated by \emph{harmonic} 
1-forms with the \emph{Morse-type 
singularities} and \emph{rational} periods. Although their topology is very different 
from the non-singular foliations, they still possess the transversal loop property 
[C], and  their leaves are \emph{near}-minimal [K] (the harmonically-generated 
singular foliations are "near-taut").
\smallskip

There is a homological version of harmonicity (described in Theorem 4.7) which 
plays an important role in our arguments. Maps with one connected fiber are 
intrinsically harmonic.
\bigskip

Computing the Thurston norm $\|[\Sigma]\|_T$  in terms 
of the topology of $M$ alone is very difficult. The idea is to employ an appropriate 
map $f : M \rightarrow S^1$ to get a handle on the problem. It is natural to 
start with maps whose fibers realize $[\Sigma]$. The geometry of $f$ allows us 
to determine the $\chi_-$-invariant of each fiber component. At least, among all the 
combinations of fiber components, we can pick a representative of  $[\Sigma]$ with 
the minimal value of $\chi_-(\sim)$. This gives rise to an "$f$-vertical" 
semi-norm $\|\sim\|_{H^f}$ on the subspace  $H_2^f \subset H_2(M, \partial M; \R)$ spanned by the 
fundamental classes of various fiber components. We call a combination of fiber 
components which delivers $\|[\Sigma]\|_{H^f}$ the \emph{best} and denote it $F_{best}$.

Our main goal is to understand 
the relation between the  "incomputable" norm $\|[\Sigma]\|_T$ and the "computable" $\|[\Sigma]\|_{H^f}$. 
For instance, how to tell when a map $f$ has the property $\|[\Sigma]\|_T = \|[\Sigma]\|_{H^f}$?
A  somewhat different  question can be investigated: "When $\Sigma$ 
delivering $\|[\Sigma]\|_{H^f}$ is  realizable by a genuine $f$-fiber, and
not by a union of fiber components?" Figure 1 shows a map which is not intrinsically harmonic 
and for which $F_{best}$---a union of two spherical fiber components---is distinctly different from
any fiber. Answering both questions will allow
us  to characterize maps for which a fiber delivers the Thurston norm.
\bigskip 

This article is a by-product of my unsuccessful attempts to establish an analog of 
the Thurston Theorem for harmonically-generated foliations \emph{with singularities}.  
An important case of such foliations is provided by generic harmonic maps to a circle, 
that is, by generic rational harmonic 1-forms on $M$. \smallskip

For some time, I believed that, for a harmonic map 
$f : M \rightarrow S^1$, $\|[\Sigma]\|_T = \|[\Sigma]\|_{H^f}$ --- the best union of fiber 
components realizes 
the Thurston norm of its homology class. All my efforts to prove this very naive conjecture 
(by employing the theory of minimal surfaces) failed, until I found a simple 
counter-example  (cf. Example 4.10 ---the \emph{Harmonic Twister}).  Although too 
weak on its own, some form of  harmonicity seems to be a valuable ingredient in any 
"Best Fiber Component Theorem" (cf. Corollaries 8.3, 8.8): we always assume that our maps 
$f$ have no local extrema.\smallskip 

In fact, the reality is as far from what I 
conjectured as it could be: there are harmonic  maps $f : M \rightarrow S^1$ with 
very few singularities and with the $\chi_-(F_{best})$ arbitrary 
distant from the Thurston norm $\chi_-([F_{best}])$ (cf. Example 4.10). 

The phenomenon occurs because  maps can have arbitrary big "twists" 
$\rho_{\chi_-}(f)$.  Crudely, the twist invariant $\rho_{\chi_-}(f)$ measures 
the minimal complexity of the intersection patterns of  surfaces $\Sigma \subset M$,  
delivering the Thurston norm, with a generic
fiber component.  When $[\Sigma]$ is in the homology class of a fiber, $\rho_{\chi_-}(f)$
can be estimated  from above by the minimal $\tilde f$-\emph{height}  of such 
a $\Sigma$, being appropriately lifted to the $f$-induced cyclic covering 
$\tilde M \rightarrow M$.
Here $\tilde f$ stands for a function on $\tilde M$ covering the map $f$. 

In a sense, one  can also think about $\rho_{\chi_-}(f)$ as the $S^1$-controlled 
size of a homotopy, which takes a given map $f$ to a map $f_1$ with one of its fibers 
delivering the Thurston norm of its cohomology class (cf. Corollary 6.15).
\smallskip  

\smallskip

Similar invariants can be introduced for any probe surface $\Sigma \subset M$ 
in a vertical homology class  $[\Sigma] \in H_2^f$ (cf. Section 6). They are based on the twists  
$\rho(\Sigma, F)$ which measure the complexity of  
the intersection pattern $\mathcal C := \Sigma \cap F$ inside a generic fiber component $F$. 
The number $\rho(\Sigma, F) + 1$ does not exceed the number of components, in which $\mathcal C$ divides 
the surface $F$. By ignoring the components of $\mathcal C$ which bound a disk in $F$, a modification 
$\rho^\circ(\Sigma, F)$ of $\rho(\Sigma, F)$ is introduced.

In the Introduction, we use $\rho^\circ(\Sigma, f)$ to denote 
$max_F\, \{ \rho^\circ(\Sigma, F)\}$, where $F$ runs over all possible
$f$-fiber components.

When $[\Sigma]$ is chosen to be the homology class of a fiber $F$, 
the quantity $\rho(\Sigma, F)$ admits an interpretation as 
the \emph{breadth} $b(F, \Sigma)$  of a  lifting $\hat F \subset \tilde M$ relative to 
a special lifting $\hat\Sigma \subset \tilde M$ of $\Sigma$ (cf. Definition 6.7).  It has an 
upper bound $h(\Sigma, f)$ which is defined to be the integral part 
of the minimal $\tilde f$-\emph{height} of $\hat\Sigma$ plus one. 
\bigskip

Given a Morse map $f : M \rightarrow S^1$,
we consider a finite distribution of the values 
$\{{\chi_-}(f^{-1}(\theta))\}_{\theta \in S^1}$  
 along the circle.  A number ${\chi_-}(f^{-1}(\theta))$ can  
jump only when $\theta$ crosses an $f$-critical value. 
We define the ${\chi_-}$-\emph{variation} of the function $\{\theta \rightarrow {\chi_-}(f^{-1}(\theta))\}$ 
by the formula  (3.1).
For maps $f$ with no local extrema, the variation counts, so called, 
\emph{non-bubbling}  $f$-critical points (cf. Definition 3.2).

Of course, for any fibration $f$, $var_{{\chi_-}}(f) = 0$. 
However, if a harmonic $f$ is not a fibration, then $var_{\chi_-}(f) = 0$ 
implies that the Thurston \emph{semi}-norm is not a norm: some non-trivial 
classes in $H_2(M; \partial M;\, \Z)$ are represented by 2-spheres or 
2-disks.
\smallskip

We promote here a slogan: \hfil\break  
\emph{"Maps $f$ with the $0$-variation are like fibrations over the circle"}.
\bigskip

Finally, a  ghost of our original \emph{Best Fiber Conjecture} 
"${\chi_-}([F_{best}])  \; = \; {\chi_-}(F_{best})$" pays a visit:

\begin{thm}
Consider a Morse map $f : M \rightarrow S^1$ with all its 
fibers being connected\footnote{Such maps $f$ are intrinsically harmonic.}.  Then
\begin{itemize}
 
\item the ${\chi_-}$-invariant of the best fiber $F_{best}$ attains the 
minimal value among all the surfaces homologous to a fiber, if and only if,  
the twist $\rho_{{\chi_-}}(f) = 0$. 

\item In fact,   $var_{{\chi_-}}(f) = 0$ implies $\rho_{{\chi_-}}(f) = 0$.
\end{itemize}
Moreover, the same conclusions are  valid for connected sums of maps 
with connected fibers.  
\end{thm}

This theorem is a very special case of our main results---Theorems 8.2, Corollary 8.3, Theorem 8.6, 
Corollary 8.7 and Corollary 8.14. 
To avoid technicalities, let us state these propositions for another special, but important class 
of maps to a circle---for the \emph{self-indexing} harmonic maps (in fact, one can deform any 
map $f : M \rightarrow S^1$ into an intrinsically harmonic self-indexing map).

Given such a map $f$ and any "probe" surface 
$\Sigma \subset M$, homologous to a fiber and well-positioned (cf. Definition 7.1) with respect to the 
"worst" fiber $F_R$,
\begin{equation}
{\chi_-}(\Sigma)\; \geq \; {\chi_-}(F_{best}) - 
var_{{\chi_-}}(f)\cdot \rho^\circ(\Sigma, f), 
\end{equation}
where $var_{{\chi_-}}(f)$ is the number of non-bubbling  $f$-critical points. \smallskip
 
Evidently, unless the defect $var_{\chi_-}(f)\cdot \rho^\circ(\Sigma, f)$ is smaller than 
${\chi_-}(F_{best})$, this inequality is not very informative. \smallskip 

For harmonic self-indexing maps with a non-zero variation, (1.1)  can be also
viewed as  giving some grip of the invariants 
$\rho_{{\chi_-}}(f), h_{{\chi_-}}(f)$ (cf. Section 6): 
\begin{equation}
h_{{\chi_-}}(f)\; \geq \; \rho_{{\chi_-}}(f) \; \geq \; 
\frac{{\chi_-}(F_{best}) - {\chi_-}([F_{best}])}{var_{\chi_-}(f)}.
\end{equation}

In other words, if ${\chi_-}(F_{best}) \gg {\chi_-}([F_{best}])$
and the number of $f$-critical points is small, then the ${\chi_-}$-minimizing, 
well-positioned surface $\Sigma$
must be very tall. That is, $f$ "wraps" $\Sigma$ many times 
around the circle. Also, such $\Sigma$  must have a complex intersection patterns 
with a  fiber --- an intersection which is comprised of at least as many curves 
as the RHS of (1.2) requires. 
\bigskip

Now we describe the organization of the paper. It is  comprised of nine sections 
(including the Introduction) followed by a Notation List. \smallskip

In \emph{Section 2} we consider \emph{intrinsically harmonic} 1-forms
and maps into a circle. One can associate a finite graph $\Gamma_\omega$ with any 
closed 1-form $\omega$. The intrinsically 
harmonic forms and maps give rise to very special graphs. As a result, it possible 
to express faithfully intrinsic harmonicity in pure combinatorial, graph-theoretical terms.
\smallskip

In \emph{Section 3} we use graphs $\Gamma_f$ as book-keeping devices to
record the distribution of  $\chi_-$-invariants of the $f$-fibers and their 
connected components. The section deals with the effects of deforming a given 
map $f$ on the graph $\Gamma_f$ and these invariants.\smallskip

In \emph{Section 4} we develop further graph-theoretical manifestations 
of harmonicity (cf. Theorem 4.7). 
Examples of harmonic maps $f$, with only two singularities and with the norm
$\|[F_{best}]\|_{H^f}$ being \emph{arbitrary distant} from the  
Thurston norm $\|[F_{best}]\|_T$, conclude Section 4. In these examples, which 
we call \emph{Harmonic Twisters},  surfaces $\Sigma$, which deliver the 
Thurston norm, have arbitrary big twists and heights (cf. Figure 15).
\smallskip

In \emph{Section 5} we study a very special case of \emph{self-indexing} maps 
$f: M \rightarrow S^1$ and surfaces $\Sigma$'s, which have the
simplest  intersection pattern with the "worst" fiber $F_{worst}$. For them, we 
establish the most desirable result: $\chi_-(\Sigma) \geq \chi_-(F_{best})$, 
$g(\Sigma) \geq g(F_{best})$ (cf. Theorems 5.2). 
Section 5 indicates the main ideas of our approach in a form which is divorced 
from the combinatorial complexities of the general case (presented in Section 8). 
After 2-surgery, the resolved surface is pushed into a neighborhood of $F_{best}$---a union 
of "good" fiber components---, where it can be  effectively compared with the $F_{best}$.
\smallskip

\emph{Section 6} is the most tedious of them all. Here we develop the main technical 
tools: the notions of twist, breadth and height invariants of special surfaces 
$\Sigma \subset M$ in relation to a given map $f: M \rightarrow S^1$. Ultimately, 
the singularities of $f$ are responsible for the the non-triviality of these 
invariants.\smallskip

In \emph{Section 7} we aim to separate a generic embedded surface $\Sigma$, 
representing a given vertical homology class, from the union $F_R$ of the "worst" fiber components. 
In a sense, such a separation will permit us to reduce the case of general maps $f$
to the case, treated in Section 5. It is achieved by \emph{resolving} the 
intersections of a probe surface $\Sigma$ with the fiber components from
$F_R$ (cf. Figure 11). 
\smallskip

\emph{Section 8} contains the proofs of our main results---Theorems 8.2, 8.6, 8.13 and 
Corollaries 8.3, 8.7, 8.8, 8.14. Here we combine the strategy from Section 5 with the 
combinatorial tools developed in Sections 4 and 6 to derive generalizations 
of the inequalities (1.1) and (1.2).
After proving a variety  of "best fiber component theorems" generalizing Theorem 1.1,  we
proceed  to apply these results the  $\chi_-$-characteristic of special links 
(cf.  Corollary 8.15).
\smallskip
 
Finally, \emph{Section 9} deals with the way a surface $\Sigma$,  homologous to
a best combination $F_{best}$ of fiber components \emph{in the complement} to 
the $f$-singularities, can be tangent to the $f$-fibers.  In a sense, we connect, 
via the twist invariants, the Morse theory of $f$ with the induced Morse theory of 
$f|_\Sigma$ on a probe surface $\Sigma \subset M$.
\smallskip

To help our reader to cope with the expanding variety of notations, 
we conclude with a Notation List. \bigskip  

\emph{Acknowledgment:}\quad This work is shaped by numerous, thought-provoking 
conversations with Jerry Levine. I am thankful for his generosity and support. 
I am also grateful to the referee whose suggestions significantly improved the 
clarity of the original presentation.     
  
%%%%%%%%%%%%%%%%%%%%%%%%%%%%%%%%%%%%%%%   

\section{Intrinsically harmonic 1-forms and their graphs}

Next, we proceed with a description of a few facts, constructions and notations 
related to an intrinsic characterization of harmonic (rational) 1-forms. Actually, 
these facts are not specific to dimension three.  

Let $\Sigma \subset M$ be an oriented surface and 
$[\omega] \in H^1 (M; \Z )$ be the class Poincar\'{e}-dual to $[\Sigma]$.
It can be 
realized by a closed rational 1-form $\omega$ on $M$, or, equivalently, by a map 
$f: M \rightarrow S^1$ with $\Sigma$ as one of its regular fibers. 
The two realizations are linked by the formula 
$f^\ast (d\theta ) = \omega$, where $d\theta$ is the canonical 1-form on the
oriented 
circle. For a given class $[\omega ]$, one can choose its representatives $f$ 
and $\omega$ with Morse-type singularities. Furthermore, if $[\omega ] \neq 0$, 
through a deformation of $f$, the singularities of indices 0 
and 3---the local minima and  maxima of $f$---can be eliminated. 
\smallskip

By considering harmonic maps $f: M \rightarrow S^1$ or, what is the same,
harmonic 1-forms $\omega = f^\ast(d\theta)$, we exclude the singularities of indices 0 and 3---harmonic 
functions have no local maxima and minima. One might 
wonder, if there are  restrictions on the distribution of critical points of 
indices 1 and 2, imposed by the harmonicity and which are not prescribed by the 
topology. Fortunately, this question has a comprehensive answer.

\begin{figure}[ht]
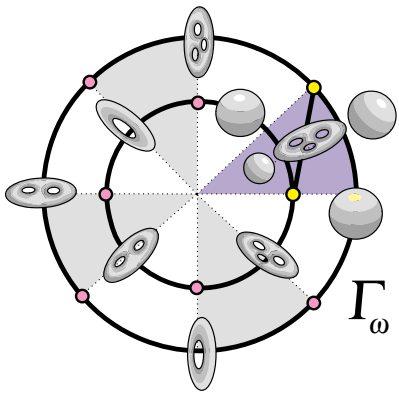

\centerline{\BoxedEPSF{Graph scaled 450}} 
\bigskip
\caption{Graph $\Gamma_\omega$ and the fiber components}
\end{figure}

In [FKL], for any closed 1-form $\omega$ on $M$, we have introduced a finite
graph $\Gamma_\omega$.  
The manifold $M$ is canonically mapped onto the graph $\Gamma_\omega$ by a map 
 $p_\omega$. When the singular foliation $\mathcal{F}_\omega$, determined by
$\omega$, has only compact leaves, 
the points of $\Gamma_\omega$ are just their \emph{connected components}.  
In particular, each connected component of every fiber of a Morse
mapping $f: M \rightarrow S^1$ corresponds to a single point in the graph
$\Gamma_f = \Gamma_\omega$. Vertices of $\Gamma_f$ correspond to critical 
points of $f$. 

An example of $\Gamma_f$ is given in Figure 1, which also 
depicts generic fiber components suspended over each edge of $\Gamma_f$. 

The map $f$ factors through a canonical 
projection $p_f: M \rightarrow \Gamma_f$  and thus, generates an equally canonical map 
$\pi_f : \Gamma_f \rightarrow S^1$.  In Figure 1, $\pi_f$ is 
induced by the radial projection.

For Morse maps $f$ with distinct critical values in $S^1$ and no local extrema, 
some vertices of 
$\Gamma_f$ mark critical points $x_\star$, such that crossing the 
critical value $\theta_\star = f(x_\star)$, causes the fiber to change 
the number of its connected components. Such singularities correspond to the 
\emph{trivalent} vertices 
in $\Gamma_\omega$. The rest of the singularities correspond to \emph{bivalent} 
vertices of the graph. 
\begin{defn}
A 1-form $\omega$ or a related map $f: M \rightarrow S^1$ are called \emph{intrinsically harmonic}, 
if they are harmonic with respect to \emph{some} metric on $M$. \smallskip
\end{defn}
In [C], Calabi established the following global criterion:
\begin{thm}
A closed 1-form $\omega$ on a \emph{closed} manifold $M$ is
intrinsically harmonic, if and only if, through any point in $M$, different from the 
singularities of $\omega$, there is a loop $\gamma$, along which $\omega$ is strictly \emph{positive},
that is, $\omega(\dot\gamma) > 0$.
\end{thm} 

In [FKL], we observed that a similar property can be formulated for the graph
$\Gamma_\omega$. In fact, for closed manifolds, the intrinsic harmonicity of $\omega$ becomes
equivalent to the following \emph{positive loop property} of $\Gamma_\omega$: 
through any point of $\Gamma_\omega$, one can draw an oriented loop, comprised of $\omega$-oriented edges. 
We call such $\Gamma_\omega$'s \emph{Calabi Graphs}.

Thus, the notion of the Calabi graphs provides us with a completely combinatorial 
description of the intrinsic harmonicity.\smallskip

In particular, if for 
some arc in $S^1$, its pre-image in $\Gamma_f$ is a \emph{single} edge, 
the graph automatically satisfies the positive loop property, and the map $f$ to 
$S^1$ is intrinsically harmonic.  

The graph in Figure 1 violates the 
positive loop property and the corresponding map $f$ is not intrinsically harmonic.
\smallskip

In [FKL], we proved  Theorem 2.3 below for two extreme model cases:   
the case when all the leaves of the foliation $\mathcal{F}_\omega$ are 
compact and the case, when none of the leaves is compact. In  [Ho], 
Honda established the general case.

\begin{thm} Let $\omega$ be a closed 1-form on a closed $n$-manifold $M$ with
the Morse-type singularities. Assume that  $\omega$ has no critical points of
indices 0 and $n$. Then one can deform $\omega$ (through the space of closed 1-forms) 
to an intrinsically harmonic form $\tilde{\omega}$, which has the same
collection of singularities. \qed 
\end{thm} 

%%%%%%%%%%%%%%%%%%%%%%%%%%%%%%%%%%%%%%%  

\section{The Genera of Fibers and the Combinatorics of Handle Moves} 

We examine some graph-theoretical descriptions of Morse maps $f: M \rightarrow S^1$ 
in connection to genera and $\chi_-$-invariants of their fiber components. We analyze how 
$f$-deformations affect these combinatorial descriptions.\bigskip 

The genus $g(\Sigma )$ of a surface $\Sigma$ is defined to be half of the rank
of the homology group $H_1 (\Sigma ; \Z )$\footnote{for closed surfaces it is an integer, for surfaces
with boundary, it might be a half of an integer.}. 
If $\Sigma$ consists of several components, 
its genus is the sum of the components' genera.  The
Euler characteristic $\chi(\sim )$ of a surface is  the sum of the
Euler characteristics of its components. Finally, the Thurston's 
$\chi_- (\sim )$-characteristic of a surface is the
\emph{absolute value} of 
the  Euler number of the union of all its components, \emph{excluding}
the 2-spheres 
and, in the case of surfaces with boundary, also excluding the 2-disks. 

For a closed oriented surface $\Sigma$,\, $\chi_-(\Sigma) = 2|\nu(\Sigma) -
g(\Sigma)|$, where 
$\nu(\Sigma)$ stands for the number of non-spherical connected components in
$\Sigma$.
\smallskip

Let $\gamma$ be a simple loop on an oriented surface $\Sigma$. Performing a
2-surgery on  $\Sigma$ along $\gamma$ has the following effect on the 
three invariants. \smallskip 

{\bf  2-Surgery List A:}
\begin{itemize}
\item When the loop $\gamma$ is null-homotopic, then 
the surgery will have the following effect:
\begin{enumerate} 
\item the  genus will remain the same, 
\item the  Euler characteristic will increase by $2$, 
\item the  $\chi_-$ will remain the same;
\end{enumerate} 

\item When the loop $\gamma$ separates  $\Sigma$ and is not null-homotopic, then 
the surgery will have the following effect:
\begin{enumerate} 
\item the  genus will remain the same,  
\item the  Euler characteristic will increase by $2$, 
\item the  $\chi_-$ will decrease by $2$;
\end{enumerate}  

\item When the loop $\gamma$ does not separate  $\Sigma$, then 
the surgery will change: 
\begin{enumerate} 
\item the  genus by subtracting $1$, 
\item the  Euler characteristic by adding $2$, 
\item the  $\chi_-$ will remain the same, if the
surface 
is a torus, and will decrease by 2 otherwise.
\end{enumerate}
\end{itemize}

Let $\gamma$ be a simple arc on an oriented surface $\Sigma$, connecting two
points on its boundary $\partial\Sigma$. Performing a  relative 2-surgery 
on  $\Sigma$ along $\gamma$ produces  changes described in 
\smallskip

{\bf 2-Surgery List B:}
\begin{itemize}
\item When the arc $\gamma$ is null-homotopic modulo $\partial\Sigma$, then 
the surgery will have the following effect:
\begin{enumerate} 
\item the  genus will remain the same, 
\item the  Euler characteristic will increase by $1$, 
\item the  $\chi_-$ will remain the same;
\end{enumerate} 

\item When the arc $\gamma$ separates  $\Sigma$ and is not null-homotopic 
modulo $\partial\Sigma$ , then 
\begin{enumerate} 
\item the  genus will remain the same,  
\item the  Euler characteristic will increase by $1$, 
\item the  $\chi_-$ will decrease by $1$;
\end{enumerate}  

\item When the arc $\gamma$ does not separate  $\Sigma$, then 
\begin{enumerate} 
\item the  genus will drop by $1$, 
\item the  Euler characteristic will increase by $1$, 
\item the  $\chi_-$ will remain the same, if the
surface is an annulus, and will decrease by 1 otherwise.
\end{enumerate}
\end{itemize}

These observations can be summarized in 

\begin{lem} Under 2-surgery, the genus of a surface and its 
$\chi_- (\sim )$-characteristic are  non-increasing quantities. The 
Euler characteristic is strictly increasing. \qed
\end{lem}
\smallskip

In the vicinity of each vertex and over a small arc of the circle centered on 
the critical value, $\Gamma_f$ looks as depicted in the four diagrams of Figure 2. 
The labels 1 and 2 indicate the Morse index of the critical 
point. Each vertex is a \emph{bivalent} or a \emph{trivalent} one. 
The corresponding critical points also are called bivalent or trivalent. 
The diagrams do not include the case of critical points of indices 0 and 3. 
They can be depicted  by an oriented edge, emanating from or terminating 
at a vertex of multiplicity 1.
%%%%%
\begin{figure}[ht]
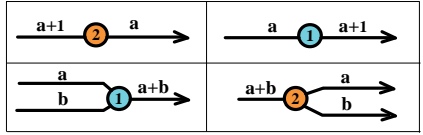

\centerline{\BoxedEPSF{3-valent.g scaled 300}}
\bigskip
\caption{Bivalent and trivalent singularities and the 1-chain $\tau_g(f)$.}
\end{figure}
%%%%
Using $f$, one can 
produce a 1-chain $\tau_g(f)$ on the graph $\Gamma_f$: just assign to each edge
the genus of 
the generic fiber component over it. Examining the surgery lists $A$ and $B$ above, 
we see that 
the chain $\tau_g(f)$ is a \emph{relative cycle} modulo the bivalent vertices.
When $f: \partial M \rightarrow S^1$ is a fibration, crossing 
a bivalent vertex of index 1 from left to right, results in an increase of the 
$\tau_g(f)$-value by 1, and crossing a bivalent vertex of index 2 from left to
right, results in a decrease of the $\tau_g(f)$-value by 1. 
 
\smallskip

\begin{figure}[ht]
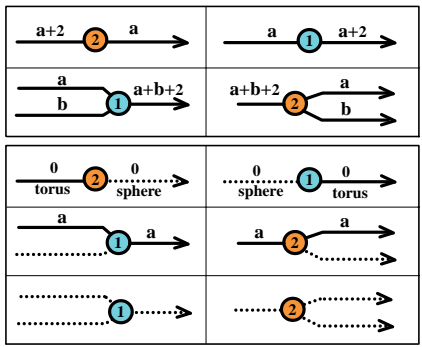

\centerline{\BoxedEPSF{3-valent.chi scaled 300}}
\bigskip
\caption{Bivalent and trivalent singularities and the 1-chain $\tau_{\chi_-}(f)$. 
The dotted lines mark the spherical fiber components.}
\end{figure}

For a given smooth map $f: M \rightarrow S^1$ 
with a finite number of critical values $\{\theta_i^\star \}_{0\leq i \leq k}$,
we define the variation  $var_{{\chi_-}} (f)$ as a \emph{half} 
of the cyclic sum
\begin{equation} 
|{\chi_-} (f^{-1}(\theta_0 )) - {\chi_-}(f^{-1}(\theta_k))| + 
\sum_{i = 0}^{k - 1}\; |{\chi_-}(f^{-1}(\theta_{i + 1})) - 
{\chi_-}(f^{-1}(\theta_i))|,
\end{equation}
where a typical point $\theta_i$ belongs to the open arc 
$(\theta_i^\star , \theta_{i + 1}^\star )$.
We introduce the oscillation $osc_{{\chi_-}}(f)$
\begin{equation} osc_{{\chi_-}}(f) := max_\theta \{{\chi_-}(f^{-1}(\theta_i))\} - 
min_\theta \{{\chi_-}(f^{-1}(\theta_i))\}\; 
\end{equation}

By definition, $var_{{\chi_-}}(f) \geq  osc_{{\chi_-}}(f)$. 

Note that, $\{var_{{\chi_-}}(f) = 0\}$ and $\{osc_{{\chi_-}}(f) = 0\}$  are 
equivalent conditions imposed on $f$.
Evidently, for a fibration $f$, 
$var_{{\chi_-}}(f) = osc_{{\chi_-}}(f) = 0$. 
\bigskip

We say that a Morse map $f: M \rightarrow S^1$ is \emph{self-indexing}, if there
is a point 
$\theta_b \in S^1$, so that, moving from $\theta_b$ along the oriented circle,
the critical 
values of critical points of lower indices precede the ones of higher indices.
Most of the time, we assume that the $f$-critical values are all distinct. 
\smallskip

Given a self-indexing map with no critical points 
of indices 0 and 3, one can find two points $\theta_b, \theta_w \in S^1$,
such that 
the oriented arc $(\theta_b, \theta_w)$ contains all critical values of
index 1 and  
the complementary arc $(\theta_w, \theta_b)$---all  critical values of index 2.
Let  $F_{best} = f^{-1}(\theta_b)$ and $F_{worst} = f^{-1}(\theta_w)$.
\smallskip

For a self-indexing map  $f: M^3 \rightarrow S^1$ with no local maxima and minima, 
the variation $var_{{\chi_-}}(f)$ equals
$osc_{{\chi_-}}(f) := {\chi_-} (F_{worst}) - {\chi_-} (F_{best})$:  the invariant 
${\chi_-}$ is non-increasing under 2-surgery.
\smallskip
\smallskip

Next, we examine the effect of deforming Morse maps $f: M \rightarrow S^1$
on the invariants ${\chi_-}(F_{worst}),  {\chi_-}(F_{best})$ and 
$var_{{\chi_-}}(f)$. 
We use  graphs $\Gamma_f$'s, equipped with a canonical map $\pi_f$ to a
circle, as book-keeping devices. 
\smallskip

We notice that the boundary $\partial \tau_g(f)$ of the 1-chain $\tau_g(f)$ is a 0-chain 
on $\Gamma_f$ supported on the bivalent vertices. Its $l_1$-norm is the  
variation $var_g(f)$. 
\smallskip

In a similar way, one can introduce an 1-chain $\tau_{\chi_-}(f)$ on  
$\Gamma_f$ by assigning to each edge
the $\chi_-$-characteristic of the corresponding fiber component. Figure 
3 deals with the case, when $M$ is closed, or when $\partial M \rightarrow S^1$ is 
a fibration. It is divided into generic and special  patterns. Special patterns 
arise when at least one of fiber components is a sphere.  
\smallskip
\begin{defn}
A critical point is called \emph{bubbling}, if there is at least one spherical
or disk fiber component in its vicinity. 
\end{defn} 
\smallskip

Note that only at the bubbling vertices the chain $\tau_{\chi_-}(f)$ satisfies 
the cycle condition. The boundary $\partial \tau_{\chi_-}(f)$ of the 1-chain 
$\tau_{\chi_-}(f)$ is a 0-chain on $\Gamma_f$ supported on the non-bubbling vertices. 
Its $l_1$-norm is the  variation $var_{\chi_-}(f)$. 
\smallskip
\smallskip

Deforming $f$ causes $\Gamma_f$ 
to go through a number of transformations that can be decomposed in a few basic
moves. Before and after deformations, all the $\pi_f$-images of the vertices in
$\Gamma_f$ are assumed to be distinct in $S^1$.  
 
The five diagrams in Figure 4 depict the effect on the chain $\gamma_g(f)$ of 
deforming a Morse function, 
so that the critical value of index 1 is placed \emph{below} the critical level
of index 2 (equivalently, of sliding a 1-handle "below" a 2-handle). 
The diagrams are produced by combining the patterns from Figure 2 in pairs. 

We notice that such an operation is always possible [M]. Unfortunately, it could
only \emph{increase} the $l_1$-norm of the chains $\tau_g(f)$, $\partial\tau_g(f)$, 
or the value $g(F_{best})$. From this perspective, the inverse operations 
(acting from the left to the right configurations) are desirable, but not always
geometrically realizable! 
%%%%%%%%%%%%
\begin{figure}[ht]
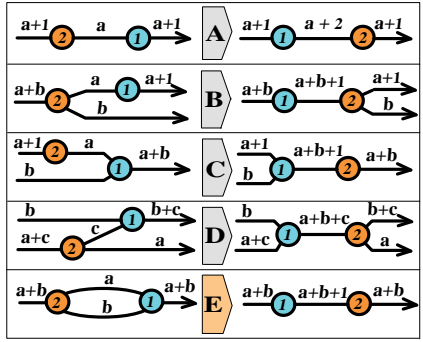

\centerline{\BoxedEPSF{5-moves scaled 300}}
\bigskip
\caption{Generic patterns of moving a critical point of index 2 
above a critical
point of index 1 and the affect of these moves on the genera of fiber components.}
\end{figure} 
%%%%%%%%%%%% 
On the other hand, any two critical points of the same index can be re-ordered. 
 
Consider a portion $W$ of $M$ represented by the diagrams in Figure 4 and view $W$ 
as a cobordism between two surfaces $\Sigma_0$ and $\Sigma_1$ represented by the 
left and right ends of the diagrams. Examining the five moves, we notice that, 
in the configurations $B$ through $E$, the two critical points can not cancel each 
other \emph{locally}, that is, by a deformation which is constant on 
$\Sigma_0 \cup \Sigma_1$---the trivial cobordism that would result
is inconsistent with: 1)  the prescribed connectivity of $\Sigma_0 \cup \Sigma_1$ (cases B and C) or 
2) with the connectivity of $W$ (case D), or with the non-triviality of $H_1(W, \Sigma_0)$ 
(cases D, E). As the diagrams testify, the local cancellation of an 1-handle and a 2-handle is 
only possible among the vertices in the diagram A. 
  
Note that, the contribution to the variation $var_g (\sim)$ in four diagrams 
$A$---$D$ remains invariant under the moves! In fact, the $var_g (\sim)$ can be 
changed only through the cancellation of singularities in pairs, or through the 
$E$-moves. In the first case it drops by 1, in the second --- rises by 1.
%%%%%
\begin{figure}[h]
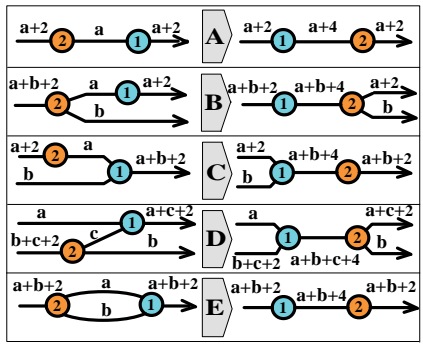

\centerline{\BoxedEPSF{5-moves.chi scaled 300}}
\bigskip
\caption{Generic patterns of moving a critical point of index 2 above a critical
point of index 1 and the affect of these moves on the $\chi_-$-characteristics of 
the fiber components.}
\end{figure}
%%%%%%%
Figure 5 shows the effect of the same deformations on the chain $\tau_{\chi_-}(f)$. 
It depicts cases, where the spherical and disk components are not involved. The diagrams are 
produced by combining in pairs the first four patterns in Figure 3.

The 14 bubbling (that is, special spherical) patterns are the result of combining in pairs 
the patterns in Figure 3. We leave their depiction to the reader.

Note, that the  variation $var_{\chi_-}(f) = \frac{1}{2}\|\partial\tau_{\chi_-}(f)\|_{l_1}$ 
is preserved under all the moves in Figure 5. It can increase by 2  only under 
special moves which involve bubbling singularities.
The variation $var_{\chi_-}(f)$ is non-decreasing under all the moves. 

The diagrams in Figures 2---5 reflect our fundamental assumption that 
$f: \partial M \rightarrow S^1$ is a fibration. 
\smallskip

These observations can be summarized in the  Lemmas 3.3---3.5 below. We assume that all 
the deformations, described in the lemmas, take place in the space $\mathcal M(M, S^1)$ 
of smooth maps $f: M \rightarrow S^1$ 
with \emph{no local maxima and minima} and with the generalized Morse singularities
(these are the usual quadratic Morse singularities and the generic cubic singularities, 
resulting from the merge of two Morse singularities). Let 
$\mathcal M^{\odot}(M, S^1) \subset \mathcal M(M, S^1)$ be the subspace of maps with 
no spherical or disk components in their fibers. Considering maps from 
$\mathcal M^{\odot}(M, S^1)$ eliminates the special bubbling patterns. 
For instance, if no non-trivial element from $H_2(M, \partial M; \Z)$ admits a 
representation by a sphere or a disk, then any harmonic map $f: M \rightarrow S^1$ belongs to 
$\mathcal M^{\odot}(M, S^1)$. 

\begin{lem} The number of bivalent, index-one vertices of the graph $\Gamma_f$ 
equals to the  variation $var_g (f)$ of the fiber genus, which, in turn, coincides 
with the half of the $l_1$-norm $\|\partial\tau_{g}(f)\|_{l_1}$ of the 0-cycle 
$\partial\tau_{g}(f)$.  

The numbers of trivalent and bivalent vertices in $\Gamma_f$ are deformation invariants 
of  the Morse maps, which
avoid the  the cancellation of singularities and the $E$-moves. Such deformations can be
decomposed into a sequence 
of the four basic moves $A$---$D$, depicted in Figure 4, and their inverses
(when those are realizable),
plus the moves which reorder the adjacent vertices of the same Morse 
index \footnote{they do not affect the 
graph $\Gamma_f$ and the chain $\tau_g(f)$.}.

Under an  $E$-move, the norms $\|\partial\tau_{g}(f)\|_{l_1}$ and 
$\|\tau_g(f)\|_{l_1}$ both increase by 2. 
\qed 
\end{lem}

\begin{lem} Under the cancellation of singularities of indices 1 and 2, the norms  
$\|\tau_g(f)\|_{l_1}$, $\|\tau_{\chi_-}(f)\|_{l_1}$ and the variations $\|\partial\tau_{g}(f)\|_{l_1}$, 
$\|\partial\tau_{\chi_-}(f)\|_{l_1}$ are  decreasing.  
\qed
\end{lem}

\begin{lem} The number of non-bubbling vertices 
of the graph $\Gamma_f$ equals to the $\chi_-$-variation 
$var_{\chi_-}(f) = \|\partial\tau_{\chi_-}(f)\|_{l_1}$.

The variation is preserved under the generic 
moves\footnote{those avoid bubbling singularities.}, as shown in Figure 5, 
together with their 
inverses. The  special moves (involving bubbling singularities) can increase 
the  variation by 2. 

In particular,  $var_{\chi_-}(f)$ is invariant under the deformations within 
the space of Morse maps with a fixed list of index 1 and 2 singularities and with 
no spherical fiber components.

Under the generic moves,  $\|\tau_{\chi_-}(f)\|_{l_1}$ increases 
at least by 4. Under the special moves, it increases by 2, or is preserved. 
\qed 
\end{lem}

%%%%%%%%%%%
%%%%%%%%%%%%

\section{Graph-theoretical manifestations of harmonicity} 

In this section we further examine harmonicity in terms of the graph theory. \smallskip 

Let $\pi_f: \Gamma_f \rightarrow S^1$ be a map of graphs corresponding to a given Morse map 
$f: M \rightarrow S^1$.  For a map $f$ with no local maxima and minima, the vertices of $\Gamma_f$
all are bivalent or trivalent.  In addition,  they come in two flavors: indexed by 1 or 2 
depending on the Morse index of the corresponding  critical point.\smallskip

There are a few restrictions on the distribution of vertices of indices 1 and 2 in $\Gamma_f$.  They are 
shown in Figure 6.  
%%%%%%
\begin{lem} Assume that $f:  M \rightarrow S^1$ has no local extrema. If a simple \emph{positive} loop 
$\gamma \subset \Gamma_f$ contains vertices only of a
particular index, then all the f-critical points corresponding to these vertices are bubbling. 

Any loop $\gamma \subset \Gamma_f$ which is not mapped by $\pi_f$ into $S^1$ in a monotone 
fashion, contains at least one  vertex of index 1 and at least one vertex of index 2. 
\end{lem}
%%%%%%%
%%%%%%
\begin{figure}[ht]
\centerline{\BoxedEPSF{rotary.same.index scaled 400}} 
\bigskip
\caption{}
\end{figure} 
%%%%%%%%
{\bf Proof.}\quad Follow the change in the $\chi_-$-values or genera of fiber components along $\gamma$ 
and use the principle: "what goes up has to come down".  Figures 2, 3 and 6 will facilitate the argument. 
Because all vertices of $\Gamma_f$ are bivalent or trivalent, each local maximum of 
$\pi_f: \gamma \rightarrow S^1$ corresponds to a trivalent vertex of index 1 and each 
local minimum --- a trivalent vertex of index 2 (cf. Figure 8). These local maxima and minima alternate along 
$\gamma$ and their cardinalities are equal.   \qed  
\bigskip

We introduce two finite subsets $A$ and $R$ of $\Gamma_f$ which will play an important role in the paper. 
The elements of $A$ will
be called \emph{attractors} and the elements of $R$ --- \emph{repellers}\footnote{The names are inspired by 
the roles these combinatorial devices will play in our method of tackling the Best Fiber Component Problem.}.  
Each  $\pi_f$-oriented edge of 
$\Gamma_f$ with its left vertex being of index 1 and its right vertex of index 2 acquires exactly one repeller;
each oriented edge  with its left vertex being of index 2 and its right vertex of index 1 acquires exactly one 
attaractor (cf.  Figure 7). 

We observe that changing $f$ to $-f\; mod.\, 2\pi$ switches the orientations of the edges in 
$\Gamma_f$ and turns points of index 1 into points of index 2. Therefore, $\Gamma_f$ and $\Gamma_{-f}$ 
share the same sets of attractors and repellers. 

It is worth noticing that the elementary moves from Figure 4 all increase the number of 
repellers by 1. So, it is easy to increase the size of $R$, to decrease it is a very different 
story. For instance, one might wonder: \emph{what is the minimal number of repellers in a given} 
2-\emph{homology class}? 
%%%%%%
\begin{figure}[ht]
\centerline{\BoxedEPSF{RAgraph scaled 520}} 
\bigskip
\caption{}
\end{figure} 
%%%%%%%%
\begin{lem} Let $\pi_f: \Gamma_f \rightarrow S^1$ be such that no positive loop in $\Gamma_f$ contains 
vertices only of a particular index and no vertex of $\Gamma_f$ is 
univalent\footnote{This is the case when $f: M \rightarrow S^1$ does not 
have extrema and bubbling singularities.}. 
Then each $r \in R$ gives rise to a pair of \emph{subtrees} $T_r^+, T_r^- \subset \Gamma_f$ with 
the common root at $r$ and their leaves belonging to $A$. The branches of $T_r^+$ ($T_r^-$) are formed by paths
emanating from $r$ in the positive (negative) direction and terminating at the first point from $A$ they 
encounter. The trees $\{T_r^+, T_r^-\}_{r\in R}$ form a cover of $\Gamma_f$.

In particular, these conclusions hold for any $\Gamma_f$ produced by a harmonic map $f$ of a manifold $M$ 
with the property: no non-trivial class in $H_2(M, \partial M; \Z)$ admits a spherical or disk representative.
\end{lem}
%%%%%%
{\bf Proof.}\quad Assume that the subgraph $T_r^+$ contains a loop. This could happen in a number of ways. 

1). There exists an $\pi_f$-positive path $\xi$ in $\Gamma_f$ which leaves $r$ and closes on itself at a vertex $x$ 
without encounter an attractor.  This generates a positive loop $\tau \subset \xi$ which contains $x$.  
If $\tau$ does not contain $r$, then $x$ must be of index 1.  By the lemma's hypotheses, $\tau$ must contain 
at least one vertex of index 2. Hence, an attaractor must exist on $\tau$.  This contradicts to assumption 1).
If  $r \in \tau$, $\xi$ is a positive closed path.  It  contains an oriented edge $[y, z]$, where the vertex 
$y$ is of  index 2, the vertex $z$ of index 1 and $r\in [y,z]$.  Therefore, the loop $\xi$ must also contain 
an oriented edge $[y', z']$, where the vertex $y'$ is of index 1 and the vertex $z'$  of index 2.  So, $\xi$ 
must contain an attractor, which  contradicts to assumption 1).
\smallskip

2). The second option for $T_r^+$ to contain a loop arises when there are two distinct positive paths 
emanating from $r$ and terminating at the same attractor $a$. We can assume that, for both paths,  $a$ is the first 
attractor after $r$. Since the two paths must first separate at 
some point $x$ of index 2 (which succeeds $r$) and then join at another point $y$ of index 1 (which precedes $a$), 
each of the paths must contain at least one attractor distinct from $a$ and which precedes it. Thus, the two paths 
must terminate at these two distinct attractors \emph{before} they reach $a$. This contradiction rules out loops of 
the second type. 

Finally, we need to show that any point $x \in \Gamma_f$  belongs to some tree  $T_r^+$ or $T_r^-$. Consider  
a positive  path $\xi$ through $x$ which does not admit any extension. Such a path must be closed or must 
close on itself in both positive and negative directions.  If $\xi$ is closed, by the lemma's 
hypothesis, it must contain at least one  vertex of index 1 and at least one vertex of index 2, unless $f$ is 
a fibration.  Therefore, $\xi$ 
will contain at least one attractor and one repeller. Thus, moving from $x$ along $\xi$ in the positive or 
negative directions we must encounter a repeller $r$. Evidently, $x \in T_r^\pm$ for the first repeller.  

The case when $\xi$ through $x$ closes on itself in both directions already has been analyzed in 1). Again, the two
loops at the "ends" of $\xi$ each must contain an attractor-repeller pair. In the worst case, at least there we 
will find the right repeller---a repeller $r$ whose trees $T_r^\pm$ contains $x$. 
\smallskip

It remains to notice that harmonic maps $f$ do not have local extrema, thus, excluding univalent vertices 
in $\Gamma_f$. Also, no bubbling singularities can occur, because the relevant spherical fiber components of 
harmonic maps must generate non-trivial elements in the 2-homology of $M$ (contrary to the hypotheses 
about $M$).
\qed 

\begin{lem}  For any loop $\gamma \subset \Gamma_f$,\, $\int_\gamma\; A = \int_\gamma\; R$, where 
$\int_\gamma\; A$ and $\int_\gamma\; R$ stand for the sums of $(\pm)$-weighted points from $A$ and $R$ along the 
loop $\gamma$. The sign of a point $x \in \gamma$ is produced by comparing the orientation of $\gamma$ with the 
$\pi_f$-induced orientation of the edge containing $x$ (cf. Figure 8).  
\end{lem}
%%%%%%
\begin{figure}[ht]
\centerline{\BoxedEPSF{int(A)=int(R) scaled 500}} 
\caption{\quad $\int_\gamma\; [A] = \int_\gamma\; [R]$}
\end{figure} 
%%%%%%%%
{\bf Proof.}\quad
The numbers of local maxima and minima of the function $\pi_f|_\gamma$ along any loop $\gamma$ are equal.  
Maxima occur at vertices of index 2 and minima occur at vertices of index 1. The signs attached 
to  repellers and attractors along each arc of $\gamma$ between two consequent extrema are the 
same and alternate as one crosses from an arc to an arc. Also, since along each arc attractors 
and repellers alternate and since the cardinality of attractors exceeds the cardinality of repellers by 1, 
the attractors contribute to the integral 1 more (less) than the repellers along any accenting (descending) arc. 
Because the number of accenting and descending arcs are equal, the total contributions of $R$ and $A$ to 
the integral are equal as well.  
\qed

\begin{defn} A homology class in $H_2(M, \partial M; \Z)$ is called $f$-\emph{vertical} if it 
can be represented as a $\Z$-linear combination of the fundamental classes of the fiber components
(which do not necessarily belong to the same $f$-fiber). 

We denote the subgroup of $f$-vertical classes
by $H_2^f$. It contains a positive cone $H_2^{f +}$ generated by non-negative linear combinations of the fundamental 
classes of the fiber components\footnote{The orientations of the fiber components is determined  via $f$ by the 
preferred orientations of $M$ and $S^1$.}.
\end{defn}

\begin{lem} The elements of $H_2^f$ are detected by their intersection numbers with loops 
$\{C_k \subset M\}_k$ whose images in $\Gamma_f$ form a basis of $H_1(\Gamma_f; \Z)$. 
In fact, $H_2^f$ is Poincar\'{e}-dual to the image of $H^1(\Gamma_f; \Z)$ in $H^1(M; \Z)$
induced by the canonic projection $p_f: M \rightarrow \Gamma_f$. 
\end{lem}
{\bf Proof.}\quad By the Poincar\'{e} duality, 
an element of $H_2(M, \partial M; \Z)$ is determined by its intersection numbers with loops in $M$. 
The intersection number
of a loop $C \subset M$ with a vertical 2-cycle $[\Sigma]$ equals to the intersection number of its 
weighted finite support $\sigma$ in $\Gamma_f$ with the image $C'$ of $C$ under the map $p_f: M \rightarrow \Gamma_f$.  
Therefore,  $\Sigma \circ C = \sigma \circ C' = \sigma^\ast(C')$, where $\sigma^\ast$ stands for the  
1-cocycle in $\Gamma_f$ dual 
to the 0-chain $\sigma$. Thus, $\Sigma \circ C$ reduces to the natural non-degenerated pairing between 
$H_1(\Gamma_f; \Z)$ and $H^1(\Gamma_f; \Z)$. Therefore, $H_2^f \approx p_f^\ast(H^1(\Gamma_f; \Z))$. \qed
\bigskip

Let $\Z[A]$ be the free $\Z$-module generated by the attracting set $A = \{a\}$. Elements of $\Z[A]$ can
be  viewed as functions $\kappa : A \rightarrow \Z$ or, equivalently, as formal sums $\sum_{a\in A}\; \kappa_a \cdot a$
with integral coefficients $\{\kappa_a\}$. Combinations with non-negative coefficients generate a positive cone 
$\Z_+[A] \subset \Z[A]$.

Under the Lemma 4.2 hypotheses, each point in $\Gamma_f$ which is not a vertex serves as a root of a subtree in 
$\Gamma_f$ with its leaves in $A$. As a result, 
any fiber component is cobordant to a union of fiber components indexed by elements of $A$. Therefore,
%%%%%
\begin{lem} 
Let $\pi_f: \Gamma_f \rightarrow S^1$ be such that no positive loop in $\Gamma_f$ contains 
vertices only of a particular index. Then the obvious homomorphisms 
$\Z[A] \rightarrow H_2^f$, \, $\Z_+[A] \rightarrow H_2^{f +}$ are onto. \qed
\end{lem}
%%%%%%
\begin{thm} A map $f: M \rightarrow S^1$ with no local extrema and no bubbling singularities is intrinsically harmonic, 
if and only if, 
$Ker\{\Z_+[A] \rightarrow H_2(M, \partial M; \Z)\} = 0$. That is, no positive combination of the $f$-oriented fiber 
components  from $A$ is homologous to zero in $M$ modulo $\partial M$.

In particular, if no non-trivial class of $H_2(M, \partial M; \Z)$ has a spherical or disk 
representative\footnote{Equivalently, when the Thurston semi-norm is a norm.} and $f$ is 
harmonic, then $Ker\{\Z_+[A] \rightarrow H_2(M, \partial M; \Z)\} = 0$. 
\end{thm}
%%%%%%
{\bf Proof.}\quad  
The Calabi's positive loop property is equivalent to the intrinsic harmonicity of $f$ ([C]). It implies that, 
for any $f$-oriented fiber component $F$, there is a positive loop $C$, so that
$F \circ C > 0$. Furthermore, for every other fiber component $F'$,\; $F' \circ C \geq 0$. Therefore, by Lemma 4.5,
$Ker\{\Z_+[A] \rightarrow H_2(M, \partial M; \Z)\} = 0$. \smallskip

On the other hand, if the Calabi positive loop property fails for a point $x \in M$, then the upper world 
$U_x$ of $x$ --- the set of points in $M$ which can be reached from $x$ following an $f$-positive path ---  
is bounded by several fiber components (one of which contains $x$ and the rest contain some singularities of 
$f$ ([C]). Along these components the gradient of $f$ is directed inwards $U_x$. Hence, the  
union of these $f$-oriented components produces a trivial element in $H_2(M, \partial M; \Z)$. Since, by Lemma 4.6,
any $f$-oriented component is cobordant to a union of a few components indexed by elements of $A$, we have 
produced a nontrivial element in the kernel $Ker\{\Z_+[A] \rightarrow H_2(M, \partial M; \Z)\}$.

To validate the last statement of the proposition, we notice that for harmonic maps spherical 
fiber components must be homologically non-trivial, which contradicts to the postulated nature of $M$. 
\qed

\begin{lem} For an intrinsically harmonic $f$, $H_1(\Gamma_f; \Z)$ admits a basis represented by 
$\pi_f$-positive loops. 

If no non-trivial class of $H_2(M, \partial M; \Z)$ has a spherical or disk 
representative, then the Calabi graph $\Gamma_f$ has no positive loops with vertices only 
of a particular index.
\end{lem}
{\bf Proof.} Since through any point $x \in \Gamma_f$ there exists a positive loop, the statement 
follows by the induction on the number of edges in the complement to a maximal tree in $\Gamma_f$. 
The second statement follows from Lemma 4.1. \qed 
\bigskip

Now we introduce two quantities which (in Section 8) will play  an important role in our arguments. The 
first one is the difference ${\chi_-}(F_R) - {\chi_-}(F_A) = \|F_R\| - \|F_A\|$.
\begin{lem} For any $\Gamma_f$ as in Lemma 4.2, 
\begin{eqnarray} var_{{\chi_-}}(f) = {\chi_-}(F_R) - {\chi_-}(F_A) 
\end{eqnarray}
\end{lem} 

{\bf Proof.}\quad We notice that 
${\chi_-}(F_R) - {\chi_-}(F_A) = \sum_{r\in R}\;[ {\chi_-}(F_r) - \sum_{a\in T_r^+}{\chi_-}(F_a)]$. 
Note that $[ {\chi_-}(F_r) - \sum_{a\in T_r^+}{\chi_-}(F_a)]$ is 
twice the number of non-bubbling vertices of index 2 on the tree $T_r^+$. Hence, 
$\chi_-(F_R) - \chi_-(F_A)$ equals the total number of non-bubbling $f$-singularities of index 2, while
$g(F_R) - g(F_A)$ equals the total number of "bivalent" (i.e. locally non-separating) $f$-singularities of index 2.
(cf. Lemmas 3.3---3.5.) 

Similar counting which employs the negative trees $\{T_r^-\}$ will reveal  $\chi_-(F_R) - \chi_-(F_A)$ 
as the total number of non-bubbling $f$-singularities of index 1, and 
$g(F_R) - g(F_A)$ as the total number of "bivalent"  $f$-singularities of index 1.  Again, using Lemmas 3.4, 3.5, 
we see that  formula-definition (3.1) gives still another count of 
non-bubbling  singularities.  \qed
\bigskip

Denote by $F_{best}^{[F_A]}$ a surface  $\coprod_{a \in A} \kappa_a\cdot F_a$---an oriented union 
of fiber components--- which represents the homology 
class $[F_A] := \sum_{a \in A} [F_a]$ and delivers the minimal value  of 
$\sum_{a \in A} |\kappa_a|\cdot {\chi_-}(F_a)$ among 
such representatives.\smallskip

We introduce $Var_{{\chi_-}}(f)$ --- a modification of $var_{{\chi_-}}(f)$ --- via the formula
\begin{eqnarray} Var_{{\chi_-}}(f) = {\chi_-}(F_R) - {\chi_-}(F_{best}^{[F_A]}) =
\|F_R\| - \|[F_A]\|_{H^f}
\end{eqnarray}
Evidently, $Var_{{\chi_-}}(f) \geq var_{{\chi_-}}(f)$. Although $var_{{\chi_-}}(f)$ 
is a more pleasing invariant (it connects in a more direct way with the $f$-singularities), 
actually, it is $Var_{{\chi_-}}(f)$ which will participate more often in our estimates. 

%%%% 
\begin{example}
{\bf The Harmonic Twister}
\end{example}
%%%%
Let us consider the case when the graph $\Gamma_f$ is an oriented loop with  
two vertices $a$ and $b$ of indices 1 and 2 respectively. The map 
$\pi_f : \Gamma_f \rightarrow S^1$ is the obvious one. The 1-chain $\tau_g(f)$ 
on $\Gamma_f$ is defined to take the value $n$ on the oriented arc $(b, a)$ and 
the value $n + 1$ on the oriented arc $(a, b)$.  It is easy to realize these 
data by a map $f: M^3 \rightarrow S^1$, $M^3$ being a closed manifold. By
applying 
a move $A$ from  Figure 4 we send vertex $b$ on a "round trip" and homotop $f$ 
to a new map $f_1$. The new map will generate a new      
chain $\tau_g(f_1)$, taking the value  $n + 1$ on $(b, a)$  and the value $n + 2$ on
$(a, b)$.
This deformation can be repeated again and again to produce maps with arbitrary
big 
genera $n + k$ of the best fiber. Of course, the original best fiber of genus $n$ 
is still residing in $M^3$; however, it is \emph{invisible} on the level of the
new graphs $\Gamma_{f_k}$ (in this case, identical with the original one) and new
chains $\tau_g(f_k)$.

A similar argument applies to the $\chi_-$-invariants and the chain 
$\tau_{\chi_-} (f)$.
 
As Corollary 8.9 and Example 8.12 imply, the original $f$-fiber $F_0$ will intersect the new
$f_k$-fiber along a complex pattern of loops, none of which bounds a disk in the 
new fiber(cf. Figure 15). Furthermore, these intersections cannot be 
removed even by an isotopy of $F_0$. 

Since all the graphs $\Gamma_{f_k}$ satisfy the Calabi positive loop
property, 
all the maps $f_k$ are intrinsically harmonic [FKL], and all the
$f_k$-fibers 
are near-minimal surfaces [K]. This means that, for any choice of two
disjoint 
3-disks $D_1$ and $D_2$, centered on the two singularities, and any 
$\epsilon > 0$,
there exists a riemannian metric on $M^3$ with the following properties: 1) the 
map $f_k$ is harmonic; 2) the $f_k$-fibers are minimal surfaces outside of the
disks; 
3) the area of the portion of each fiber inside the disks is smaller than
$\epsilon$ 
(in other words, by the choice of metric, the deviation of fibers from the
minimality 
can be localized around the singularities and made numerically insignificant). 

In contrast with the fibrations over a circle, as the Twister example 
demonstrates, this "near-tautness" of the singular foliation $\mathcal F_{f_k}$,
does not imply the minimality of the $\chi_-$-characteristic of the best fiber 
in its homology class. 
\qed
\smallskip 

A challenging problem is how to "untwist" a given map $f: M \rightarrow S^1$ and 
to lower the $l_1$-norms of the characteristic chains $\tau_{\chi_-}(f)$, 
$\partial\tau_{\chi_-}(f)$, or the the value $\chi_-(F_{best})$.  
The twist invariant $\rho_{\chi_-}([\Sigma], F_R)$ from Section 6 does measure the 
"twist" of $f$ (relative to $[\Sigma]$). Regrettably, I do not know how to produce 
maps with the zero twist out of a variational principle.

Deforming a generic map $f: M \rightarrow S^1$ to a harmonic map, 
while preserving the list of its singularities,  as it is done in the proof 
of Theorem 1, pages 474-475, in [FKL], requires elementary moves $D$ and $E$ from 
Figure 5. Unfortunately, they have 
the potential to increase $\chi_-(F_{best})$ and the variation. At the same time,
some form of harmonicity seems to be an essential ingredient in our arguments, especially 
if one expects a \emph{fiber}  to deliver the Thurston norm
(cf. Figure 1 depicting a non-harmonic map whose fibers fail to deliver the Thurston norm). 
This tension between our desire to lower the value $\chi_-(F_{best})$ and the need to 
harmonize the map $f$ calls for an investigation beyond the scope of this paper. 

%%%%%%%%%%%%%%%%%%%%%%%%%%%%%%%%%%%%%%%
 
\section{Self-indexing maps to $S^1$ and the $\chi_-$-minimizing 2-cycles}   

To avoid combinatorial complications and to make the future arguments more transparent, 
first we treat the case of  \emph{self-indexing}  maps 
$f: M \rightarrow S^1$ with no local extrema. Automatically, such maps  are intrinsically harmonic.  
For  a self-indexing $f$, unless it is a fibration,  $\Gamma_f$ is a union of two trees which 
share the same root and the same set of leaves.  More general maps are considered in Section 8. 
\smallskip
\smallskip 

Let $F$ and $\Sigma$ be  oriented embedded surfaces which intersect
transversally in $M$. When 
$\partial M = \phi$, we assume the surfaces are closed; otherwise, their
boundaries are
contained in $\partial M$. The intersection $\mathcal C = F \cap \Sigma$
consists of a 
number of oriented loops and arcs. Their orientations are induced by the
orientations of 
$F$, $\Sigma$ and $M$. As we modify the intersection, we still call it $\mathcal
C$.

\begin{defn} The oriented intersection pattern $\mathcal C = F \cap \Sigma$ 
is \emph{totally reducible}, if it is comprised of curves which bound disks in $F$ 
or of arcs  which bound relative disks in $(F, \partial F)$.
\end{defn}   

\begin{thm} 
Let $f: M \rightarrow S^1$ be a map from an oriented  3-manifold $M$ to an oriented circle. 
Assume that: 
\begin{itemize}
\item $f$ has no critical points of indices 0 and 3;
\item if $\partial M \neq \phi$, then $f: \partial M \rightarrow S^1$ is a
fibration;
\item  along the circle, the critical values of critical points of index 1 belong to an 
oriented arc $(\theta_b , \theta_w )$, and the critical values of points of
index 2 --- to a complementary oriented arc $(\theta_w , \theta_b)$.\footnote
{In other words, $f$ is a \emph{self-indexing} map.} 
\end{itemize}

Let $(\Sigma,  \partial \Sigma) \subset (M, \partial M)$ be an embedded oriented
surface, homologous $rel. \;\partial M$ to a fiber.  Assume that $\Sigma$
has a \emph{totally reducible} intersection $\mathcal C$ with the fiber  
$F_{worst} = f^{-1}(\theta_w)$. 

Then $\chi_- (\Sigma) \geq \chi_- (F_{best})$ and $g(\Sigma) \geq g(F_{best})$, 
where $F_{best} = f^{-1}(\theta_b)$. 
\end{thm}

{\bf Proof.}\quad Let $(\Sigma, \partial\Sigma) \subset (M, \partial M)$ be an
oriented regularly embedded surface, 
homologous to a fiber modulo $\partial M$, and such that $\mathcal C =
\Sigma \cap F_{worst}$ 
is totally reducible. If a loop $\gamma$ from $\mathcal C$ bounds a disk $D \subset 
F_{worst}$ which 
is free of any other intersection curves, then we can perform a 2-surgery on
$\Sigma$ along 
$\gamma$ with the $D$ as the core of a 2-handle. The resulting surface $\Sigma'$ is
homologous to 
$\Sigma$ and has $\mathcal C \setminus \gamma$ as its intersection set with
$F_{worst}$.  
According to Surgery List A, $g(\Sigma ') \leq  g(\Sigma)$ and 
$\chi_- (\Sigma') \leq \chi_- (\Sigma)$.
A similar argument, based on Surgery List B, leads to a similar conclusion for 
a surface resulting 
from the surgery along an arc $\gamma$,  bounding a disk modulo $\partial
F_{worst}$. When $\gamma$ bounds in $F$ a disk which contains 
other curves from $\mathcal C$, we perform 2-surgery on $\Sigma$ starting with 
the "most interior" disks and gradually "moving outwards".

Because $\Sigma \cap F_{worst}$ is totally reducible, we can produce a new
surface $\Sigma'$ 
in the homology class of the original $\Sigma$ having an \emph{empty} 
intersection with $F_{worst}$. 
Moreover, $g(\Sigma ') \leq  g(\Sigma)$ and $\chi_- (\Sigma') \leq \chi_-
(\Sigma)$. So, if we will be able to prove the desired inequalities for $\Sigma'$,
then they will be valid for the original $\Sigma$ as well.
\smallskip

Now we revert to a generic notation $\Sigma$ for this new surface $\Sigma'$.
Starting at $\theta_{b}$ and moving along the oriented circle, first we cross 
the critical values of \emph{all} index-one critical points. Similarly, 
starting at 
$\theta_{w}$ and moving along the oriented circle,
we first  meet the critical values of \emph{all} index-two critical points. 

We cut the manifold $M$ open along $F_{best} = f^{-1}(\theta_{b})$. The boundary
of the 
resulting manifold $\hat{M}$ consists of two copies of $F_{best}$---the surfaces
$F_{best}^0$ and 
$F_{best}^1$. 

Denote by $\tilde{M}$
the cyclic covering space of $M$, induced by the map $f$.
We can regard the $\hat{M}$ as a fundamental domain of the cyclic action on
$\tilde{M}$. 
The map $f$ generates a Morse function $\hat{f}: \hat{M} \rightarrow [0, 1]$, 
so that $\hat{f}^{-1}(0) = F_{best}^0$ and $\hat{f}^{-1}(1) = F_{best}^1$. 
The original $f$ induces a $\Z$-equivariant Morse function $\tilde{f}: \tilde{M}
\rightarrow \R$ with no singularities along $\partial\tilde{M}$.
Here $\partial\tilde{M}$ denotes the preimage of $\partial M$, under the covering map 
$\tilde{M} \rightarrow M$. We can view $\hat f$
as a restriction of $\tilde f$ to the fundamental domain.
\smallskip

By the construction, all $\hat{f}$-critical points of index 1 lie below the
surface    
$F_{worst} = f^{-1}(\theta_{w})$, and all critical points of index 2 ---above
it. \smallskip 

We consider an $f$-gradient-like vector field $X$ on $M$ and its lifting $\tilde{X}$
on the 
covering space $\tilde{M}$. When necessary, these fields will be adjusted.
Since $f: \partial M \rightarrow S^1$ is a fibration, we always can assume that $X$ 
is \emph{tangent} to the boundary $\partial M$ and does not vanish there.

Let index $\alpha$ enumerate the critical points $\{x_\alpha \}$ of index 1, and 
index $\beta$ --- the critical points $\{x_\beta \}$ of index 2 in $M$. Denote by 
 $\{\hat x_\alpha \}$ and $\{\hat x_\beta \}$ the corresponding $\hat{f}$-critical
points in $\hat{M} \subset \tilde{M}$. 

We denote by $\tilde{D}_{\alpha}^1$ the two descending trajectories of the field $\tilde{X}$, 
which emanate from the singularity $\{\hat x_\alpha \}$, and by
$\tilde{D}_{\alpha}^2$ --- 
the union of all ascending trajectories.   
Similarly, we denote by $\tilde{D}_{\beta}^1$ the two ascending trajectories and
by $\tilde{D}_{\beta}^2$ --- the union of all descending trajectories, which
emanate from the singularity $\{\hat x_\beta \}$. 

Let $\hat{D}_{\alpha}^1 = \tilde D_\alpha^1\cap \hat M$,\; 
$\hat{D}_{\alpha}^2 = \tilde D_\alpha^2\cap \hat M$,\; $\hat{D}_{\beta}^1 = \tilde D_\beta^1\cap \hat M$,\;
$\hat{D}_{\beta}^2 = \tilde D_\beta^2\cap \hat M$.

Since $\tilde{X}$ is tangent to the boundary $\partial\tilde{M}$ and does not vanish there, 
all the sets $\{\tilde{D}_{\alpha}^i\}$ and  $\{\tilde{D}_{\beta}^i\}$ 
($i = 1, 2$) have an \emph{empty intersection} with $\partial\tilde{M}$.

Note that the portions of the sets $\hat{D}_{\beta}^2$ lying \emph{above} the surface
$F_{worst}$, and of the sets $\hat{D}_{\alpha}^2$ lying \emph{below} the surface $F_{worst}$
are diffeomorphic to two-dimensional disks.  
\smallskip

Denote by $\hat{\Sigma}$ the preimage of
$\Sigma$  under the natural map $p: \hat{M} \rightarrow M$. Clearly,
$\hat{\Sigma }$ has an empty 
intersection with the surface $F_{worst} \subset \hat{M}$ and, therefore,
$\hat{\Sigma }$ is 
divided into two disjoint pieces: $\hat{\Sigma }^0$ lying below $F_{worst}$ and  
$\hat{\Sigma }^1$ lying above it. We aim to push $\hat{\Sigma }^0$ towards
$F_{best}^0$ and 
\emph{below} any critical value, produced by the singularities of index 1. 
At the same time, we will try to push $\hat{\Sigma }^1$ towards $F_{best}^1$ and 
\emph{above} any critical value, produced by the singularities of index 2. 
In general, both 
desired isotopies are obstructed by the critical points $\{\hat x_\alpha \}$ and  
$\{\hat x_\beta \}$; however, after 2-surgery on $\hat{\Sigma }$, as Figure 6
indicates, the two deformations will become possible. (Note that,  Figure 9
depicts the simplest case, when the intersection of $\hat{D}_{\beta}^2$ with 
$\hat{\Sigma}$ consists of a single loop.)  

%%%%%%
\begin{figure}[ht]
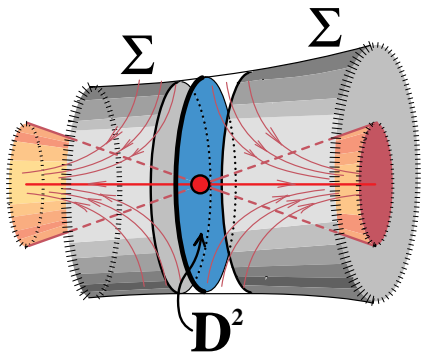

\centerline{\BoxedEPSF{2-surgery scaled 450}}
\bigskip
\caption{Moving  $\Sigma$ above the critical level of
index 2 after a 2-surgery,}
\end{figure} 
%%%%%%%%
Consider the intersection of each set $\tilde{D}_{\beta}^2$ with the surface
$\hat{\Sigma}^1$. 
Above the level of $F_{worst}$, the set $\tilde{D}_{\beta}^2$ is a smooth disk.
Since $\hat{\Sigma}^1$ lies above 
$F_{worst}$, by a choice of the gradient vector field $X$ (equivalently, by a
small 
isotopy of $\Sigma$) the intersection of $\hat{\Sigma}^1$ with the disk can be
made transversal.
 It will consist of a number of closed simple curves $S^1_{\beta , k}$, shared
by the disk and 
the surface. The curves are closed because the disks $\tilde{D}^2_\beta$ can not
reach the 
boundary $\partial\tilde{M}$---the gradient field has been chosen to be tangent
to the 
boundary and has no zeros there. We will use these loops to perform
2-surgery on
$\hat{\Sigma}^1$ inside 
$\hat{M}$.  We start with the most "inner" loop in the disk, say, with
$S^1_{\beta , 1}$. Inside of the set  
$\tilde{D}_{\beta}^2$ it bounds a disk $D^2_{\beta , 1}$, which is embedded in
the ambient $\hat{M}$.
The normal frame to $S^1_{\beta , 1}$ in $\hat{\Sigma}^1$ extends to a normal
frame of 
$D^2_{\beta , 1}$ in $\hat{M}$. 
We perform a surgery on $\hat{\Sigma}^1$ along $S^1_{\beta , 1}$ by attaching a
2-handle 
$D^2_{\beta , 1}\times [0, 1]$. The surgery eliminates the intersection
$S^1_{\beta , 1}$ 
and does not produce new intersections with the "disk" $\tilde{D}_{\beta}^2$.
This procedure 
is repeated, until all the intersection loops $\{S^1_{\beta , k}\}_{\beta, k}$
are eliminated. 
The resulting surface $\hat{\Sigma}^1_\star$ has an empty intersection with all
the sets 
$\{\tilde{D}_{\beta}^2 \}$, and therefore, can be isotoped along the gradient
flow to 
a location \emph{above} all critical points $\{\hat x_\beta \}$. In the process, the
boundary 
$\partial\hat{\Sigma}^1_\star$ glides along the boundary $\partial\hat{M}$. 
The isotopy can be chosen to be fixed in a small neighborhood of $F^1_{best}$.

A similar treatment can be applied to the surface $\hat{\Sigma}^0$ and the
portions of 
$\{\tilde{D}_{\alpha}^2 \}$ lying below the surface $F_{worst}$. It will produce
a new 
surface $\hat{\Sigma}^0_\star$ located below all the critical points $\{\hat x_\alpha
\}$,  
in a neighborhood of $F^0_{best}$. Because the surgery did not touch a portion
of 
the original $\hat{\Sigma}^0$ in a neighborhood of $F^0_{best}$ and 
a portion of the original $\hat{\Sigma}^1$ in a neighborhood of $F^1_{best}$, 
the new $\hat{\Sigma}^0_\star$ and $\hat{\Sigma}^1_\star$ still define a surface 
$\Sigma_\star$ in $M$, located in a regular tubular neighborhood of the fiber
$F_{best}$.
\smallskip

Since $\Sigma$ and $\Sigma_\star$ are linked by 2-surgery inside $M$, followed
by a 
regular isotopy, they define the same class in $H_2 (M;\partial M; \Z)$.

On the other hand, by Lemma 3.1, their genera  and $\chi_-$-characteristics 
satisfy the inequalities: $g(\Sigma) \geq g(\Sigma_\star)$ and  
$\chi_- (\Sigma) \geq \chi_- (\Sigma_\star)$.
\smallskip

Applying Lemma 5.3 below, with $\Sigma = \Sigma_\star$, $F = F_{best}$ and $d = 1$,
we complete the proof of Theorem 5.1.\qed

\begin{lem} Let $U$ be a regular neighborhood of a connected oriented surface $F$ and 
$V$ a regular neighborhood of its boundary $\partial F \subset \partial U$. 
Let  $(\Sigma, \partial\Sigma) \subset (U, V)$ be an embedded oriented surface. Then  
$\chi_- (\Sigma ) \geq  |d| \cdot \chi_- (F)$ and $g(\Sigma ) \geq |d| \cdot g(F)$,
where $d$ is a total degree of the map $\Sigma \rightarrow F$ induced by the retraction 
$p: (U, V) \rightarrow (F, \partial F)$.   
\end{lem}

{\bf Remark.}\quad   
The assumption that $\Sigma$ is an embedded surface is important: it is 
easy to construct examples of immersed surfaces which violate the conclusion 
of the lemma. 
Take for instance, a double cover of a torus $T$ by another torus $T_1$,
immersed in 
$T \times [0, 1]$ with the projection $T \times [0, 1] \rightarrow T$ inducing
the covering degree-two map $T_1 \rightarrow T$. Here $g(T_1) = g(T)$ and not 
twice $g(T)$. 
\smallskip

{\bf Proof.} Let $\Sigma = \coprod_j \Sigma_j$, where each $\Sigma_j$ is connected and 
let $p_j: \Sigma_j \rightarrow F$ be a map of degree $d_j$ induced by the retraction 
$U \rightarrow F$.  Recall that if $d_j \neq 0$, then  
the homomorphism $(p_j )_\ast : H_\ast (\Sigma_j ; \R ) \rightarrow H_\ast (F_0 ; \R )$ 
is an epimorphism [W]. In particular,  
$(p_j)_\ast : H_1 (\Sigma_j ; \R ) \rightarrow H_1 (F_0 ; \R )$ is an epimorphism. 
Hence, for such a component, $g(\Sigma_j) \geq g(F)$ and 
$|\chi(\Sigma_j)| \geq |\chi(F)|$.
\smallskip

Put $d = \sum_j\, d_j$.
We claim that $\Sigma$ must divide $U \approx F\times[0, 1]$ into at least  $|d| + 1$ 
connected regions $\{U_i\}$. If the complement to $\Sigma$ in $U$  consists 
of less than $|d| + 1$ components, then it is possible to construct a path 
$\gamma \subset U$, which connects a point $a \in F \times \{0\}$ with a point $b  \in F \times \{1\}$
 and which has less than $|d|$ intersections with $\Sigma$. 
Existence of such a $\gamma$ contradicts with the fact that $[\Sigma] = \sum_j\;[\Sigma_j]$ and $d[F]$ 
are homologous in $(U, V)$. 

Take the  region $U_1$ adjacent to $F \times \{0\}$ and let $V_1 :=  U_1 \cap V$. 
The projection $p: U_1 \rightarrow F \times \{0\}$ maps at least 
one component of the surface $\partial U_1 \setminus V_1$, distinct from $F\times \{0\}$,
by  a  degree 1 map. Indeed, consider the intersections of a generic segment 
$I = x \times [0, 1], \; x\in F$, with  $\partial U_1 \setminus V_1$. Among them pick the highest intersection,
say $a$. Let $S_1$ be the component of $\Sigma$ which contains $a$.
At $a$ the path $I$ leaves the domain $U_1$ "forever". Take a point $b \in I$ just below $a$ and 
connect $b$ by a path  $\gamma \subset U_1$ with the base of $I$. One can construct $\gamma$ in such 
a way that its $p$-projection is a loop homologous to zero in $F$: just add an appropriate 
kick in $F$ to any candidate for $\gamma$. Denote by $J$ the portion of
$I$  above $b$. The new path $K$ ---the union of $J$ with $\gamma$---intersects with $S_1$ at a single point 
$a$ and shares its beginning and end with the $I$. Furthermore, the loop $K \cup I$ is 
null-homologous in $F\times[0,1]$. Since the algebraic intersection number of $K$ with 
with $S_1$ is 1, so must be the intersection number with $I$, that is, $deg(p|_{S_1}) = 1$.
 
Next, consider a region $U_2$ adjacent to $U_1$ along $S_1$. 
The same reasoning now applies to $U_1 \cup U_2$.  We conclude
that it must have a boundary component $S_2 \neq S_1$, 
which projects into $F \times \{0\}$ by the map $p$ of degree 1.    
This inductive process provides us with at least $|d|$ components $\{\Sigma_j = S_j\}$ of $\Sigma$, 
each of which  maps onto $F \times \{0\}$ by a  degree 1 map.
Thus, $g(\Sigma) \geq |d| \cdot g(F)$. By the same token, each of the $|d|$ components
$\Sigma_j$ has the property $|\chi(\Sigma_j)| \geq |\chi(F)|$. Furthermore,
$\chi_-(\Sigma_j) \geq \chi_-(F)$---a sphere can not be mapped by a 
non-zero degree map onto an orientable surface, different from a sphere. Hence, 
$\chi_-(\Sigma) \geq |d|\cdot\chi_-(F)$.  \qed 
\smallskip

The following statement is very much in line with the Harmonic Twister example: 
it shows that harmonicity \emph{alone} is too weak and too flexible 
to insure the "best fiber theorem". For instance, the hypothesis in Theorem 5.2, 
requiring the probing surface $\Sigma$ to have a totally reducible intersection 
with the worst fiber, is essential.
%%%%% 
\begin{lem}
Let $M$ be an oriented closed 3-manifold.
Any \emph{connected}, orientable, embedded and non-separating
surface $\Sigma \subset M$  can be viewed as the best fiber of a self-indexing map 
$f: M \rightarrow S^1$ with all the fibers being connected.  Such an $f$ is intrinsically harmonic. 
\end{lem}
%%%%%
{\bf Proof.}\quad   
%In fact, if $\partial \Sigma$ does not separate each boundary torus it meets,
%such an $f$ can be assumed 
%to produce fibrations of the tori over the circle.  
Apply the Thom-Pontryagin construction to a regular neighborhood $U$ of
$\Sigma \subset M$.
Perturb $f$ away from $U$ to covert it into a Morse map. Then cut $M$ open along
$\Sigma$ to 
get a Morse function $\hat f: \hat M \rightarrow [0, 1]$ which
maps one copy of $\Sigma$, $\Sigma_0 \subset \partial \hat M$, to $0$, and the
other copy, 
$\Sigma_1 \subset \partial \hat M$, --- to $1$. Through a standard deformation of
$\hat f$, 
fixed on $\partial \hat M$, one can eliminate all local maxima and minima. 
This leaves us only with critical points of indices 1 and 2. 
Then we can deform $\hat f$ into a \emph{self-indexing} Morse function $\hat f'$
without changing it at $\partial \hat M$ ([M]). Since, $\Sigma$ is connected, 
the self-indexing $\hat f'$ must have only connected fibers (cf. Figure 2), and 
therefore, is intrinsically harmonic. By Lemma 3.1, $\Sigma$ is the best fiber of $f'$. \qed

%%%%%%%%%%%%%%%%%%%%%%%%%%%%%%%%%%%%%%%

\section{The twist of a map $f: M \rightarrow S^1$ and the 
$\tilde f$-breadth of surfaces}  

In this section we introduce a number of invariants characterizing the intersection complexity 
of two (hyper)surfaces $\Sigma, F \subset M$ which are specially positioned 
with respect to a given map $f: M \rightarrow S^1$. In Section 8 these invariants will contribute 
to our estimates of the Thurston norm. We also introduce the "twist" of
$f$ in terms of the intersection complexity of a surface $\Sigma$,
which  delivers the Thurston norm, with a generic fiber component. Ultimately, it is the presence 
of the $f$-singularities which is 
responsible for  the non-triviality of these invariants.
\smallskip

Let $(\Sigma, \partial \Sigma) \subset (M, \partial M)$  be an oriented surface 
representing an $f$-vertical class $[\Sigma]$ (cf. Definition 4.4). Let 
$F$ be a finite union of fiber components.
For  such a pair $(\Sigma, F)$, we introduce a non-negative integer $\rho(\Sigma, F)$. 
It will measure the complexity of the transversal intersection 
$\mathcal C = \Sigma \cap F$ inside $F$. The fact that $M$ is 3-dimensional is 
not important here: a similar invariants make sense for any map $f: M \rightarrow S^1$ and 
any pair of vertical hypersurfaces in $M$.  

The pattern $\mathcal C = \cup C_i \subset F$ is comprised of oriented simple 
curves (arcs and loops). Each curve $C_i$ is equipped with a  normal framing, 
induced by the preferred normal framing of $\Sigma$.\smallskip 

It is crucial to notice that the algebraic intersection number of any loop 
$\gamma \subset F$ with $\mathcal C$ is zero. Indeed, the algebraic intersection 
$\gamma \circ \mathcal C = \gamma \circ \Sigma$ of
such a $\gamma$ with any surface $\Sigma$ representing a vertical class $[\Sigma]\in H_2^f$ 
is zero: just consider $\Sigma'$ homologous to $\Sigma$ and comprised of fiber components 
distinct from those of $F$ to conclude that $\gamma \cap \Sigma' = \emptyset$.  
\smallskip

We consider an oriented graph $K_{\mathcal C}$ whose vertices correspond to connected
components of $F \setminus \mathcal C$ 
and edges --- to the connected components of $\mathcal C$.  The
orientation of the edges is prescribed by 
the preferred normal frames to the intersection curves. Because any loop in $F$
has a trivial algebraic intersection with
$\mathcal C$, each loop in $K_{\mathcal C}$ will have an equal number of
"clockwise" and "counter-clockwise" oriented edges.
   
Consider a 1-cochain $c$ on $K_{\mathcal C}$, which takes value 1 at each
oriented edge of the graph. Since, 
by the argument above, $c$ takes the zero value on every loop in $K_{\mathcal C}$, 
it is a coboundary: $c = \delta u$ \, for some 0-cochain 
$u$. The potential $u$ is a function on the vertices $C_0(K_{\mathcal C})$ of 
$K_{\mathcal C}$, which prescribes the flow $c$ through 
the edges. For each connected component of $F$, $u: C_0(K_{\mathcal C})
\rightarrow \Z$ is well-defined, up to 
a choice of a constant. We can synchronize  these choices  by equating
all the maximum values of $u$ on 
different connected components of $K_{\mathcal C}$. 

Denote by $\rho(\Sigma, F)$, or by $\rho (\mathcal C)$ for short, 
\emph{one less} the number of \emph{distinct values} 
taken by the synchronized function $u$. This 
integer will be our measure of complexity of the intersection $\mathcal C$. 

For an oriented graph with all vertices 
being  sources and sinks, $u$ takes only two values and $\rho (\mathcal C) = 1$.
In general, $\rho (\mathcal C)$ does 
not exceed the number of connected components in $F \setminus (F \cap \Sigma)$ 
minus one. 

One can think of the potential $u$ as an integral-valued function on 
$F \setminus \mathcal C$, constant on its components. In this interpretation, 
the curves from $\mathcal C$ can be imagined as dams erected on $F$, and the 
$u$-values---as the water levels for each of the fields from 
$F \setminus \mathcal C$. Depending on the orientations, crossing a 
dam results in a change of the water level by $\pm 1$. In this model,
$\rho (\mathcal C)$ is the integral variation of the water level across 
the irrigation system $\mathcal C$. In other words, for each connected 
component $F_\alpha$ of $F$, we consider an integral 2-chain $E_\alpha$ whose boundary 
is the 1-cycle $\Sigma \cap F_\alpha$, and define  $\rho (\mathcal C)$ as 
$max_\alpha \{osc (E_\alpha)\}$, where $osc (E_\alpha)$ denotes the 
oscillation in the values of the coefficients in the chain $E_\alpha$.

Assume that two curves $C_1$ and $C_2$ from 
$\mathcal C$ can be linked by an oriented arc $\gamma$ which has a single transversal 
intersection $x$ with $C_1$, a single transversal intersection $y$  with $C_2$, the two 
intersections being of \emph{opposite} signs. Also, assume that $\gamma$ misses the rest 
of the curves from $\mathcal C$. Then the water level $u_\star$ along $\gamma$ before it 
hits $C_1$ and after it hits $C_2$ must be equal. So, we can connect the corresponding 
fields by a canal following $\gamma$ and fill it with water up to the level $u_\star$. 
This irrigation construction will merge the two fields into a single one, replaces 
$C_1$ and $C_2$ with their connected sum $C_1 \# C_2$, but will not change the value of 
$\rho(\sim)$.
\smallskip 

As we isotop the surface $\Sigma$ in $M$, its transversal intersections $\mathcal C$ 
with $F$ are subjected to an isotopy in $F$ and occasional surgery of the  
types $C_1 \sqcup C_2 \Rightarrow  C_1 \# C_2$, $C_1 \# C_2 \Rightarrow  C_1 \sqcup C_2$, 
or of the birth-annihilation types $C \Rightarrow \emptyset$,\, $\emptyset \Rightarrow C$.  
Here the loop $C$ bounds a disk in $F$ (or in $(F, \partial F)$) \emph{and} in $\Sigma$. 
 Also, a different type of 
surgery can occur: it corresponds to  connecting two points $x$ and $y$ on the 
\emph{same} curve $C$ by an oriented
arc.  It has an effect of separating $C$ into two components $C_1$ and
$C_2$.   Under the transformations $C_1 \sqcup C_2 \Rightarrow  C$,
$C \Rightarrow  C_1 \sqcup C_2$, the value of $\rho(\sim)$ is preserved. 
Only the birth-annihilation surgery can change it.
\smallskip 

A modified  definition of $\rho(\Sigma, F)$ will be useful. In the 
modification, from the very beginning, we  \emph{exclude} all loops from 
the intersection $\Sigma \cap F$ which \emph{bound a disk} in $F$. 
We also exclude arcs of $\Sigma \cap F$ which bound a disk in $(F, \partial F)$. 
This gives us a simpler intersection pattern $\mathcal C^\circ$. 
Then we employ $\mathcal C^\circ$ to 
define $\rho^\circ(\Sigma, F)$ as $\rho(\mathcal C^\circ)$. This 
quantity can also change under the surgery of the type 
$C_1 \sqcup C_2 \Rightarrow  C_1 \# C_2$.\smallskip

The definition below  introduces  new invariants which depend only on the  
homology class $[\Sigma] \in H_2^f$, a value ${\chi_-}(\Sigma)$, and a  surface $F \subset M$
which is a union of fiber components (alternatively, whose fundamental class $[F]$ is proportional
to $[\Sigma]$).
%%%%%
\begin{defn} Let 
$\rho_{{\chi_-}}(\Sigma, F)$ be the \emph{minimum} of 
$\{\rho(\Sigma', F)\}_{\Sigma'}$, where $\Sigma' \subset M$ is  
\emph{homologous} to $\Sigma$ and ${\chi_-}(\Sigma') \leq {\chi_-}(\Sigma)$.
\end{defn}
%%%%%%%
\begin{lem}
$\rho_{{\chi_-}}(\Sigma, F)$ is the minimum of 
$\{\rho^\circ(\Sigma', F)\}_{\Sigma'}$, where $\Sigma' \subset M$ is  
\emph{homologous} to $\Sigma$ and ${\chi_-}(\Sigma') \leq {\chi_-}(\Sigma)$.
\end{lem}
%%%%%
{\bf Proof.}\quad By performing 2-surgery on any given $\Sigma$ along 
curves from $\Sigma \cap F$ which bound disks in $F$, we 
can replace $\Sigma$ with a new surface $\Sigma'$, such that $[\Sigma'] = [\Sigma]$,
${\chi_-}(\Sigma') \leq {\chi_-}(\Sigma)$, and $\rho^\circ(\Sigma', F) = \rho(\Sigma', F)$. \qed
\smallskip

\smallskip

\begin{defn} We fix a vertical homology class $[\Sigma] \in H_2^f$ and a repeller set 
$R \subset \Gamma_f$. Consider all oriented surfaces $\{\Sigma \subset M\}$ which deliver the 
minimal value of ${\chi_-}(\sim)$ in the homology class $[\Sigma]$. Among them pick 
$\Sigma$'s with the minimal value of the twist $\rho(\Sigma, F_R)$\footnote{Due to Lemma 6.3, 
this is equivalent to minimizing $\rho^\circ(\Sigma, F_R)$.}. We denote this optimal 
value by $\rho_{{\chi_-}}([\Sigma], F_R)$ and call it the $R$-\emph{twist} of 
$f$ relative to the class $[\Sigma]$.
\end{defn}
Since any two sets of repelling components are isotopic, $\rho_{{\chi_-}}([\Sigma], F_R)$ 
does not depend on a particular choice of $R \subset \Gamma_f$.  

For technical reasons, Definition 6.3  employs a special union
$F_R$ of  fiber components.  Replacing the $F_R$ in Definition 6.4
with "any fiber component $F$", one can introduce a modified 
definition which makes sense for any Morse map $f$. 

When $[\Sigma] = [F]$ is the homology class of a fiber, we also will use the abbreviation 
$\rho_{{\chi_-}}(f)$ for $\rho_{{\chi_-}}([F], F_R)$.\smallskip

We shall see that the $S^1$-controlled size of a homotopy, which links a given map $f$ to 
the one whose best fiber delivers the Thurston norm, gives an upper bound for $\rho_{\chi_-}(f)$.
\smallskip
%%% 
\begin{example}
(Jerome Levine)
\end{example}
%%%%
Let $\tilde M \rightarrow M$ be a cyclic cover induced by $f: M  \rightarrow S^1$. 
Not for every $f$-vertical homology class $[\Sigma]$, each surface
$\Sigma$ which  realizes  $[\Sigma]$ admits a lift to $\tilde M$. Take, for instance, a connected  sum
$M = (F_0 \times S^1) \# (F_1 \times S^1)$,
where $F_0$ and $F_1$ are oriented surfaces of your choice. Let $f: M \rightarrow S^1$  be a 
connected sum of obvious projections.  Consider $\Sigma = F_0 \# (\partial D^2 \times S^1)$, where 
$D^2 \subset F_1$ is a 2-disk. Note that $[\Sigma] = [F_0]$---clearly, a vertical class---, but $\Sigma$ 
can not be lifted to $\tilde M$: it contains a loop $\gamma$ which is mapped by the degree 1 map $f|$ onto $S^1$. 
We notice that $\gamma \circ \Sigma = 0$, which is  in agreement with the lemma above. The immersed 
surface $\Sigma \cup F_1$, realizing 
the homology class $[F_0] + [F_1]$ (which satisfies the lemma's hypotheses), demonstrates 
that Lemma 6.6 below can not be generalized for immersed surfaces.\smallskip

\begin{lem} 
If $[\Sigma]$ is proportional to the homology class $[F]$ of an $f$-fiber, then
\emph{any} surface $\Sigma$ realizing $[\Sigma]$ admits a lifting to the space $\tilde M$.
When $\Sigma$ is connected, the lift $\hat \alpha: \Sigma \subset \tilde M$ is unique up 
to the deck transformations in $\tilde M$. 
\end{lem}

{\bf Proof.}\quad By Poincar\'{e} duality, any 2-homology class is characterized by its 
intersections with loops forming a basis in $H_1(M; \Z)$. In particular, the homology 
class $[F]$ of a fiber is characterized by the property $[F]\circ[\gamma] = f_\ast([\gamma])$,
where $f_\ast([\gamma])$ stands both for the image of $[\gamma]$ --- an integer --- under 
$f_\ast: H_1(M; \Z) \rightarrow H_1(S^1; \Z) \approx \Z$. Thus, any class $k[F]$ is 
determined by the property $k[F]\circ[\gamma] = k\cdot f_\ast([\gamma])$. By linearity,
this is equivalent to the proposition:  $[\Sigma]$ is proportional to $[F]$, iff,  
$[\Sigma]\circ[\gamma] = 0$ for any $[\gamma] \in Ker(f_\ast)$. Evidently, any class $[\Sigma]$ which
satisfies this  criteria has the property: 
$\{f_\ast([\gamma]) \neq 0\} \Rightarrow \{[\Sigma]\circ[\gamma] \neq 0\}.$ 
 
Now, in order to prove the lemma, it suffices to show that the image of any loop 
$\gamma \subset \Sigma$ under $f$ 
is null-homotopic in $S^1$. Since $\Sigma$ is two-sided, $[\gamma] \circ [\Sigma] = 0.$ On the other 
hand, if the loop $f([\gamma]) \neq 0$, then by the argument above,  
$[\gamma] \circ [\Sigma] \neq 0$. This contradiction completes the proof.   \qed
\bigskip

Put $\hat \Sigma := \hat\alpha(\Sigma)$. If $\Sigma$ consists of many components $\{\Sigma_j\}$, 
each of them admits its own lift $\hat \alpha_j: \Sigma_j \subset \tilde M$. However, 
not any combination $\hat \Sigma$ of $\{\hat \Sigma_j \subset \tilde M\}_j$ will serve our goals.
We will be especially interested in liftings $\hat\alpha$ such that the (relative) 2-cycle 
$\hat \Sigma$ is the \emph{boundary} 
of an integral 3-chain $C$ in $\tilde M$ modulo  
$\tilde{\partial M} \cup \tilde M^{+\infty} \cup \tilde M^{-\infty}$. Here 
$\tilde M^{+\infty}$ ($\tilde M^{-\infty}$) stands for the positive (negative) ends 
of $\tilde M$.\footnote{When $f$ is primitive, $\tilde M$ is connected and has a single positive 
and a single negative end.} We denote the set of such special liftings by $\mathcal
B(\Sigma)$.

The surface $\hat \Sigma$ divides $\tilde M$ into a finite number connected domains $\{U_l\}_l$.
When $\hat \Sigma \in \mathcal B(\Sigma)$, the (infinite) 3-chain $C$ can be chosen so that it  
attaches the same integral multiplicity $u_l$ to every 3-simplex from the domain $U_l$. 
This  type of observation is already familiar from the beginning of Section 6, where it was employed 
(in dimension 2) to introduce the invariant $\rho(\Sigma, F)$. Indeed, if $\hat \Sigma =
\partial C$, then $\hat \Sigma \circ \gamma = 0$ for any loop $\gamma \subset \tilde M$. This allows us to define 
$u_j$ as $\hat \Sigma \circ \beta$, where $\beta$ is a positively oriented path which connects the 
appropriate negative end of $\tilde M$ with a generic point $x \in U_j$. It follows that no component of 
$\hat \Sigma$ has a domain $U_j$ on both sides. 

For $\hat \Sigma \in \mathcal B(\Sigma)$, let 
\begin{eqnarray}
\rho(\hat \Sigma, \tilde M) = osc\{u_l\} := max_l\; \{u_l\} - min_l\; \{u_l\}
\end{eqnarray}

We introduce a subset $\mathcal B_k(\Sigma) \subset \mathcal B(\Sigma)$ based on liftings $\hat\Sigma$,
subject to $\rho(\hat\Sigma, \tilde M) = k$.
\smallskip

For a lifting $\hat \Sigma$ which separates those components of $\tilde M$ where it resides, 
$\rho(\hat \Sigma, \tilde M) = 1$. Furthermore, we have the following lemma which 
improves upon Lemma 6.6.

\begin{lem} Let $\Sigma \subset M$ realize $k[F]$---the $k$-multiple of the homology class of 
a fiber $F$. 
Then there exists a special lifting $\hat \Sigma \in \mathcal B(\Sigma)$, so that 
$\rho(\hat\Sigma, \tilde M) = k$. 

Conversely, if a surface $\Sigma \subset M$ admits a lifting $\hat \Sigma \in \mathcal B(\Sigma)$,
then $[\Sigma]$ is proportional to $[F]$ with the coefficient $\rho(\hat\Sigma, \tilde M)$.
\end{lem} 

{\bf Proof.}\quad We start with the case $k = 1$. Applying the Thom-Pontryagin construction 
to any oriented surface $\Sigma \subset M$,
produces a map $f_\Sigma : M \rightarrow S^1$.  Here $f_\Sigma$ stands for an  approximation 
of the Thom-Pontryagin map $P_\Sigma$ by a Morse function which coincides with $P_\Sigma$ in the 
vicinity of $\Sigma$. Denote by $\tilde{M}_\Sigma \rightarrow M$ 
the cyclic cover induced by $f_\Sigma$,  and by 
$\tilde f_\Sigma : \tilde{M}_\Sigma \rightarrow \R$ ---an appropriate lift of $f_\Sigma$. 
We denote by $\hat\Sigma$ the surface of constant level 
$\tilde f^{-1}_\Sigma(f_\Sigma(\Sigma)) \subset \tilde M_\Sigma$. 

Since $F$ and $\Sigma$ are homologous in $M$, 
the Thom-Pontryagin maps $f_{\Sigma} : M \rightarrow S^1$,  $f_{F} : M \rightarrow S^1$, 
produced by $\Sigma$ and $F$, are homotopic. Thus, they induce equivalent cyclic 
coverings of $M$.  The spaces of these coverings are $\Z$-equivariantly homeomorphic with
the homeomorphism $\phi : \tilde M_\Sigma \rightarrow \tilde M_F$ covering the identity map. 
We employ $\phi$ to identify the two spaces, and use the notation $\tilde M$ for
both of them. We also identify $\hat\Sigma$ with $\phi(\hat\Sigma)$. 

We denote by $t$ the upward deck translation --- a generator of the cyclic group $\Z$ 
acting on $\tilde M$. Then $\tilde M$ can be represented  as a union  
$\cup_{n\in\Z}\; t^n(\hat M_\Sigma)$. Here the fundamental region $\hat M_\Sigma$ 
is bounded by $t(\hat\Sigma)$ and $\hat\Sigma$. 
Therefore, the 2-cycle $\hat \Sigma$ is a boundary of the 3-chain $\tilde{M}^{+\infty}(\hat\Sigma) := 
\cup_{n\geq 0}\; t^n(\hat M_\Sigma)$. Hence, $\rho(\hat\Sigma, \tilde M) = 1$.\smallskip 

When  $[\Sigma] = k[F]$, a similar argument applies to $\Sigma$ and $k$ parallel 
copies of a fiber $F$. Consider the covering $\tilde M_{kF} \rightarrow M$, induced by  
the Thom-Pontryagin map $f_{kF}$. As before, there exists a lifting 
$\hat\Sigma \subset  \tilde M_{kF}$ which separates the positive and negative ends of $\tilde M_{kF}$.
The covering $\tilde M_{kF} \rightarrow M$ factors through $\tilde M \rightarrow M$. 
The fiber of $p: \tilde M_{kF} \rightarrow \tilde M$ consists of $k$ points; furthermore, 
$\tilde M_{kF}$ is homeomorphic to $k$ copies of $\tilde M$ (this point will be explained later). 
Since $\hat \Sigma \subset  \tilde M_{kF}$ separates the positive and negative ends of
$\tilde M_{kF}$,  the portion of $\hat\Sigma$, residing in each of the $k$ copies of $\tilde M$, 
separates the ends of the relevant copy. Applying the transfer
$p_\ast$  to the 3-chain $\tilde M_{kF}^{+\infty}(\hat \Sigma)$ produces 
a 3-chain which bounds $p(\hat\Sigma) \subset \tilde M$. We notice that $p(\hat\Sigma)$ 
consists of $k$ surfaces, each of which separates the ends of $\tilde M$. 
This proves the first claim.\smallskip

Assume that $[\Sigma]$ is not proportional to a fiber, but still admits a lifting from 
$\mathcal B(\Sigma)$.   Then there exists a loop $\gamma$ so that 
$[\gamma] \in Ker f_\ast$ and $[\Sigma]\circ[\gamma] \neq 0$ (cf. the proof of Lemma 6.5). 
Denote by $\hat\gamma \subset \tilde M$ 
a lift of $\gamma$. It is a loop. When $\hat \Sigma \in \mathcal B(\Sigma)$, then,  for all $n$, \, 
$t^n([\hat\gamma])\circ [\hat \Sigma] = 0$. Hence, $[\Sigma]\circ[\gamma] = 0$. 
This contradiction proves the second claim of the lemma. \qed
\smallskip

For a given $\Sigma \subset M$ comprised of several components, a lifting $\hat \Sigma \in \mathcal B(\Sigma)$
is not unique, even up to  deck translations (although, for each component, it is unique). 
For example, if a union $\hat \Sigma_0$ of a few components of $\hat \Sigma \in \mathcal B(\Sigma)$ is 
a boundary of a 3-chain in $\tilde M$, then we can apply any deck transformation
to $\hat \Sigma_0$, while leaving $\hat \Sigma \setminus \hat \Sigma_0$ intact, to produce a new lifting 
from $\mathcal B(\Sigma)$. Evidently, some restrictions on $\Sigma \subset M$
must be in place in order to claim the uniqueness, up to deck translations, of the lifting
$\hat \Sigma \in \mathcal B(\Sigma)$.\smallskip

We need to spell out the argument which we already used in the proof of Lemma 6.6. 
Let  $z^k: S^1 \rightarrow S^1$ be the canonical map of degree
$k$. For a time being, we  choose the circle with the radius $1/2\pi$. We denote by $f_k$  the composition 
of $f: M \rightarrow S^1$ with $z^k$ (so that $f_1 = f$) and let $\tilde M_k \rightarrow M$ be the 
cyclic covering induced by the $f_k$.  One can view $\tilde M_k$ as a balanced product 
$\tilde M \times_{\{T^n\}} \Z$. The cyclic $T$-action on the product $\tilde M \times \Z$ is 
defined by the formula $T(x, q) = (\tau(x),\, q - k)$,  where $x \in \tilde M$,  $q \in \Z$, and 
$\tau$ is the preferred generator of the cyclic action on $\tilde M$. The transformation 
$t: (x, q) \rightarrow (x, q + 1)$ commutes with $T$, and thus, gives rise to a cyclic $t$-action on 
$\tilde M_k$. The obvious $t$-equivariant map 
$\tilde M_k \rightarrow \Z/k\Z$ is onto and divides $\tilde M_k$ into $k$ disjoint copies of $\tilde M$.
Therefore, the natural $k$-to-1 map $p_k: \tilde M_k \rightarrow \tilde M$ splits.

Any $t$-equivariant function on $\tilde M_k$ is generated by a function 
$\tilde h: \tilde M \times \Z \rightarrow \R$,
subject to two properties: 1) $\tilde h(x, q + 1) = \tilde h(x, q) + 1$ (equivariance), and 
2) $\tilde h(\tau(x), q - k) = \tilde h(x, q)$ (being a well-defined function on the balanced product). 
In particular, the $t$-equivariant function 
$\tilde f: \tilde M \rightarrow \R$ (which covers $f$) produces a function $\tilde h_k$ with the properties
1) and 2) above: just put $\tilde h_k(x, q) = k\cdot \tilde f(x) + q$. We denote by $\tilde f_k$ 
the function on $\tilde M_k$ generated by $\tilde h_k$. 
\bigskip

Now, at least for the important case when $[\Sigma]$ is $k$-proportional to the 
homology class $[F]$ of a fiber, we will give a more conceptual 
interpretation of the twist numbers $\rho(\Sigma, \sim)$  in terms of 
the cyclic cover $\tilde M_k \rightarrow M$, induced by the map $f_k$. 
\smallskip

For each surface $\Sigma$ representing a homology class $k[F]$, 
we can measure the number of times it is "wrapped" by $f$  
around $S^1$: let
\begin{eqnarray} h(\Sigma; f) :=  1 +  min_{\hat \Sigma_k \in \mathcal B_k(\Sigma)}\; 
\big\lceil osc(\tilde f_k\big|_{\hat \Sigma_k})\big\rceil.  
\end{eqnarray}
Here $\lceil r \rceil$ stands for the integral part of a real number 
$r$. Abusing previous notations, $\mathcal B_k(\Sigma)$ in (6.2) denotes 
the set of $\Sigma$-liftings which separate the positive and negative ends of $\tilde M_k$. \smallskip 

The following definitions aim to introduce  notions of \emph{breadth} and \emph{height} 
of a given surface $\Sigma \subset M$ relative to a given map $f: M \rightarrow S^1$. 
They rely on Lemma 6.6. 
We always assume that surfaces $\Sigma$ and $F$ are in general position. 

%%%%
\begin{defn} 
For a given map $f: M \rightarrow S^1$, let $\Sigma \subset M$ be a surface representing the $k$-multiple of the 
homology class of a fiber. Let $F$ be a union of a few fiber components  (not necessarily belonging to the same fiber).  
Denote by $\mathcal A(F)$ the set of \emph{all} liftings $\{\hat F\}$ of $F$ to the space $\tilde M_k$.
Each surface $\hat F$ intersects only with finitely many copies 
$\{t^n(\hat \Sigma)\}_n$ of the surface $\hat \Sigma_k \subset \tilde M_k$, where $\hat \Sigma_k \in \mathcal B_k(\Sigma)$. 
We minimize the number of such copies over the set $\mathcal A(F)$, denote the minimum by $b(F, \Sigma)$, 
and call it the \emph{breadth}  of $F$ relative to $\Sigma$. 
\end{defn}
%%%%%%
Note that, for any $f$-fiber $F'$ in general position with $F$,\, $F \cap F' = \emptyset$. 
Thus, $b(F, F') = 0$.\smallskip

At this point, it is not 
clear why $b(F, \Sigma)$ does not depend on the choice of the special lifting 
$\hat\Sigma_k \in \mathcal B_k(\Sigma)$. Definition 6.7 will be justified by  
linking directly  $b(F, \Sigma)$ with the quantity $\rho(\Sigma, F)$ which is independent 
on the liftings of $\Sigma$. 

%%%%
\begin{prop} Assume that $\Sigma \subset M$  represents a $k$-multiple of the homology class of a fiber
and let $F$ be any union of fiber components. Then $\rho(\Sigma, F) = b(F, \Sigma)$.
\end{prop}
%%%%%

{\bf Proof.}\quad For a lifting $\hat\Sigma_k \in \mathcal B_k(\Sigma)$, let $u_{\hat \Sigma_k}$ 
denote a step function which takes  value 0 at the points 
of the domain $\tilde{M}^{-\infty}_k(\hat\Sigma_k)$ and value 1 at the points of 
$\tilde{M}^{+\infty}_k(\hat\Sigma_k)$. Recall that $\hat\Sigma_k$ is the boundary of the 3-chain 
$\tilde{M}^{+\infty}_k(\hat\Sigma_k)$.\smallskip

Consider all the surfaces $\{t^n(\hat\Sigma_k)\}_n$ having 
a non-empty intersection with a particular lifting $\hat F_k \subset \tilde M_k$ of $F$.  
Restrict the potential function 
$$u\; := \;\sum_{\{n|\; t^n(\hat\Sigma_k)\, \cap \, \hat F_k \;\neq\; \emptyset\}}\; u_{ t^n(\hat \Sigma_k)}$$ 
to $\hat F_k$. 
Crossing in $\hat F_k$ an oriented curve from $\hat{\mathcal C}_n = \hat F_k\cap t^n(\hat{\Sigma_k})$ 
in the positive normal direction is  the same as crossing in $\tilde M_k$ the 
corresponding component of $\hat\Sigma_k$ in the preferred normal direction: both have 
the effect of increasing the potential $u$ by 1. Therefore, $u$ gives rise to 
a 0-cochain on the oriented graph dual to the pattern $\hat{\mathcal C}_n$ in $\hat F_k$.
It takes at most two values.

We observe that, since $p: \hat \Sigma_k \rightarrow \Sigma$ and 
$p: \hat F_k \rightarrow F$ are 1-to-1 maps, the covering map 
$p: \tilde M_k \rightarrow M$ defines a diffeomorphism
of pairs $\coprod_{n}\; \hat{\mathcal C}_n \subset \hat F_k$ and  $\Sigma \cap F \subset F$, 
where the disjoint union employs all the non-empty $\hat{\mathcal C}_n$'s. 
In particular,  the images ${\mathcal C}_n := p(\hat{\mathcal C}_n)$ of distinct intersections 
$\hat{\mathcal C}_n$'s are disjoint in $M$.  
  
Thus, $u$ also produces a 0-cochain on the oriented graph, dual to the pattern 
$\Sigma \cap F$ in $F$ (equivalently, an integral 2-chain on $F$ whose boundary is $\Sigma \cap F$). 
It takes exactly as many consecutive values as the number of patterns $\{\hat{\mathcal C}_n \neq \emptyset\}$.
Hence, $osc(u|_F)  \geq \rho(\Sigma, F) + 1$.  

For each component $F^\beta$ of $F$, the potential $u$ is determined, up to a choice of a constant,
by the oriented intersection $F^\beta \cap \Sigma \subset F^\beta$. Therefore, the oscillation 
$osc(u|_{F^\beta}) = \rho(\Sigma, F^\beta) + 1$. 

Next, we minimize $osc(u|_F)$ by independently applying deck translations in $\tilde M_k$ to various
components $\hat F^\beta_k$ to achieve the equality 
$osc(u|_F) - 1 = min_\beta\, \{\rho(\Sigma, F^\beta)\} := \rho(\Sigma, F)$. 
This can be done by moving each component $\hat F^\beta_k$ above $\hat\Sigma_k$ and 
so that it has a non-empty intersection with the fundamental domain bounded by 
$\hat\Sigma_k$ and $t(\hat\Sigma_k)$. 
As a result, $b(F, \Sigma) = \rho(\Sigma, F)$. \qed \bigskip

Let $\mathcal S_{{\chi_-}}^k$ be the set of surfaces $\Sigma \subset M$ which realize the Thurson norm of  
$k[F]$, where $[F]$ stands for the homology class of an $f$-fiber. Among $\Sigma \in \mathcal S_{{\chi_-}}^k$ 
consider  surfaces with the minimal breadth  $b(F_R, \Sigma)$. We denote  this optimal breadth by 
$b_{{\chi_-}}(F_R, k[F])$. When $k = 1$,
we also use an abbreviation $b_{{\chi_-}}(f)$  for $b_{{\chi_-}}(F_R; [F])$. \smallskip 

Combining Proposition 6.8 with Lemma 6.2, we get 
%%%%
\begin{cor} Let $R$ be the repeller set for a map $f: M \rightarrow S^1$ 
and $[F]$ the homology class of a fiber.
Then \hfil\break \centerline{\qquad \qquad $b_{{\chi_-}}(F_R, k[F]) = \rho_{{\chi_-}}(k[F], F_R) = 
min_{\Sigma \in \mathcal S_{\chi_-}^k} \;\{\rho_{{\chi_-}}^\circ(\Sigma, F_R)\}$.\qed} 
\end{cor}
%%%%

%%%%
\begin{lem} Let $\Sigma \subset M$ represent a $k$-multiple of the  homology class of a fiber
and let $F$ be any finite union of fiber components.
Then $\rho(\Sigma, F)  \leq h(\Sigma, f)  - \epsilon$, where $\epsilon = 0, 1$ depending on a 
particular location of  $F$ in $M$.  

When $F$ is a fiber and $\Sigma$ is connected (hence, $k = 1$), 
then $\rho(\Sigma, F)  = h(\Sigma; f) - \epsilon$, which implies a very weak dependence of $\rho(\Sigma, F)$ 
on the fiber $F$.
\end{lem}
%%%%%
{\bf Proof.}\quad Take $\hat\Sigma_k$ which separates the positive and negative ends of $\tilde M_k$ 
and delivers $h(\Sigma; f)$.
Let $\hat F_k \subset \tilde M_k$ be a lifting of $F$ which, together with $\hat\Sigma_k$, delivers 
$b(F, \Sigma) = \rho(\Sigma, F)$ (as described in the proof of Proposition 6.8). Denote by $F^\beta$  
a typical connected component of $F$.

Using the cyclic $t$-action on $\tilde M_k$, the set   
of $n$'s for which $t^n(\hat\Sigma_k) \cap  \hat F^\beta_k \neq \emptyset$ is a reflection with respect to 0 
of the set of $n$'s for which $\hat \Sigma_k \cap t^n(\hat F_k^\beta) \neq \emptyset$. Therefore,
the cardinality of such $n$'s does not exceed $osc\{n:\; \hat \Sigma_k \cap t^n(\hat F_k^\beta) \neq \emptyset\}$.
Since, for each $n$,\, $t^n(\hat F_k^\beta)$  belongs to a constant level set of the function $\tilde f_k$, and, 
for distinct $n$'s, these levels are integrally spaced,  
$$osc\{n:\, \hat \Sigma_k \cap t^n(\hat F_k^\beta) \neq \emptyset\} \leq \epsilon + 
\lceil osc(\tilde f_k|_{\hat\Sigma_k})\rceil := h(\Sigma, f) - \epsilon.$$ Therefore, 
$\#\{n:\,t^n(\hat\Sigma_k) \cap  \hat F^\beta_k \neq \emptyset\} \leq  h(\Sigma, f)$.
Due to Proposition 6.8 and its proof, $\rho(\Sigma, F)$  is the maximum over all $\beta$'s of the LHS of the 
previous inequality. Hence, $\rho(\Sigma, F) \leq h(\Sigma, f) - \epsilon$. 

When $\Sigma$ is connected and $F$ is a fiber, the same 
arguments show that $\rho(\Sigma, F) = h(\Sigma, f) - \epsilon$. Indeed, the connectivity of $\Sigma$ forces 
it to cross all the "intermediate floors" $\{t^n(\hat F_k)\}$ between the top and the bottom one. \qed
\bigskip

%%%%%%%%
\begin{defn} Employing (6.2), put \; 
$h_{\chi_-}(k[F], f) = min_{\Sigma \in \mathcal S_{{\chi_-}}^k}\; \{h(\Sigma, f)\}$. 
We also use abbreviation $h_{\chi_-}(f)$  for $h_{\chi_-}([F], f)$.
\end{defn} 
%%%%%%%%

Crudely, the difference between height and breadth is like the difference 
between the degree and the number of non-zero monomials in a Laurent polynomial 
from the ring $\R[t, t^{-1}]$.\smallskip

%%%%%%
\begin{figure}[ht]
\centerline{\BoxedEPSF{slanted.A scaled 550}}
\bigskip
\caption{}
\end{figure} 
%%%%%%%%
When $\Sigma$ admits an $f$-positive loop $\gamma$ which hits it only once, 
then  one can deform $f$ in such a way that
$\Sigma$ will be of an arbitrary big height and breadth with respect to the deformed 
map (cf. Example 8.12). 
\bigskip

The following proposition is a tautology worth mentioning because its converse will be a focus 
of our efforts in Section 8. It follows from the observation that, if one can find an 
${\chi_-}$-optimal surface $\Sigma$ among fibers or fiber components of a given map $f$, then 
you can find it being disjoint from any other fixed union $F$ of fibers or fiber components.
Thus, $\rho(\Sigma, F) = 0$.  When  $b(F, \Sigma)$  makes sense, it vanishes as well.

\begin{cor} 
Let  $f: M \rightarrow S^1$ be any Morse map as in Lemma 4.2. 
\begin{itemize}
\item if a finite union of fiber components delivers the minimal value of 
${\chi_-}(\sim)$ in a vertical homology class $[\Sigma]$, then 
$\rho_{{\chi_-}}([\Sigma], F_R) = 0$;
\item if a finite union of fiber components delivers the minimal value of 
${\chi_-}(\sim)$ in the homology class $k[F]$, then 
$b_{{\chi_-}}(F_R; k[F]) = 0$; 
\item  if a finite union of fiber components delivers the minimal value of 
${\chi_-}(\sim)$ in the homology class $k[F]$, then 
$h_{{\chi_-}}(k[F], f) = 1.$ \qed 
\end{itemize} 
\end{cor}
\smallskip

By [T], compact leaves of taut (that is, no generalized Reeb components) foliations 
deliver the Thurston norm of their homology classes. Therefore,

\begin{cor} 
Let a compact oriented surface $\Sigma \subset M$ be a union of leaves of a \emph{taut}  
smooth foliation $\mathcal F$ whose leaves are transversal to $\partial M$.  
Consider a Morse map $f: M \rightarrow S^1$  homotopic to the Thom-Pontryagin map 
$f_\Sigma : M \rightarrow S^1$, the homotopy being an identity on $\Sigma$.
Then $\rho_{{\chi_-}}(f) = b_{{\chi_-}}(f) =  0$ and $h_{{\chi_-}}(f) = 1$. \qed
\end{cor}
\smallskip

Although an effective computation of the invariants  
$b_{{\chi_-}}(F_R; k[F])$, $b_{{\chi_-}}(f)$, $h_{{\chi_-}}(f)$ 
seems to be as difficult as the computation of the norm $\|[F]\|_T$, one has a good 
grip on how these invariants might change under an $S^1$-controlled homotopy of the map $f$. 
 
%%%%%%%
\begin{lem} Let $S^1$ be a circle with the circumference 1.  Assume that a homotopy 
$\{f_\tau: M \rightarrow S^1\}_{0 \leq \tau \leq 1}$  
is such that each $\tau$-path $f_\tau(x)$, $x\in M$, winds less than $q$ times around the circle. 
If $h_{{\chi_-}}(k[F], f_0) \leq l$, then 
$h_{{\chi_-}}(k[F], f_1) \leq l + 2kq$. 

When the image $f_\tau(x)$ of any point $x \in M$ moves clockwise, then 
$h_{{\chi_-}}(k[F], f_1) \leq l + kq.$ Hence, $\rho_{{\chi_-}}(k[F], F_{R_1}) \leq l + kq$ as well.
\end{lem}
%%%%%%
{\bf Proof.}\quad 
Let $\Sigma^\star \subset M$ be a surface which delivers the minimal value of  ${\chi_-}(\Sigma^\star)$  
in the  $k$-multiple of the homology class of a fiber. Due to Lemma 6.6, $\Sigma^\star$  
admits special lifting $\hat\Sigma^\star_k \in \mathcal B_k(\Sigma^\star)$. In
addition, assume that $\hat\Sigma^\star_k$ is chosen so that \hfil\break
\centerline{$1 + \lceil osc(\tilde f_{0,\, k}\big |_{\hat\Sigma^\star_k})\rceil = h(k[F], f_0)  \leq l$.} 

For any $\Sigma$ in the homology class $k[F]$, consider the function
$\Phi(\tau, \hat \Sigma_k) := \hfil\break osc (\tilde f_{\tau,\, k}\, \big|_{\hat\Sigma_k})$, where $\hat \Sigma_k \in 
\mathcal B_k(\Sigma)$ and $\tilde f_{\tau,\,k}: \tilde M_k \rightarrow \R$ covers the homotopy
$f_\tau$. Note that, for any $x \in \tilde M_k$, 
$|\tilde f_{1,\,k}(x) - \tilde f_{0,\, k}(x)| < kq$. Hence, for any lifting 
$\hat\Sigma_k \in \mathcal B_k(\Sigma)$,\,
$|\Phi(1, \hat \Sigma_k) - \Phi(0, \hat \Sigma_k)| <  2kq$, implying 
$\Phi(1, \hat \Sigma_k) <  \Phi(0, \hat \Sigma_k) + 2kq$. 
Thus, $h(\Sigma^\star, f_1) := $ $$1 + min_{\hat \Sigma_k^\star \in \mathcal B_k(\Sigma^\star)} 
\lceil  \Phi(1, \hat\Sigma_k^\star) \rceil \leq 
1 + \lceil\Phi(0, \hat \Sigma_k^\star)\rceil + 2kq = h(\Sigma^\star, f_0) + 2kq.$$ 
Therefore, $h(k[F], f_1) \leq h(\Sigma^\star, f_0) + 2kq : = h(k[F], f_0) + 2kq$.

A similar argument is valid for a "clockwise" homotopy. \qed
%%%%%%
\begin{cor}  
The $S^1$-controlled size of 
a homotopy, which links a given map $f$ to a map with a ${\chi_-}$-minimizing fiber or 
a union of fiber components, gives an upper bound on $h_{{\chi_-}}(f)$, and thus, on 
$\rho_{{\chi_-}}(f)$. 
\end{cor}
%%%%%%%%
{\bf Proof.}\quad Let  $\Sigma \subset M$ be a surface minimizing ${\chi_-}(\sim)$
in the homology class of an $f$-fiber. Consider any Morse approximation $f_1$ of 
the Thom-Pontryagin map $f_\Sigma$ (it is homotopic to $f$) such that 
$f_1 |_{\Sigma} = f_\Sigma |_{\Sigma} = pt$. By Corollary 6.13, $h_{{\chi_-}}(f_1) = 1$. 
By compactness argument, there exists a minimal natural number $q$, so that the image of 
any point in $M$, under the homotopy linking $f$ with $f_1$, winds less than $q$ times 
around the circle. By Lemma 6.14, $\rho_{{\chi_-}}(f) \leq h_{{\chi_-}}(f) \leq 1 + 2q$.  
\qed
\bigskip

%%%%%%%%%%%%%%%%%%%%%%%%%%%%%%%%%%%%%%%

\section{Resolving intersections with fibers} 

The main thrust of the arguments below is influenced by the proofs of Theorem 1
in [T] and of Theorem 2.3 in [H]. \smallskip

Given a Morse map $f: M \rightarrow S^1$, consider an embedded oriented surface 
$(\Sigma, \partial\Sigma )  \subset (M, \partial M)$ representing 
a vertical homology class. Let $F$ be a finite union of fiber components in 
general position with respect to $\Sigma$.
Denote by  $\{C_i\}$ the components of the intersection $\mathcal C = \Sigma \cap F$. 
As before, each $C_i$ is given an orientation by its preferred normal in 
$F$, which coincides with the preferred  normal of $\Sigma \subset M$. 
\smallskip

In a small neighborhood of each component $C_i$, the surfaces $\Sigma$ 
and $F$ divide $M$ in \emph{four} regions.
The orientations of $F$ and $\Sigma$ pick a unique pair of non-adjacent
quadrants, say I and III, along a typical intersection curve $C_0$ as shown in
Figures 11 and 12.

Along $C_i$, we can \emph{resolve} the intersection of $F$ and $\Sigma$ in a
unique way, as shown in  Figure 12, diagram A. 
The resolution $F \bowtie_i \Sigma$ will occupy a pair of non-adjacent 
quadrants. The resolved surface $F \bowtie_i \Sigma$ inherits the normal frames
of $\Sigma$ and $F$.
These local resolutions $F \bowtie_i \Sigma$ can be pasted into a well-defined
oriented surface 
$F \bowtie \Sigma$, homologous to $[\Sigma] + [F]$ (cf. Figure 11). 
%%%%%%%
\begin{figure}[ht]
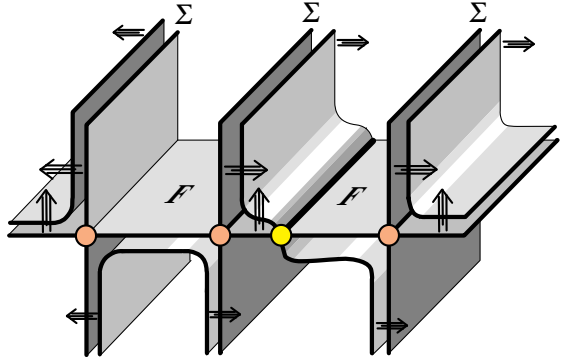

\centerline{\BoxedEPSF{resolution scaled 400}}
\bigskip
\caption{Resolving the intersection of $\Sigma$ and $F$. Note, the new
intersection---the bold line---of
the resolution $F \bowtie \Sigma$ with $F$.}
\end{figure}
%%%%%%% 
As we are trying to paste $F \bowtie_i \Sigma$ together, we eliminate all the old
intersection curves from $F \cap \Sigma$ 
and often are forced to introduce new intersections of $F \bowtie \Sigma$ with $F$.
This happens because, 
over each connected component $F_j^\circ$ of $F \setminus (\Sigma \cap F)$, some
of the "germs" of 
$F \bowtie_i \Sigma$ will reside \emph{above} $F$, and some \emph{below} it.

Consider, for instance resolving $m$ 
coherently oriented meridians on a torus $T^2 = F$. The new surface $F \bowtie_i
\Sigma$ still will have $m$ 
intersection loops with the torus. In this example, the resolution does not help
to simplify the intersection 
pattern. However, if at least two meridians have opposite orientations, the
simplification becomes possible. 

We intend to show that the new intersections are simpler than the original ones,
and that, through iterations of resolutions, they can be eventually eliminated, 
thus, producing a new embedded surface which does not intersect $F$ at all. 

In addition to the canonical resolution (diagram $A$ in Figure 12), we
also will use its modifications, shown in diagrams $B$, $C$.   

%%%%
\begin{figure}[ht]
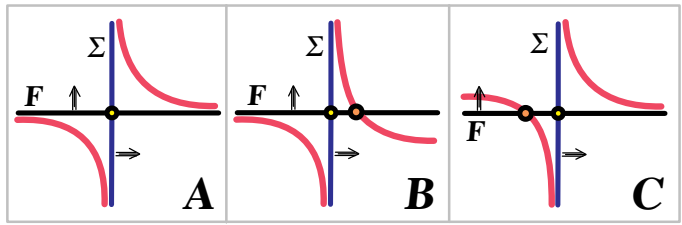

\centerline{\BoxedEPSF{ABC scaled 480}}
\bigskip
\caption{Three ways of resolving the intersection} 
\end{figure} 
%%%%%%

Each time we resolve the intersection $\mathcal C = F \cap \Sigma$, along a 
particular component $C_i$ of $\mathcal C$, we delete a regular neighborhood 
$U_i \approx C_i\times [-\epsilon, +\epsilon]$ of $C_i$ from $\Sigma$ and a 
regular neighborhood $V_i \approx C_i\times [-\delta, +\delta]$ of $C_i$ from $F$ 
and repaste the four components $\{C_i \times \pm \epsilon, C_i \times \pm \delta\}$ 
in a new way. Comparing the Euler numbers of $\Sigma \setminus \coprod_i U_i$ and 
$F \setminus \coprod_i V_i$ with the ones of $F \bowtie \Sigma$ and $F \coprod \Sigma$, 
we see that $\chi(F \bowtie \Sigma) = \chi(F) + \chi(\Sigma)$.

To describe the relation between  
$\chi_-( F \bowtie \Sigma )$ and $\chi_-(F) + \chi_-(\Sigma)$ we will use a few 
ideas from [T], pages 103-104. In general, the desired additivity 
$\chi_-( F \bowtie \Sigma ) = \chi_-(F) + \chi_-(\Sigma)$  is 
upset by the new spherical or disk components generated as a result of the resolution.
However, there are situations where the emergence of new spherical components 
can be prevented.\smallskip

Recall, that it is possible to perform 2-surgery on $\Sigma$ using the 2-disks bounding 
in $F$ the loops and arcs from $\Sigma \cap F$. If we perform the surgery starting 
with the most "inner" disks in $F$ and gradually moving "outwards", the resulting surface 
$\Sigma^\odot$ is, up to the obvious isotopies, unique.  Its intersection 
with $F$ is free of loops which bound disks in $F$. 

\begin{defn} We say that $\Sigma \subset M$ is \emph{well-positioned with respect to} 
$F \subset M$, if the transversal intersection $\Sigma^\odot \cap F$ has no components which 
bound a disk in $\Sigma^\odot$.

In particular, if $F \subset M$ is an \emph{incompressible} surface, then any $\Sigma$ 
is well-positioned with respect to it.
\end{defn}

Note that if $\Sigma \cap F$ contains a loop which bounds a disk in $\Sigma$, but not in $F$, 
then $\Sigma$ is not well-positioned with respect to $F$. However, even if no such loop exists, it 
is still possible that the intersection $\Sigma^\odot \cap F$ will contain a  loop which 
bounds a disk in $\Sigma^\odot$, but not in $F$.

\begin{lem} Let $\Sigma$ represent an $f$-vertical homology class and $F$ be a union of 
fiber components. If $\Sigma$ is well-positioned with respect to $F$,
then the  surface $\Sigma^\odot \subset M$ has the following properties: 1) $\Sigma^\odot$ is cobordant 
to $\Sigma$; 2) $\chi_-(\Sigma^\odot) \leq \chi_-(\Sigma)$; 3) 
$\chi_-( F \bowtie \Sigma^\odot) = \chi_-(F) + \chi_-(\Sigma^\odot)$. 

In particular, if $\rho^\circ(\Sigma, F) = 0$\footnote{that is, $\Sigma^\odot \cap F = \emptyset$}, 
then $\Sigma$ is well-positioned with respect to $F$.
\end{lem}

{\bf Proof.}\quad The new surface $\Sigma^\odot$ has properties 1) and 2) (cf. Lemma 3.1). Its intersection 
$\Sigma^\odot \cap F$ is free of components bounding a disk in $F$.  We claim that no \emph{new} spherical 
or disk component $S$ are present in 
$F \bowtie \Sigma^\odot$. Indeed, such an $S$ would be a union $S_{\Sigma^\odot} \cup S_F$ of two surfaces
whose common boundary is a subset of $\Sigma^\odot \cap F$. Here  $S_F$ 
is homeomorphic to a union of some domains in which $\Sigma^\odot \cap F$ divides $F$, and $S_{\Sigma^\odot}$ 
is homeomorphic to a union of some domains in which $\Sigma^\odot \cap F$ divides $\Sigma^\odot$. 
Note that, if a sphere is divided into two complementary domains, at least one of them must contain a 
boundary component which bounds in that domain a disk---a connected domain in the plane has a non-positive 
Euler number, unless it is a disk. Since $\Sigma^\odot \cap F$ does not bound a disk in 
$F$ (disk bounding components have been eliminated by the 2-surgery) and in $\Sigma^\odot$ 
(by being well-positioned), it is impossible to generate a new spherical or disk component $S$. \qed
\smallskip 

{\bf Remark.}\quad 
Here is the only point, where a parallel program for the genus invariants faces similar
but more serious difficulties.  Unless the number of components in $F \bowtie \Sigma$
is less than or equal to the number of components  in $F \sqcup \Sigma$, the desired 
equality  $g(F \bowtie \Sigma ) \leq  g(F) + g(\Sigma)$  is not valid. Unfortunately, 
in general, we do not know how to control 
the number of components in $F \bowtie \Sigma$. However, if $F \cap \Sigma$ is
connected, then $g(F \bowtie \Sigma ) =  g(F) + g(\Sigma)$.  \qed 
\smallskip 

Now, we use the dual graph $K_{\mathcal C}$ of $\mathcal C \subset F$ and its 
modifications as  bookkeeping
devices for recording  the resolutions of the intersection $\Sigma \cap F$.
We will subject $K_{\mathcal C}$ to elementary modifications, which
will mimic particular ways (see Figure 12) 
of resolving the intersection $\mathcal C$. Some modifications will erase a 
curve from $\mathcal C$, will eliminate the  
corresponding edge in $K_{\mathcal C}$ and will merge the two vertices 
it connects into a single one. 
\smallskip

First, we eliminate the edges which 
correspond to loops bounding a disk in $F$, or to arcs bounding a disk in 
$(F, \partial F)$. This elimination  
corresponds to 2-surgery, as described in the beginning of the proof of 
Theorem 5.2 where we have eliminated a totally 
reducible  pattern $\mathcal C$, while keeping the invariants $g(\Sigma)$,
$\chi_-(\Sigma)$ on a decline. 

After that, we consider the curves from the modified $\mathcal C = \mathcal C^\circ$, which 
correspond to the edges emanating from the vertices 
with the maximal value $m$ of the potential $u$. Then, we can perform the
$A$-type resolutions along them. 
Next, we perform the $B$-type resolutions along the rest of the intersection
curves. This has an effect 
on $K_{\mathcal C}$ of eliminating the edges of the maximal type, merging all
vertices of the level 
 $u = (m - 1)$ with the appropriate vertices of the $m$-level, and keeping the
rest of 
the graph unchanged. The modified graph $K_{\tilde{\mathcal C}}$ will have a
"truncated" level function 
$\tilde u$ of its own. Moreover, $u$, being restricted to the portion $K_{\mathcal
C}^{< m}$ of $K_{\mathcal C}$ 
below $m$, is the pull-back of $\tilde u$ under the obvious map  
$K_{\mathcal C}^{< m} \rightarrow K_{\tilde{\mathcal C}}$. Hence, $\rho
(\tilde{\mathcal C}) = \rho (\mathcal C) - 1$.
The resulting surface $\Sigma^\odot \bowtie F$ is homologous to $[\Sigma] + [F]$ and has an
intersection with $F$ described by the graph 
$K_{\tilde{\mathcal C}}$. 

Since the intersection $(\Sigma^\odot \bowtie F) \cap F$ 
consists of curves which are isotopic to the original curves from 
$\Sigma \cap F$ (cf. Figures 11, 12), and since we have excluded all the curves 
which bound a disk from our original intersection $\mathcal C$, the new 
intersection $(\Sigma^\odot \bowtie F) \cap F$ is free of disk-bounding curves as well. 

This recipe can be repeated again and again until, after $\rho^\circ(\mathcal C)$
iterations, we eliminate the 
intersection with $F$ completely. In this algorithm, the potential $u$ helps to 
paste the local resolutions (of the $A, B$ and $C$-types) together. The final
surface $\Sigma_{\check{F}}$ resides in 
the homology class of $[\Sigma] + \rho^\circ(\mathcal C) \cdot [F]$.

After $\rho^\circ(\Sigma , F)$ resolutions,
$\chi_-(\Sigma_{\check{F}}) \leq \chi_-(\Sigma) + \rho^\circ(\Sigma , F)\cdot \chi_-(F)$,
provided that $\Sigma$ was well-positioned with respect to $F$.

When $\Sigma$ is not well-positioned, we need to add a correction term to the RHS of the 
inequality above. This correction term $\mu^\circ(\Sigma, F)$ equals twice the 
number of "new" spheres plus the number of "new" disks present in $\Sigma_{\check{F}} = 
\Sigma^\odot \bowtie \{\rho^\circ(\Sigma , F)\cdot F\}$,  
but not in $\Sigma^\odot$ or in $\rho^\circ(\Sigma, F) \cdot F$. In other words,
%%%
\begin{eqnarray}
\quad \mu^\circ(\Sigma, F) = \chi_+(\Sigma^\odot \bowtie \{\rho^\circ(\Sigma, F)\cdot F\}) - 
\chi_+(\Sigma^\odot) - \chi_+(\rho^\circ(\Sigma, F)\cdot F)
\end{eqnarray}

Each new sphere will consume at least two disks bounding a curve from $\Sigma^\odot \cap F$ which 
bounds a disk in $\Sigma^\odot$ but does not bounds a disk in $F$. Each new relative disk 
will require at least one such curve. Let us denote by $\nu^\circ(\Sigma, F)$
the number of such curves. Then $\mu^\circ(\Sigma, F) \leq 2 \nu^\circ(\Sigma, F)$. 

In special cases we can rule out the emergence of new spheres just from observing the intersection 
pattern $\Sigma \cap F$ in $F$.  Let $U$ be one of the connected domains in which 
$\Sigma^\odot \cap F$ divides $F$. Imagine that $U$ contains a handle or, what is the same, that 
$d < 2 - \chi(U)$, where $d$  is the number of components in $\partial U$. Evidently, such an $U$ can 
not contribute to a sphere in $\Sigma^\odot \bowtie F$. Moreover, further resolutions can only enlarge
$U$, thus, preserving the handle inside $U$. Such a case is described in a model Example 8.12 and 
depicted in Figure 15.

We assemble these observations in
%%%%
\begin{lem} For any finite union $F$ of fiber components and any oriented surface $\Sigma \subset M$, 
representing a vertical 2-homology class, 
there exists an embedded surface $\Sigma'$ with the properties:
\begin{itemize}
\item $\Sigma'$ is homologous to the cycle $[\Sigma] + \rho^\circ(\Sigma, F) [F]$  
\item $\chi_-(\Sigma') \leq \chi_-(\Sigma) + \rho^\circ(\Sigma , F)\cdot
\chi_-(F) + \mu^\circ(\Sigma, F)$
\item $\Sigma' \cap F = \emptyset$
\end{itemize}
If $\Sigma$ is well-positioned with respect to $F$, then $\mu^\circ(\Sigma, F) = 0$. This is the case 
when $\rho^\circ(\Sigma , F) = 0$ (i.e. when $\Sigma^\odot \cap F = \emptyset$), or when $F$ is an incompressible  
surface, or when $\Sigma^\odot
\cap F$  divides $F$ into domains, each of which contains a handle. \qed
\end{lem}

%%%%%%%%%%%%%%%%%%
%%%%%%%%%%%%%%%%%%%%%

\section{Twist, variation and the $\chi_-$-optimization}  

We are in position to derive our main results. Basically, we follow the train of thought 
presented in Section 5, but now we will bring to the game the twist and height invariants from Section 6 
and the graph-theoretical considerations from Sections 3 and 4.
\smallskip

Consider an attractor set $A$ and a repeller set $R$ in $\Gamma_f$. \smallskip

Given a formal  combination $\sum_i \kappa_i F_i$ ($\kappa_i \in \R$), of $f$-oriented fiber components $\{F_i\}$, define its 
norm $\|\sum_i \kappa_i F_i\|$ by the formula $\sum_i |\kappa_i| \cdot \chi_-(F_i)$.\smallskip
 
In a similar way,  we introduce a $\chi_-$-weighted $l_1$-semi-norm $\|\sim\|_A$ on $\R[A]$ by the formula 
%%%%%
\begin{equation}
\|\kappa\|_A = \sum_{a \in A}\; |\kappa(a)|\cdot\chi_-(F_a)
\end{equation}
%%%%%
Here $\kappa\in \R[A]$ and $F_a$ denotes the fiber component corresponding to a
point $a \in A \subset \Gamma_f$. When all $\chi_-(F_a) \neq 0$, the unit ball in this norm is a convex
hull   spanned by the vectors $\{\pm \chi_-(F_a)^{-1}a\}_{a \in A}$. 

When  $\Gamma_f$ admits a tree decomposition $\sqcup_{r \in R}\; T_r^{\pm}$ as in Lemma 4.2, 
thanks to Lemma 4.6, we  have an epimorphism $P: \R[A] \rightarrow H_2^f\otimes\R$.  
One can combine (8.1) with the $P$ to define a "vertical" semi-norm 
$\|\sim\|_{H^f}$ on $H_2^f \subset H_2(M, \partial M; \R)$  
by the formula\footnote{It is easy to verify that the LHS of (8.2) satisfies the triangle inequality.}:
\begin{equation}
\|[\Sigma]\|_{H^f} = min_{\{\kappa \, \in\,  \R[A]|\;\; P(\kappa) \;=\; [\Sigma]\}}\; \{\|\kappa\|_A\} 
\end{equation}\smallskip
Alternatively, consider
%%%%%% 
\begin{defn}
The semi-norm $\|[\Sigma]\|_{H^f}$ of a vertical homology class $[\Sigma]$ can be defined as the 
minimum of $\chi_-(\sim)$, taken over all finite oriented unions of distinct fiber components which 
represent $[\Sigma]$.
\end{defn} \smallskip

The equivalence of this less technical definition with the one given by (8.2) follows 
from Lemma 4.2 and 4.6 coupled with a familiar  observation:  replacing any fiber component by 
a union of components indexed by $A$ can be accomplished via 2-surgery --- an operation which 
decreases the value of $\chi_-(\sim)$. 
\smallskip

\smallskip 

Fix an $f$-vertical homology class $[\Sigma] \in H_2(M, \partial M; \Z)$ represented by a $\Z$-linear combination
of the cycles $\{[F_a]\}_a$ with coefficients $\{\alpha_a\}$. Let $\alpha_\Sigma$ 
stand for a function $A \rightarrow \Z$ defined by the formulas $\{\alpha_\Sigma(a) = \alpha_a\}$.

For any probe surface $\Sigma$ representing $[\Sigma]$, denote by  $\rho_\Sigma : R \rightarrow \Z_+$ a 
function which assigns to each element $r \in R$ the value $\rho^\circ(\Sigma, F_r)$, where $F_r$ 
stands for the fiber component corresponding to the point $r \in \Gamma_f$. \smallskip

We pick a basis $\{C_k\}$ of 1-cycles in $H_1(\Gamma_f; \Z)$. 
Given a function $\delta: \Gamma_f \rightarrow \Z$ with a finite support located in the complement to the vertices 
of $\Gamma_f$, denote by $\int_{C_k} \delta$ the sum 
$\sum_{y \in supp(\delta)} \delta(y)\cdot [F_y \circ C_k]$, where $[F_y \circ C_k]$ stands for the
(algebraic) intersection number between the loop $C_k$ and  the fiber component $F_y$ indexed by $y \in \Gamma_f$.
\smallskip

%%%%
\begin{figure}[ht]
\centerline{\BoxedEPSF{norm scaled 520}} 
\bigskip
\caption{} 
\end{figure} 
%%%%%%
Let  $\mathcal V_\Sigma \subset \Z[A]$ be an affine sublattice defined by a system linear equations: 
\begin{equation}
\mathcal V_\Sigma:= \bigg \{ \kappa\in \Z[A]:\quad \int_{C_k} \kappa = 
\int_{C_k} \alpha_\Sigma + \int_{C_k}\rho_\Sigma \bigg \}, 
\end{equation}
where the loops $\{C_k\}$ form a  basis in $H_1(\Gamma_f; \Z)$. 
We notice that the first integral on the RHS of (8.3) depends 
only on the homology class $[\Sigma]$, while the second integral depends on its 
particular representative $\Sigma$, or rather on the function $\rho_\Sigma$ which $\Sigma$ defines
on $R$. 

Similarly, put
%%%%%
\begin{equation}
\mathcal W_\Sigma:= \bigg \{ \kappa\in \Z[A]:\quad \int_{C_k} \kappa = \int_{C_k} \alpha_\Sigma \bigg \} 
\end{equation}
%%%%%
By Lemmas 4.5, 4.6, the subgroup $H_2^f \subset H_2(M, \partial M; \Z)$ of vertical homology classes is 
isomorphic to the quotient $\Z[A]/\mathcal T_A$, where 
$$\mathcal T_A := \Big\{\kappa \in \Z[A] : \int_{C_k} \kappa = 0\Big\}.$$

We will be interested in special elements $\kappa^{\mathcal V_\Sigma} \in \mathcal V_\Sigma$ and 
$\kappa^{\mathcal W_\Sigma} \in \mathcal W_\Sigma$ 
which minimize the $\|\sim\|_A$-norm. By its  very definition, $\|\kappa^{\mathcal W_\Sigma}\|_A = 
\sum_a |\kappa^{\mathcal W_\Sigma}_a| \cdot \chi_-(F_a)$ is the minimal 
$\chi_-$-characteristic among all $\Z$-combinations of fiber components $\{F_a\}$ which realize the given 
homology class $[\Sigma]$, that is, $\|\kappa^{\mathcal W_\Sigma}\|_A = \|P(\kappa^{\mathcal W_\Sigma})\|_{H^f}$.
Recall that each fiber
component $F$ is cobordant to a union of a few $F_a$'s. This union is produced by performing 2-surgery on $F$. 
Hence,  $\chi_-(F) \geq \sum_a \chi_-(F_a)$. Therefore, the surface 
$F_{best}^{[\Sigma]} := \coprod_{a\in A} \;\kappa^{\mathcal W_\Sigma}_a\cdot F_a$ also delivers $\|[\Sigma]\|_{H^f}$ ---
the minimum of the $\chi_-$-characteristic among \emph{all} combinations of fiber components which realize the 
homology class $[\Sigma]$. Here $\kappa^{\mathcal W_\Sigma}_a\cdot F_a$ stands for 
a disjoint union of $|\kappa^{\mathcal W_\Sigma}_a|$ 
fiber components residing in a regular neighborhood of the $f$-oriented surface $F_a$; their orientations are prescribed 
by the sign of $\kappa^{\mathcal W_\Sigma}_a$.

\begin{thm} Let $f: M \rightarrow S^1$ be a Morse map with no local extrema, no bubbling 
singularities\footnote{For example, when no homology class in $H_2^f$ admits a spherical or disk 
representative, any intrinsically harmonic $f$ will do.} and with $f: \partial M \rightarrow S^1$ 
being a fibration.  Then, for any oriented surface $\Sigma$ 
representing a vertical 2-homology class $[\Sigma] = \sum_a \alpha_a [F_a]$ 
%and well-positioned with respect 
%to $F_R$\footnote{cf. Definition 7.1},
\begin{equation}
\chi_-(\Sigma) \; \geq \;  \sum_{a\in A} |\kappa^{\mathcal V_\Sigma}_a|\cdot \chi_-(F_a) - 
\sum_{r\in R} \rho^\circ(\Sigma, F_r)\cdot \chi_-(F_r) - \mu^\circ(\Sigma, F_R).
\end{equation}
Here $\kappa^{\mathcal V_\Sigma} \in \Z[A]$ is the  vector of the affine lattice 
(8.3)\footnote{For a harmonic $f$, one can pick a basis of $f$-positive loops $\{C_k\}$ in the formulas 
(8.3), (8.4).}  which  minimizes the norm (8.1). Its norm depends only on the homology class $[\Sigma]$ 
and the twist function $\rho_\Sigma : R \rightarrow \Z_+$. The number $\mu^\circ(\Sigma, F_R)$, 
counting the new spherical and disk components in the resolution $\Sigma^\odot \bowtie F_R$, is defined 
by the formula (7.1) (with $F = F_R$).

Formula (8.5) can be expressed in terms of the vertical norms:
\begin{equation}
\chi_-(\Sigma) \; \geq \;  \|[\Sigma] + \sum_{r \in R} \rho^\circ(\Sigma, F_r) [F_r]\|_{H^f} - 
\|\sum_{r\in R} \rho^\circ(\Sigma, F_r)\cdot F_r\| - \mu^\circ(\Sigma, F_R). 
\end{equation}
 \end{thm}

{\bf Proof.}\quad  As in Section 6, we consider transversal intersections $\Sigma \cap F_r$ 
giving rise to the twist invariants $\rho_r(\Sigma) := \rho(\Sigma,\; F_r)$.  As before, 
special attention is paid to the curves from $\Sigma \cap F_r$ which 
bound disks in $F_r$. They  help us to perform 2-surgery on $\Sigma$, which only can diminish 
the value of $\chi_-(\sim)$ and $g(\sim)$. The surgery produces a surface $\Sigma^\odot$.

Note that the resolutions along distinct $\Sigma^\odot \cap F_r$'s are completely independent and can be 
performed in any order.
It will require $\rho^\circ_r(\Sigma) := \rho(\Sigma^\odot,\; F_r) = \rho^\circ(\Sigma,\; F_r)$ resolutions 
to separate $\Sigma$ and $F_r$. 
Therefore, after at most $\rho^\circ_R(\Sigma) := max_{r\in R} \{\rho^\circ_r(\Sigma) \}$ 
resolutions, the resolved surface 
$\Sigma^\odot \bowtie  F_R$ will be separated from the fiber union
$F_R := \cup_r F_r$. It will reside in the homology class $[\Sigma] + \sum_r \rho^\circ_r(\Sigma) [F_r]$.
\smallskip

By Lemma 7.3,
\begin{equation}
\chi_-(\Sigma^\odot \bowtie F_R) \leq 
\chi_-(\Sigma) + \sum_{r\in R} \; \rho^\circ_r(\Sigma)\cdot\chi_-(F_r) + \mu^\circ(\Sigma, F_R).
\end{equation}
Now, as in the proof of Theorem 5.2, using the gradient and minus the 
gradient flows, we can push $\Sigma^\odot \bowtie F_R$ away 
from  $F_R$ and towards the union  $F_A := \coprod_{a\in A} F_a$ (here we relay on Lemma 4.2).
This push becomes possible only after a number of 2-surgeries using 
the descending and the ascending 2-disks of critical points 
as  their cores. The 2-disks are facing the fibers
$F_r$ (see Figure 9 and a similar argument in the proof 
of Theorem 5.2).  The surgery will produce a new surface 
$\Sigma^\star$, homologous to $\Sigma^\odot \bowtie F_R$. By Lemma 3.1,
$\chi_-(\Sigma^\star) \leq \chi_-(\Sigma^\odot \bowtie F_R)$. Therefore, 
$\Sigma^\star$ resides in a regular neighborhood of $F_A$ and its 
$\chi_-$-characteristic is less than or equal to the RHS of (8.7)  
Denote by $\Sigma^\star_a$ the union of $\Sigma^\star$-components $\Sigma^\star_{a, i}$ residing in 
a regular neighborhood $U_a$  of $F_a$. The retraction $U_a \rightarrow F_a$ induces  maps 
$\Sigma^\star_{a,i} \rightarrow F_a$ of  degrees $\kappa_{a,i}$. Put $\kappa_a = \sum_i \kappa_{a,i}$.
By an argument as in Lemma 5.3, $\chi_-(\Sigma^\star_a) \geq \sum_a |\kappa_a|\cdot \chi_-(F_a)$. 
Therefore, we get
%%%%%
\begin{equation}
\chi_-(\Sigma) + \sum_{r\in R} \; \rho^\circ_r(\Sigma)\cdot\chi_-(F_r)  + \mu^\circ(\Sigma, F_R) \geq 
\sum_{a\in A} |\kappa_a|\cdot \chi_-(F_a)
\end{equation}
%%%%%%

On the other hand, since $\Sigma^\star$ and  $\Sigma^\odot \bowtie F_R$ are cobordant,
%%%% 
\begin{equation}
\sum_{a\in A} \kappa_a \cdot [F_a] = [\Sigma] + \sum_{r\in R} \rho^\circ_r(\Sigma)\cdot[F_r].
\end{equation}
%%%%%
This equation implies that $\kappa \in \Z[A]$ defined by $\{\kappa(a) = \kappa_a\}$ belongs to the affine 
sublattice $\mathcal V_\Sigma$ (cf. (8.3)) Indeed, just consider the intersection numbers of the  
basic 1-cycles $\{C_k\}$ with the LHS and RHS of (8.9). 

Since we have little control over the $\kappa$, we safely minimize the RHS of (8.8) ---
the $\|\sim \|_A$-norm of $\kappa$ --- over the set $\mathcal V_\Sigma$ to get the 
desired inequality (8.5).  \qed
\smallskip

\smallskip

When all the twists $\{\rho^\circ_r(\Sigma) = 0\}$, the spaces (8.3) and (8.4) coincide. Furthermore,
for such a $\Sigma$,  $\mu^\circ(\Sigma, F_R) = 0$. Therefore, for 
any $\Sigma$ in the homology class $[\Sigma]$, $\chi_-(\Sigma) \; \geq \; \|[\Sigma]\|_{H^f}$. 
On the other hand, by definition, $\|[\Sigma]\|_T \leq \|[\Sigma]\|_{H^f}$. Hence, the vanishing 
$\{\rho^\circ_r(\Sigma) = 0\}$ implies $\|[\Sigma]\|_T = \|[\Sigma]\|_{H^f}$---a union of 
fiber components delivers $\|[\Sigma]\|_T$. \smallskip  

Let $\Sigma \subset M$ be the best combination of $f$-fiber components which delivers the Thurston norm of $[\Sigma]$.
If we deform $f$ in such a way that the new map $f_1$ and the old $f$ share the same set of the "repelling" 
fiber components, then vanishing of $\{\rho^\circ_r(\Sigma)\}$ for $f$ implies vanishing of 
$\{\rho^\circ_r(\Sigma)\}$ for $f_1$.
 
In combination with Corollary 6.12, these observations leads to
%%% 
\begin{cor} {\bf (The Best Fiber Component criterion)} \; 
Let $f: M \rightarrow S^1$ be a Morse map as in Theorem 8.2.  Then, for any oriented surface $\Sigma$ 
representing a vertical 2-homology class $[\Sigma] = \sum_a \alpha_a [F_a]$

\begin{itemize} 
\item If $\rho^\circ(\Sigma, F_R) = 0$, then $\Sigma$ delivers the Thurston norm $\|[\Sigma]\|_T$. 

\item For such maps $f$, the Best Fiber Component Theorem 
$\{\|[\Sigma]\|_T = \|[\Sigma]\|_{H^f}\}$ is equivalent to the property  
$\rho_{\chi_-}([\Sigma], F_R) = 0$.  

\item Let $f_1$ be a map  as in Theorem 8.2 and homotopic to $f$. Assume 
that the two maps share the same set of fiber components indexed by their repeller sets $R$ and $R_1$.   
Then $\{\|[\Sigma]\|_T = \|[\Sigma]\|_{H^f}\}$ implies $\{\|[\Sigma]\|_T = \|[\Sigma]\|_{H^{f_1}}\}$. 

\item When $[\Sigma] = [F]$---the homology class of a fiber---, then 
$b_{\chi_-}(F_R, [F]) = 0$ 
implies  $\|[F]\|_T = \|[F]\|_{H^f}$.\qed 
\end{itemize}
\end{cor}
%%%%%
\begin{example}
\end{example}
%%%%%
Often the best union of fiber components and the best fiber are quite different. 
The relation between them could be non-trivial, but it can  be described in pure 
combinatorial terms involving $\pi_f: \Gamma_f \rightarrow S^1$,  the 
1-cochain $\tau_{\chi_-}(f)$ from Section 3 and a marking of edges corresponding to the spherical
components.

Let us examine Figure 1. The fiber $F$  corresponding to a ray from the center which intersects with the 
slanted edge of $\Gamma_f$ (in the dark shaded sector) is comprised of two spherical
components and a surface $F_0$ of genus 3. Note that $F_0$  is holmologically trivial in $M$. 
Therefore, the two spheres $F_1 : = F \setminus F_0$ are homologous to $[F]$. Since 
$\chi_-(F_1) = 0$, $F_1$ is the best combination of fiber components. On the other 
hand, $\{\chi_-(f^{-1}(\theta))\}$ take values 4 and 2 only. We notice that the map $f$ 
violates the Calabi positive loop property and, hence, is not intrinsically 
harmonic. \qed
\bigskip  

In general, the relation between  solutions $\kappa^{\mathcal V_\Sigma}$ and $\kappa^{\mathcal W_\Sigma}$ 
of the two optimization problems (cf. (8.3) and (8.4)) is subtle. However, for special, so 
called, $f$-\emph{balanced} $[\Sigma]$'s, we can get a handle on  the relation between the optimal 
norms $\|\kappa^{\mathcal V_\Sigma}\|_A$ and $\|\kappa^{\mathcal W_\Sigma}\|_A$.

\begin{defn} We say that a homology class $[\Sigma] \in H_2(M, \partial M; \Z)$ is $f$-\emph{balanced}, 
if it is proportional, over the positive rationals, to  the vertical class $[F_R]$.
\end{defn}
In combinatorial terms,  the proportion between the numbers of $(\pm)$-weighted repellers and the 
$[\Sigma]$-supporting attractors  along any loop $C$ in $\Gamma_f$ is positive and $C$-independent. The sign
attached to each singleton depends on the orientation of the loop and the orientation of the singleton 
induced by the map $\pi_f: \Gamma_f \rightarrow S^1$.

For example, when all the $f$-fibers are connected, the homology class of a fiber or its multiples are 
$f$-balanced. Also, the class $[F_A] = [F_R]$ and hence, is $f$-balanced.\bigskip  

We notice that the inequalities (8.7) can be relaxed by replacing each twist $\rho^\circ_r(\Sigma)$ by 
their maximum $\rho^\circ(\Sigma, F_R) := max_{r\in R} \{\rho^\circ_r(\Sigma)\}$. In other words, one can 
employ the same number $\rho^\circ(\Sigma, F_R)$ of resolutions at each $F_r$ to create a "less optimal" 
surface $\Sigma^\odot \bowtie F_R$ which might contain a few extra-copies of some $F_r$'s.  
Hence, (8.8) and (8.9) will be modified: 
\begin{equation}
\chi_-(\Sigma) + \rho^\circ(\Sigma, F_R)\cdot\sum_{r\in R} \; \chi_-(F_r) + \mu^\circ(\Sigma, F_R) \geq 
\sum_{a\in A} |\kappa_a|\cdot \chi_-(F_a)
\end{equation}
\begin{equation}
\sum_{a\in A} \kappa_a \cdot [F_a] = [\Sigma] + \rho^\circ(\Sigma, F_R)\cdot\sum_{r\in R} \;[F_r].
\end{equation} 

Employing (8.11), we define affine sublattices in $\Z[A]$ --- modified versions of (8.3), (8.4), --- by prescribing 
the intersection numbers of the 2-cycle $\sum_{a\in A} \kappa_a \cdot [F_a]$ with a basis of loops $\{C_k \subset M\}_k$ 
in $H_1(\Gamma_f, \Z)$:
\begin{eqnarray}
\tilde{\mathcal V}_\Sigma &:= &\bigg \{ \kappa\in \Z[A]:\quad \int_{C_k} \kappa = \; ([\Sigma]\circ C_k) + 
\rho^\circ(\Sigma, F_R)\cdot([F_R]\circ C_k) \bigg \}\\ 
\;\;\quad \mathcal W_\Sigma &:= & \bigg \{ \kappa\in \Z[A]:\quad \int_{C_k} \kappa = \; ([\Sigma]\circ C_k)  \bigg\}\\
\mathcal U &:= & \bigg \{ \kappa\in \Z[A]:\quad \int_{C_k} \kappa = \; ([F_R]\circ C_k) \bigg \}
\end{eqnarray}

As before, we are interested in vectors 
$\kappa^{\tilde{\mathcal V}_\Sigma} \in \tilde{\mathcal  V}_\Sigma,\; \kappa^{\mathcal W_\Sigma} \in \mathcal W_\Sigma,\;
 \kappa^{\mathcal U} \in \mathcal U$ 
which will minimize the $\|\sim\|_A$-norm.
When $[\Sigma]$ is $f$-balanced (i.e. proportional to $[F_R]$), the vectors $\{[\Sigma]\circ C_k\}_k$ and 
$\{[F_R]\circ C_k\}_k$ are \emph{proportional} with a positive coefficient of proportionality. As a result, we can assume that 
$\kappa^{\tilde{\mathcal  V}_\Sigma} = \kappa^{\mathcal W_\Sigma} + \rho^\circ(\Sigma, F_R)\cdot\kappa^{\mathcal U}$, where 
$\kappa^{\mathcal W_\Sigma}$ and $\kappa^{\mathcal U}$ being proportional with the positive proportionality coefficient. By definition,
$\|\kappa^{\mathcal W_\Sigma}\|_A = \|[\Sigma]\|_{H^f}$. We notice that, 
 $\|\kappa^{\mathcal U}\|_A = \|[F_A]\|_{H^f}$. Indeed, use Lemma 4.3 to replace  $C_k \circ F_R$ with 
$C_k \circ F_A$  in (8.12) and (8.14).

Minimizing the RHS of (8.10), subject to (8.11), and using that  
$$\|\kappa^{\tilde{\mathcal V}_\Sigma}\|_A = \|\kappa^{\mathcal W_\Sigma}\|_A + 
\rho^\circ(\Sigma, F_R)\cdot\|\kappa^{\mathcal U}\|_A 
= \|[\Sigma]\|_{H^f} + \rho^\circ(\Sigma, F_R)\cdot\|[F_A]\|_{H^f}$$
$$=  \|[\Sigma] + \rho^\circ(\Sigma, F_R)\cdot [F_A]\|_{H^f},$$
we get our main result:

%%%%%
\begin{thm} Let $f: M \rightarrow S^1$ be a Morse map with $f: \partial M \rightarrow S^1$ being a fibration. 
Assume that $f$ has no local extrema and no bubbling singularities.   Let $[\Sigma] \in H_2(M, \partial M; \Z)$  be an 
$f$-balanced  class. Put $\|F_R\| := \chi_-(F_R)$.  

Then, for any oriented surface 
$(\Sigma, \partial\Sigma) \subset (M, \partial M)$ representing $[\Sigma]$,
\begin{equation}
\chi_-(\Sigma) \;\geq\; \|[\Sigma]\|_{H^f} - \rho^\circ(\Sigma, F_R)\cdot Var_{\chi_-}(f) - \mu^\circ(\Sigma, F_R),
\end{equation} 
where $Var_{\chi_-}(f) := \|F_R\|  - \|[F_A]\|_{H^f} = \|F_R\|  - \|[F_R]\|_{H^f}$.
\smallskip

When $\Sigma$ is well-positioned with respect to $F_R$, then $\mu^\circ(\Sigma, F_R) = 0$. 

For a self-indexing map $f$, the variation $\|F_R\|  - \|[F_R]\|_{H^f}$ is equal to the number of (non-bubbling) 
$f$-singularities.
\qed 
\end{thm}
%%%%%%%% 
\begin{cor}The statements of Theorem 8.6 are valid when $f$ is intrinsically harmonic and no non-trivial class in 
$H_2^f$ admits a representation by  spheres and disks. For instance, this is the case, if the Hurewicz homomorphism 
$\pi_2(M, \partial M) \rightarrow H_2(M, \partial M; \Z)$ is trivial. In particular, Theorem 8.6 is valid for 
harmonic Morse maps of irreducible 3-manifolds. \qed
\end{cor}
%%%%%%

Let $\mu^\circ([\Sigma], F_R)$ denote the minimum  of (7.1) taken over all surfaces $\Sigma$ which deliver the 
Thurston norm of $[\Sigma]$ and the value $\rho^\circ([\Sigma], F_R)$.\smallskip

Minimizing the RHS of (8.15) over the set of all surfaces which deliver the Thurston norm of $[\Sigma]$ 
and employing Definition 6.3, we get  
%%%%% 
\begin{cor} Let $f: M \rightarrow S^1$ be as in Theorem 8.6 and $[\Sigma] \in H_2^f$ be a balanced class.
Then its Thurston semi-norm $\|[\Sigma]\|_T$ can be compared to the $f$-vertical semi-norm 
$\|[\Sigma]\|_{H^f}$ on $H_2^f$:
\begin{equation}
\|[\Sigma]\|_{H^f} \; \geq \; \|[\Sigma]\|_T \; \geq \; \|[\Sigma]\|_{H^f}  - 
\rho_{\chi_-}([\Sigma], F_R)\cdot Var_{\chi_-}(f) - \mu^\circ([\Sigma], F_R)
\end{equation}
In particular, if the twist $\rho_{\chi_-}([\Sigma], F_R)$ or the variation $Var_{\chi_-}(f)$ vanish, then 
$\mu^\circ([\Sigma], F_R) = 0$, and we get $\|[\Sigma]\|_T = \|[\Sigma]\|_{H^f}.$ \qed
\end{cor}
%%%%%%%%%
Vanishing of the variation $Var_{\chi_-}(f) \geq \chi_-(F_R) - \chi_-(F_A)$ is a rare event: 
it can only happen when 
$\chi_-(F_R)  = \chi_-(F_A)$. In such a case,  By Lemma 3.5, all the singularities of $f$ must be of 
the bubbling type. If all the singularities are bubbling, then all the fiber 
components are incompressible. Indeed, let $S^2_i$ be a spherical fiber component residing in the 
vicinity of  a critical point $x_i$. If we cut $M$ open along $\sqcup_i S^2_i$, it will decompose 
into a number of spherical rings $S^2_i \times [0,1]$ and the rest, which we denote by $M^\odot$. 
We can attach a disk $D^3_i$ to each $S^2_i \subset M^\odot$ and  extend $f$ in an obvious 
way across the disk. For each $i$, the new map $f': M' \rightarrow S^1$ will have exactly one new critical 
point $y_i$---a local extremum---located at the center of $D^3_i$. A standard deformation $f''$ of $f'$ will 
cancel each pair $(x_i, y_i)$, thus, producing a fibration $f''$.  Its fibers are incompressible. 
Now, if a disk $D^2 \subset M$ bounds a loop $\gamma$ residing in an $f$-fiber $F \subset M^\odot$, then there is 
disk $D^2 \subset M'$ which also bounds
$\gamma$. Since fibers of $f''$  are incompressible, $\gamma$ must bound a disk in $F$. 
Therefore, if $\chi_-(F_R) = \chi_-(F_A)$, then for any $\Sigma$,
$\mu^\circ(\Sigma, F_R) = 0$.\smallskip

Although the hypotheses of Theorem 8.6
exclude the bubbling singularities, its generalization---Theorem 8.13,---allows them.   
A model example can be produced by attaching 1-handles to a disjoint union
of  fibrations over the circle and extending the fibering maps across the handles as depicted in Figure 16.

%%%%%%%
\begin{cor} Let $f: M \rightarrow S^1$  and $\Sigma$ be as in Corollary 8.8. Assume that  
a surface $\Sigma$ which delivers the Thurston norm $\|[\Sigma]\|_T$ together with the twist 
$\rho_{\chi_-}([\Sigma], F_R)$ is well-positioned with respect to $F_R$.
If $Var_{\chi_-}(f) \neq 0$, then
\begin{equation}
\rho_{\chi_-}([\Sigma], F_R) \;\geq \; \frac{\|[\Sigma]\|_{H^f} \; - \; \|[\Sigma]\|_T}{\; \|F_R\|
 \;\; - \;\;\|[F_R]\|_{H^f}}
\end{equation} 
When $\Sigma$ represents the \emph{balanced} homology class $[F]$ of a fiber, then the RHS of (8.17) 
also gives a lower bound 
on the size of the breadth $b_{\chi_-}(f)$ and height $h_{\chi_-}(f)$. \smallskip

Therefore, if the variation $\|F_R\|  - \|[F_R]\|_{H^f}$ is relatively small and 
$\|[F]\|_{H^f} \gg \|[F]\|_T$, then any well-positioned surface $\Sigma$ which delivers the Thurston norm 
of the homology class $[F]$ must be $\tilde f$-tall.\smallskip
\qed
\end{cor}

Now, we consider maps $f$ whose graphs $\Gamma_f$ are very special. For them, the homology class of a 
fiber is balanced. The self-indexing maps evidently fall in that  category.

Theorem 8.6 and Corollary 8.8, in combination with Proposition 6.8 and Lemma 6.10, imply a relaxing 
progression of  inequalities:
%%%%%%
\begin{cor} Let $f: M \rightarrow S^1$ be as in Corollary 8.8. Assume that the homology 
class $[F]$ of a fiber is balanced. 
Let $\Sigma \subset M$ represent  $k[F]$ and is well-positioned with respect to $F_R$. 
Then (8.15) implies 
\begin{eqnarray}
\chi_-(\Sigma) &\;\geq\;& \|[F]\|_{H^f} - b(F_R, \Sigma)\cdot Var_{\chi_-}(f)\\
\chi_-(\Sigma) &\;\geq\;& \|[F]\|_{H^f} - h(\Sigma; f)\cdot Var_{\chi_-}(f)
\end{eqnarray}
Similarly, (8.16) implies 
\begin{eqnarray}
\|[F]\|_T \;& \geq \;& \|[F]\|_{H^f}  - b_{\chi_-}(F_R,\, k[F])\cdot Var_{\chi_-}(f)\\
\|[F]\|_T \;& \geq \;& \|[F]\|_{H^f}  - h_{\chi_-}(k[F],\, f)\cdot Var_{\chi_-}(f) \qquad \qed
\end{eqnarray}
\end{cor}

The corollary below and Figure 14 depict a special case of the  balanced fiber.

\begin{cor} Let $f: M \rightarrow S^1$ be a Morse map with no local extrema and 
whose restriction to the boundary $\partial M$ is
a fibration.  Assume that all the fiber components
$\{F_r\}_{r\in R}$,  actually, are \emph{fibers} (cf. Figure 14). Let $\Sigma$  be an oriented surface which 
represents the homology class $[F]$ of a fiber and is well-positioned with respect to $F_R$. 
Let $F_{best}$ denote an oriented  union of fiber components which delivers $\|[F]\|_{H^f}$. Then 
\begin{equation}
\chi_-(\Sigma)\; \geq \; \chi_-(F_{best}) - \rho^\circ(\Sigma, F_R) \cdot 
\sum_{r\in R} \big [\chi_-(F_r) - \chi_-(F_{best})\big]
\end{equation}
As a result, 
\begin{equation}
\|[F]\|_T\; \geq \; \|[F]\|_{H^f} - \rho_{\chi_-}([F], F_R) \cdot 
\sum_{r\in R} \Big [\chi_-(F_r) - \|[F]\|_{H^f}\Big ]
\end{equation}
In particular, (8.22) and (8.23) are valid for any map $f$ with all its fibers being connected, 
provided that a well-positioned $\Sigma$ delivering $\|[F]\|_T$ and $\rho_{\chi_-}([F], F_R)$ exists.  
\qed
\end{cor}
%%%%%
\begin{figure}[ht]
\centerline{\BoxedEPSF{fibers.repellers scaled 450}} 
\bigskip
\caption{}
\end{figure} 
%%%%%%%
\begin{example}
\end{example}
%%%%%
We would like to recycle the algorithm from
the Harmonic Twister example. We will see that this algorithm  
produces maps $f$ with few singularities, a small variation, but arbitrary big twists
and heights. 

Let us start with the fibration $f_0: T \rightarrow S^1$ of the solid torus $T = D^2 \times S^1$ 
over the circle. The $\chi_-$-number of its fiber $D^2$ is zero. 
By deforming $f_0$ slightly inside $T$, we can introduce a pair $\{a, b\}$ of index 1 and index 2 
critical points, while keeping all the fibers \emph{connected}: 
for a homological reason, the pattern depicted in the left hand side of the diagram $E$ of 
Figure 5 cannot be realized by a trivial cobordism. 

%%%%%
\begin{figure}[ht]
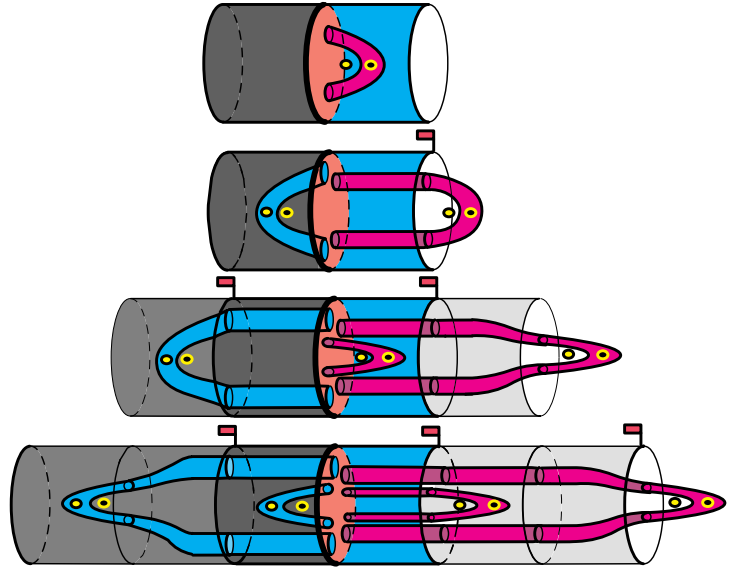

\centerline{\BoxedEPSF{twister scaled 550}} 
\bigskip
\caption{A lift of the worst $f'_k$-fiber to the space of the cyclic cover 
$D^2 \times \R \rightarrow D^2 \times S^1$. The product structure is 
defined by $f_0$.}
\end{figure} 
%%%%%%%

Now, as in the Harmonic Twister 
example, we deform $f'_0$ into a new map $f'_k$ while keeping the deformation fixed on 
the boundary $\partial T$. This is done by applying $k$ times 
the move $A$ from Figure 5. In the process,  $b$ overcomes 
$a$ exactly $k$ times in the race around the circle.  By Lemma 6.14, 
$h_{\chi_-}(f'_k) \leq 1 + k$. 
As Figure 5, $A$, testifies, all the fibers of $f'_k$ still must be connected. Also, 
because of the same diagram $A$, $\chi_-(F'_{k, best}) = 2k$. 
Evidently, $Var_{\chi_-}(f'_k) = 2$. 

In fact, the original fiber $D^2$ is well-positioned 
with respect to the surface $F_{r, k}$ --- the worst fiber of $f'_k$. To validate this fact requires 
a more careful analysis of the geometry of $F_{r, k}$ relative to $f_0$, as depicted in 
Figure 15 for $k = 0, 1, 2, 3$. The figure shows a lift $\hat F_{r, k}$ of  
$F_{r, k}$ to the space $D^2 \times \R$ of the cyclic cover. The surface $\hat F_{r, k}$ 
is comprised of a number of left and right 1-handes attached to a disk marked with a bold circle. 
Deforming $f'_{k - 1}$ into $f'_k$ results in attaching a new right ($k \equiv 0 (2)$) or left handle 
($k \equiv 1 (2)$) to $\hat F_{r, k - 1}$ and stretching the old handles.
The dots with dark (light) centers show the locations of the $\tilde f'_k$-critical points of index 1 (2).
 
As we resolve intersections of $\hat F_{r, k}$ with multiple translates of $\{t^n(\hat D^2)\}$ of the 
$\tilde f_0$-fiber $\hat D^2$ (they are marked with flags in Figure 15), we see 
that each translate cuts through $\hat F_{r, k}$ in a way
that  leaves at least one handle of $\hat F_{r, k}$ to the left and one handle to the right of $t^n(\hat D^2)$. 
Hence, no new 2-disks or spheres are produced as a result of the resolutions, i.e. $\mu^\circ(D^2, F_{r, k}) = 0$.   
By Theorem 8.6, $h_{\chi_-}(f'_k) \geq \rho_{\chi_-}(f'_k) \geq k$. Furthermore, Figure 15 
testifies that $\rho_{\chi_-}(f'_k) = k$.

Since the deformed maps $f'_k$ coincide with $f_0$ on the boundary $\partial T$, we can use them 
for twisting any given map $f: M \rightarrow S^1$:  just take a regular neighborhood of an $f$-positive 
loop $\gamma \subset M$ for the role of $T$.
\qed
\bigskip 

Let maps $\{f_\alpha : M_\alpha \rightarrow S^1\}_\alpha$ be as in Theorem 8.2. As we attach 
1-handles $\{T_\beta \approx S^2_\beta\times D^1_\beta\}_\beta$ to $\coprod_\alpha M_\alpha$ 
and extend the maps across the handles as shown in
Figure 16, we form a new manifold $M$ and a new map $f : M \rightarrow S^1$. Its graph $\Gamma_f$ is obtained
from $\coprod_\alpha \Gamma_{f_\alpha}$ by attaching to it a few new edges, each of which contributes 
a new pair of trivalent vertices of indices 1 and 2.  Moreover, the orientation of these new edges 
is such that each of them must contain an attractor. In other words, 1-surgery does not 
change the repeller set! New branches will be added to the original trees $\{T_r^{\pm}\}$ covering 
$\coprod_\alpha \Gamma_{f_\alpha}$, and the new trees with the old roots will cover the new 
graph $\Gamma_f$ (as in Lemma 4.2). \smallskip
%%%
\begin{figure}[ht]
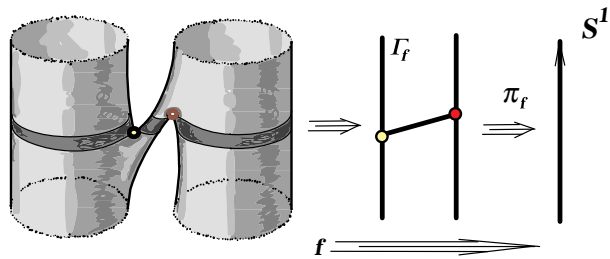

\centerline{\BoxedEPSF{handle scaled 450}} 
\bigskip
\caption{Attaching slanted 1-handle}
\end{figure}
%%%%%

From our point of view, the "simplest" case of fibrations $f_\alpha$ is a bit subtle: formally, fibrations do not 
satisfy the basic assumptions of Lemma 4.2. Recall that our techniques rely on the decomposition of
$\Gamma_{f_\alpha}$  into the trees $\{T^{\pm}_{\alpha, r}\}$. Attaching handles to the space of a fibration 
might produce a loop in the new graph only with bubbling vertices of the same index. In such a case, we 
need first to deform slightly the fibering map $f_\alpha$, so that a pair of canceling critical points 
of indices 1 and 2 is introduced. We can choose a deformation which will not disturb the majority  of
fibers. So, the best fiber of the deformed map remains the winner in its homology class. Effectively,
the deformation introduces one repeller along the loop which represents the graph $\Gamma_{f_\alpha}$.
After such conditioning of $f_\alpha$, we are ready to attach 1-handles any way we want.

Consider the natural epimorphism
\begin{eqnarray} 
\Phi_2: H_2(M, \partial M;\, \Z) \stackrel{j_2}{\longrightarrow} H_2(M, \partial M \sqcup(\sqcup_\beta S^2_\beta);\, \Z) \\
\approx \oplus_\alpha\; H_2(M_\alpha, \partial M_\alpha \sqcup(\sqcup_{\beta\in B_\alpha} D^3_{\alpha\beta});\, \Z)
\approx \oplus_\alpha\; H_2(M_\alpha, \partial M_\alpha; \Z), \nonumber
\end{eqnarray}
where $\{D^3_{\alpha\beta}\}_{\beta\in B_\alpha}$ denote 3-disks in $M_\alpha$ bounding  the bases of 
the 1-handles attached to $M_\alpha$. For every class $[\Sigma] \in H_2(M, \partial M;\, \Z)$, let  
$\Phi_2([\Sigma]) = \oplus_\alpha [\Sigma]_\alpha$. The epimorphism $\Phi_2$ splits. We pick 
a splitting homomorphism $\Psi_2: \oplus_\alpha\; H_2(M_\alpha, \partial M_\alpha; \Z) \rightarrow 
H_2(M, \partial M;\, \Z)$. Geometrically, $\Psi_2$ can be defined by picking a basis of 2-cycles in \hfil\break 
$\oplus_\alpha\; H_2(M_\alpha, \partial M_\alpha; \Z)$ and isotoping 
them away from the 3-disks $\{D^3_{\alpha\beta}\}$.   There is nothing canonical about our
choice of $\Psi_2$. 

Note that $\Phi_2([S^2_\beta]) = 0$. Therefore, vertical classes are mapped by $\Phi_2$ to vertical classes: 
$\Phi_2(H_2^f) = \oplus_\alpha\, H_2^{f_\alpha}$. \smallskip 

Our goal is to describe the relation between the norms $\|[\Sigma]\|_T$ and $\{\|[\Sigma]_\alpha\|_T\}$.
In the process, we also will describe the relation between the twists $\rho_{\chi_-}([\Sigma], F_R)$ and 
$\{\rho_{\chi_-}([\Sigma]_\alpha, F_{R_\alpha})\}$  of the maps $f$ and $\{f_\alpha\}$.\bigskip

Let $\{\Sigma_\alpha \subset M_\alpha\}$ be surfaces which deliver the norms $\|[\Sigma_\alpha]\|_T$ of their 
$f_\alpha$-vertical homology classes $[\Sigma_\alpha]$ and  the 
values of the twist invariants $\rho([\Sigma_\alpha], F_{R_\alpha})$ for the maps $\{f_\alpha\}$. 
Let $[\Sigma]$ be a vertical class, so that $\Phi_2([\Sigma]) = \oplus_\alpha [\Sigma_\alpha]$. 
Consider the surface $\sqcup_\alpha \Sigma_\alpha \subset M$. Since the kernel of $\Phi_2$ is 
spanned by the spheres $\{S_\beta^2\}$, for some integers $\{\kappa_\beta\}$, we have
$[\Sigma] = \sum_\beta \kappa_\beta [S^2_\beta] + \sum_\alpha \Psi_2([\Sigma_\alpha])$. Therefore,
$[\Sigma]$ admits a representative $\Sigma$ which is a disjoint union of 
$\sqcup_\alpha \Sigma_\alpha \subset M \setminus (\sqcup_{\alpha, \beta} D^3_{\alpha, \beta})$ with 
a bunch of vertical spheres. Recall, that the repeller set $R$ for $\Gamma_f$ is the disjoint union of 
the repeller sets $\{R_\alpha \subset \Gamma_f\}$.
The spheres do not contribute to $\rho(\Sigma, F_R)$:  they are disjoint 
from the surface $F_R$. Therefore, $\rho_{\chi_-}([\Sigma], F_R) \leq  
max_\alpha\, \rho_{\chi_-}([\Sigma_\alpha], F_{R_\alpha})$.
Also, $\chi_-(\Sigma) = \sum_\alpha \chi_-(\Sigma_\alpha) = \sum_\alpha \|[\Sigma_\alpha]\|_T$. So, 
$\|[\Sigma]\|_T \leq \sum_\alpha \|[\Sigma_\alpha]\|_T$. 

A similar argument is applicable to an estimation of $h_{\chi_-}(f)$.
Consider surfaces $\{\Sigma_\alpha \subset M_\alpha \setminus (\sqcup_\beta D^3_{\alpha, \beta})\}$ 
which deliver the norm $\|[F_\alpha]\|_T$ of the
$f_\alpha$-fiber homology class $[F_\alpha]$. Among them pick $\Sigma_\alpha$ with the minimal value of 
$h(\Sigma_\alpha, f_\alpha)$ (cf. (6.2)). One can align their special liftings 
$\{\hat \Sigma_\alpha \subset \tilde M_\alpha\}$, so that they will be located above $\tilde f_\alpha^{-1}(0)$
and will have a non-empty intersection with $\tilde f_\alpha^{-1}((0, 1])$. The space $\tilde M$ is built 
by performing an equivariant 1-surgery on $\sqcup_\alpha \tilde M_\alpha$. Although 
$\{\hat \Sigma_\alpha \in \mathcal B_1(\Sigma_\alpha)\}$, $\sqcup_\alpha \hat \Sigma_\alpha$ 
is not automatically in $\mathcal B_1(\sqcup_\alpha  \Sigma_\alpha)$ in $\tilde M$: one needs to add to 
$\sqcup_\alpha  \hat\Sigma_\alpha$ a few spheres $\{\hat S^2_\beta\}$ from 
$\tilde f^{-1}(1) \cap (\sqcup_\beta \tilde T_\beta)$ to get a surface
$(\sqcup_\alpha  \Sigma_\alpha) \sqcup (\sqcup_\beta S^2_\beta)$ which admits 
a special lifting bounding a 3-chain in $\tilde M$. However, this addition will not 
change the value of $\chi_-(\sqcup_\alpha  \Sigma_\alpha)$. Therefore, 
$h_{\chi_-}(f) \leq  max_\alpha\, \{ h_{\chi_-}(f_\alpha)\}$.
\bigskip 

These considerations lead to 
%%%%%
\begin{thm} Let maps $\{f_\alpha : M_\alpha \rightarrow S^1\}_\alpha$ be as in Theorem 8.2.  
Let $f: M \rightarrow S^1$ be constructed by performing 1-surgery on the map $\sqcup_\alpha f_\alpha$.
Then, for any $f$-vertical class $[\Sigma]$,
\begin{eqnarray} 
\|[\Sigma]\|_T \leq \sum_\alpha \|[\Sigma_\alpha]\|_T,
\end{eqnarray} 
where $\Phi_2([\Sigma]) = \oplus_\alpha [\Sigma]_\alpha$. 
Therefore, 1-surgery on a map which satisfies hypotheses of Theorem 8.2 does not 
increase the Thurston norm of its vertical classes (under any splitting homomorphism
$\Psi_2$). 

Also, this 1-surgery does not increase the twist and height invariants:
\begin{eqnarray} 
\rho_{\chi_-}([\Sigma], F_R) &\;\leq\;&  max_\alpha\, \{\rho_{\chi_-}([\Sigma_\alpha], F_{R_\alpha})\}\\
h_{\chi_-}(f) &\;\leq\;& max_\alpha\, \{h_{\chi_-}(f_\alpha)\}. \qquad \qquad  \qed
\end{eqnarray}
\end{thm}
%%%%%%

Although the new map $f$ has bubbling singularities, the conclusions of Lemma 4.2 are still
valid for its graph $\Gamma_f$: it admits a cover by trees rooted at $R = \sqcup_\alpha R_\alpha$ --- 
a fact central to our previous arguments. Indeed, each new edge contains an attractor.
Therefore, 

\begin{cor} If, for each $f_\alpha$ as in Theorem 8.2, a  union 
$F_\alpha \subset M_\alpha \setminus (\sqcup_\beta D^3_{\alpha, \beta})$ of fiber components
delivers the  Thurston norm, that is, if $\| [F_\alpha]\|_{H^{f_\alpha}} = \| [F_\alpha]\|_T$, then 
$\sqcup_\alpha F_\alpha \subset M$  
delivers the Thurston norm of any class $[\Sigma]$ whose $\Phi_2$-image is $\oplus_\alpha [F_\alpha]$. 

In particular, performing 1-surgery on a disjoint union of fibrations $\{f_\alpha\}$,  results in a map $f$ 
whose fiber $F$ delivers $\|[F]\|_T$.
\end{cor}
 
{\bf Proof.}\quad By Corollary 8.3, $\| [F_\alpha]\|_{H^{f_\alpha}} = \| [F_\alpha]\|_T$  is equivalent 
to $\rho([F_\alpha], F_{R_\alpha}) = 0$. Using (8.25) from Theorem 8.13 and again 
Corollary 8.3, the first claim follows. 

In the   case when $f_\alpha$ is a fibration,  first, we deform slightly the fibering map, so that a pair 
of canceling critical points of indices 1 and 2 is introduced and at least one fiber is untouched by the 
deformation.  After this, we  attach 1-handles.   
\qed 
\bigskip 

Our results can be applied to links in 3-manifolds.\smallskip 

We define the Thurston norm $\chi_-(L, M)$ of a (framed) link $L$ in a closed 
3-manifold $M$ to be the minimum of $\chi_-$-invariants of all oriented 
embedded surfaces which bound $L$. Similarly, denote by $g(L, M)$ the genus 
of the link $L$.

In many interesting cases $\chi_-(L, M)$ 
coincides with the Alexander norm $\|(L, M)\|_{Al}$ of $L\subset M$ 
(see [Mc] for the definition) and is 
delivered by the Seifert's algorithm [Cr], [Mur]. In general, $\|(L, M)\|_{Al} \leq \chi_-(L, M)$, 
provided that the first Betti number $b_1(M\setminus L) \geq 2$ [Mc]. 

At least for connected sums of fibered links, we will compute of $\chi_-(L, M)$ in terms of 
the vertical norms.  
\smallskip

Let $L$ be a framed link in an oriented, closed 3-manifold $M$ and $U$ its open tubular 
neighborhood. We denote by $M^{\circ}$ the manifold $M^{\circ} = M \setminus U$.
The framing defines an embedding $L \subset \partial M^{\circ}$.\smallskip

In view of Theorem 8.13 and Corollary 8.14,  we get the following proposition. 

\begin{cor} Let $L_\alpha$ be a framed \emph{fibered} link in an oriented, closed 3-manifold $M_\alpha$. 
Let $F_\alpha \subset \partial M^{\circ}_i$ be an oriented surface which bounds $L_\alpha$ and delivers 
$\chi_-(L_\alpha, M_\alpha)$. Denote by $f_\alpha: M^{\circ}_\alpha \rightarrow S^1$ a 
Morse map which satisfies the hypotheses of Theorem 8.2 and with the surface $F_\alpha$ as one of its fibers
\footnote{An suitable approximation to the Thom-Pontryagin map 
$f_{F_\alpha}: M^{\circ}_\alpha \rightarrow S^1$ will do.}. 

Let $M$ be a manifold constructed by performing 1-surgery on the $\sqcup_\alpha M_\alpha$, the 1-handles being 
attached to $\coprod_\alpha M_\alpha^\circ$. A new link $L = \coprod_\alpha L_\alpha \subset M$ is formed.
The maps $\{f_\alpha\}$ are extended across the 1-handles 
as depicted in Figure 16 to give rise to a new Morse map $f: M^\circ \rightarrow S^1$. \smallskip
Then,  ${\chi_-}(L, M) = \sum_\alpha\, {\chi_-}(F_\alpha).$

In particular, this additivity formula holds when all the $L_\alpha$'s are \emph{fibered} links with 
fibers $F_\alpha$.   \qed
\end{cor}

%%%%
\begin{example}
\end{example}
%%%%%% 
We illustrate this construction in the context of structures 
induced by holomorphic functions $h: \C^2 \rightarrow \C$.
Let $ h({\bf z}) := h(z_1, z_2)$ be a complex polynomial function with isolated critical points. 
Some of these points   $\{{\bf z}_\alpha^\star\}$ might reside on the complex curve $V = h^{-1}(0)$ where 
they manifest themselves as isolated singularities of $V$. 

Consider the intersection $L_{r,\alpha}$  of the $r$-sphere 
$S_{r, \alpha}^3 \subset \C^2$ centered on 
${\bf z}_\alpha^\star$ with the real surface $V$. For a generic $r$, 
$L_{r, \alpha}$ is a smooth link in $S_r^3$. 
 
The function $f({\bf z} ) = h({\bf z} )/ |h({\bf z} )|$ 
defines a smooth map $f: \C^2 \setminus V \rightarrow S^1$. In terms of the natural 
parameter $\theta$ along 
the circle,  $f$ is given by a (real multivalued) formula
%%%%
\begin{eqnarray}
\theta ({\bf z} ) = -i\, log( h({\bf z} )) + i\, log (|h({\bf z} )|) = Im\{ log (h({\bf z}))\}   
\nonumber 
\end{eqnarray}
%%%%
The map $f$ induces maps $f_{r,\alpha}: S_{r,\alpha}^3 \setminus L_{r,\alpha} \rightarrow S^1$. 
According to Milnor [M1], 
for small $r$'s, $f_{r,\alpha}$ are fibrations with a fiber $F_{r,\alpha}$. 
 
In $\C^2 \setminus V$ one can link each pair $(S_{r,\alpha}^3, S_{r,\alpha + 1}^3)$
by a path $\gamma_\alpha$ which is mapped by $f$ in a monotone fashion into $S^1$.
We can pick all $\{\gamma_\alpha\}$ in such a way that they do not intersect 
each other and  do not penetrate in the interior the spheres. 
Attaching to the spheres thin 1-handles $\{D^3\times \gamma_\alpha\}$ with 
the $\gamma_\alpha$'s as cores, we produce a new sphere $S^3$ with a link 
$L = V \cap S^3  = \coprod_\alpha L_{r,\alpha}$ and a map 
$f: (S^{3})^{\circ}  \rightarrow S^1$ of the sort described in Corollary 8.15.
\qed
%%%%%
%%%%%

\section{On tangencies of surfaces to the fibers}

In this section we are going to link the singularities  of a Morse map $f: M \rightarrow S^1$ with 
the singularities of its restriction $f|: \Sigma \rightarrow S^1$ to a probe surface $\Sigma$ 
residing in the complement to the $f$-singularities and  
realizing a particular 2-homology class of $M$.\smallskip

Denote by $\mathcal F_f$ the singular foliation in $M$, defined by the fibers 
of a Morse map $f: M \rightarrow S^1$ as in Theorem 8.2. The map $f: \partial M \rightarrow S^1$ 
is assumed to be a fibration over the circle. 
Denote by $M^\circ$ the complement 
in $M$ to the $f$-critical points, and by $M^\bullet$---the compactification 
of $M^\circ$ by the 2-spheres $\{S_\alpha^2\}$ "surrounding" the singularities 
$\{x_\alpha\}$ of $\mathcal F_f$. Let $\mathcal S^2 := \sqcup_\alpha S^2_\alpha$.  

The 2-plane field, formed by the 
tangent planes to  $\mathcal F_f$, defines an oriented 2-bundle $\xi_f$ over 
$M^\bullet$ (the orientation being induced by the orientations of $M$ and $S^1$). 
Since $f: \partial M \rightarrow S^1$ is a fibration, the bundle 
$\xi_f|_{\partial M}$ is trivial with a preferred trivialization, defined 
by the intersections of the 2-planes with the boundary $\partial M$.  
Its Euler class $\chi(\xi_f) \in H^2(M^\bullet ; \Z)$ is Poincar\'{e}-dual in 
$M^\bullet$ to a 1-cycle $\beta_f \in H_1(M^\bullet, 
\mathcal S^2 ;\, \Z)$. This 1-cycle can be viewed as the zero 
set of a generic section $\sigma$ of the bundle $\xi_f$, subject to a number of 
boundary conditions over $\partial M$ and $\mathcal S^2$.
Specifically, we require the restriction $\sigma|_{\partial M}$ to be a 
non-vanishing section of $\xi_f|_{\partial M}$, pointing in the direction of 
the fibers of $f: \partial M \rightarrow S^1$. Also, we assume 
$\sigma|_{S_\alpha^2}$ to be a generic section of $\xi_f|_{S_\alpha^2}$ 
defined by the trace of $\mathcal F_f$ on $S_\alpha^2$.

For any \emph{immersed}  surface $i: (\Sigma, \partial\Sigma) \rightarrow (M^\circ, \partial M)$, 
the integer $\langle \chi(\xi_f|_{\Sigma}),\, [\Sigma]\rangle$ equals to the 
algebraic intersection $\Sigma\circ\beta_f$. This number is an invariant of 
the homology class of $i_\ast[\Sigma] \in H_2(M^\bullet, \partial M; \Z)$.

Isolated points, where $\Sigma$ is tangent to the non-singular foliation 
$\mathcal F_f$ on $M^\circ$, come in two flavors: \emph{positive}, when the 
preferred orientations of $\Sigma$ and $\mathcal F_f$ agree, and \emph{negative}, 
when they disagree. By intersecting the surface with the fibers, the non-singular 
foliation $\mathcal F$ induces a singular  oriented foliation $\mathcal F_\Sigma$ 
on $\Sigma$. 
Isolated points of tangency also can be of an 
\emph{elliptic} and \emph{hyperbolic} types. For the elliptic type, the index of 
the vector field $X_\Sigma$ on $\Sigma$, normal to $\mathcal F_\Sigma$, is positive and, 
for the hyperbolic ones, it is negative. 

We denote by $h_+$ and $h_-$ the number of 
positive and negative hyperbolic points. Similarly, let $e_+$ and $e_-$ stand for the 
number of positive and negative elliptic tangencies. 

Let $\mathcal I_+$ denote 
the sum of indices of the vector field $X_\Sigma$ at all positive tangent points, 
and let $\mathcal I_-$ denote the sum of indices of  $X_\Sigma$ at all negative 
tangent points. When all the tangencies are of the Morse type, then, 
%%%
\begin{equation} \mathcal I_+ = e_+ - h_+,\quad  \mathcal I_- = e_- - h_-.
\end{equation} 
%%%%%
As in  [T], one can prove that the Euler characteristic of $\Sigma$ and 
its intersection with the relative 1-cycle $\beta_f$ can be calculated in terms 
of the tangencies:
%%%   
\begin{equation}\chi(\Sigma) = \mathcal I_+  + \mathcal I_-,  \qquad
 \langle \chi(\xi_f|_{\Sigma}),\, [\Sigma]\rangle = \mathcal I_+  - \mathcal I_-.
\end{equation} 
%%%%
If $\Sigma$ has no spherical or disk components, then 
$\chi_-(\Sigma) = | \mathcal I_+  + \mathcal I_- |$. 

Consider an immersed surface $(\Sigma, \partial\Sigma) \propto (M^\bullet, \partial M)$, 
which is homologous in $(M^\bullet, \partial M)$ to a combination $F_{best}$  
of fiber components which delivers $\|[\Sigma]\|_{H^f}$. For such a $\Sigma$, 
$\langle \chi(\xi_f|_{\Sigma}),\, [\Sigma]\rangle = 
\langle \chi(\xi_f|_{F_{best}}),\, [F_{best}]\rangle$. 

At the same time, 
$\langle \chi(\xi_f|_{F_{best}}),\, [F_{best}]\rangle = \chi(F_{best})$ --- 
the bundle $\xi_f$ is formed by the tangent planes to the fibers. Therefore,
we have 

\begin{lem}  
If an immersed surface $\Sigma$ is homologous to $F_{best}$  
in $(M^\bullet, \partial M)$ and both $\Sigma$ and $F_{best}$ have no spherical 
and disk components, then 
$\chi_-(\Sigma) = | \mathcal I_+  + \mathcal I_- |$ and 
$\chi_-(F_{best}) = | \mathcal I_+  - \mathcal I_- |$.
As a result, $\chi_-(\Sigma) \geq \chi_-(F_{best})$, if and only if, \;
$\mathcal I_+\cdot \mathcal I_- \geq 0$. \qed
\end{lem} 

In combination with Theorem 8.6, this leads to a theorem below which 
makes it possible to estimate the twist $\rho^\circ(\Sigma, F_R)$ of a probe 
surface $\Sigma$ in terms of the Morse data of $f$ and of its restriction to 
$\Sigma$.

We do not know if any probe surface $\Sigma \subset M$ homologous 
to $F_{best}$ in $M$ can be replaced  by an \emph{embedded} surface 
$\Sigma' \subset M^\bullet$ with $\chi_-(\Sigma') \leq \chi_-(\Sigma)$ and 
homologous to $F_{best}$ in $M^\bullet$. It is easy to verify that such a 
replacement $\Sigma'$ exists among immersed surfaces.

\begin{thm} Let $f: M \rightarrow S^1$ be a  map as in 
Theorem 8.2. Let a surface $\Sigma \subset M^\bullet$ be well-positioned with 
respect to $F_R$ and homologous in $(M^\bullet, \partial M)$ 
to a best combination $F_{best}$ of fiber  
components (i.e. $\chi_-(F_{best}) = \|[F_{best}]\|_{H^f}$). 
Assume that both $\Sigma$ and $F_{best}$ have no spherical 
and disk components. Then
$$Var_{\chi_-}(f)\cdot\rho^\circ(\Sigma, F_R)\; \geq \;  | \mathcal I_+  - \mathcal I_- | - 
| \mathcal I_+  + \mathcal I_- |.$$ 
Furthermore,
$$\frac{1}{2} Var_{\chi_-}(f)\cdot\rho^\circ(\Sigma, F_R)\; \geq \; \mathcal I_- \quad 
 and \quad \mathcal I_+ \leq 0.$$
For generic tangencies, the last two inequalities can be re-written as
\begin{eqnarray}
\frac{1}{2} Var_{\chi_-}(f)\cdot\rho^\circ(\Sigma, F_R)\; \geq \; e_- - h_- \quad 
 and \quad e_+ \;\leq\; h_+.
\end{eqnarray}
\end{thm}

The first statement is non-trivial only when the signs of $\mathcal I_+$ and 
$\mathcal I_-$ are \emph{opposite}. The last two inequalities imply that:
\begin{itemize} 
\item a surface $\Sigma$, as above, with many negative elliptic tangent points and 
few negative hyperbolic ones, must have 
a sizable twist; 
\item
the number of positive
elliptic tangent points does not  exceed the number of positive hyperbolic ones. 
\end{itemize}
\smallskip

\begin{figure}[ht]
\centerline{\BoxedEPSF{ind scaled 380}}
\bigskip
\caption{}
\end{figure}
\bigskip

{\bf Proof}\quad  Employing that $\mathcal I_+  - \mathcal I_- = 
\langle \chi(\xi_f|_{\Sigma}),\, [\Sigma]\rangle = \chi(F_{best}) \leq 0$,
and that $\mathcal I_+  + \mathcal I_- = \chi(\Sigma) \leq 0$, all the statements 
follow from the diagram in Figure 17. The grey area in the diagram  depicts the 
solution set for the inequality 
$|\mathcal I_+  - \mathcal I_- | - \break | \mathcal I_+  + \mathcal I_- | \leq 2a$, 
where $2a = Var_{\chi_-}(f)\cdot\rho^\circ(\Sigma, F_R)$. The solutions of all the three 
inequalities above is the domain shaded with dark grey. The bold diagonal lines 
represent surfaces $\Sigma$ with a fixed negative value of 
$\langle \chi(\xi_f|_{\Sigma}),\, [\Sigma]\rangle$. \qed
\bigskip

We  conclude this chapter with a few remarks 
about minimal surfaces in the homology class of a fiber. We observe that,  the 
singularities of the map $f$ can act as \emph{attractors} for families of minimal 
surfaces. Below we describe such a behavior in general terms. However, our grasp of 
this interesting phenomenon is poor.

The basic fact is that any two minimal 
(connected) surfaces  have only isolated tangencies of 
\emph{hyperbolic} type, unless they are identical. Although, such tangencies are not 
necessarily of the Morse type, they still are canonical: their smooth type 
is modeled after the tangency  at the origin of the surface $\{t = Re(z^n)\}$ 
and the surface $\{t = 0\}$ in the 3-dimensional space $\R\times\C$ 
(cf. Lemma 1.4 in [FHS]) . 

In [K], for a given intrinsically harmonic 1-form $\omega$, and thus for any harmonic map 
into the circle, we have 
constructed (two-parametric) families of Riemannian metrics $g_\mu$ on $M^d$ with the 
following properties:
\begin{itemize}
\item $\omega$ is harmonic with respect to $g_\mu$,
\item outside the disks $\{D_\mu(x_\alpha)\}_\alpha$ of radius $\mu$ surrounding the 
singularities $\{x_\alpha\}_\alpha$ of 
$\omega$, the foliation $\mathcal F_\omega$ is comprised of minimal hypersurfaces,
\item the deviation of the leaves inside the disks from the minimality is \break
$\sim\mu^{d -1}$-small. 
\end{itemize}

Therefore, for a given intrinsically harmonic map  $f: M^3 \rightarrow S^1$ and 
any  collection of \emph{non-singular} fibers $\{F\}$, whose closure does 
not contain the $f$-singularities, there is a 
family of metrics $g_\mu$ on $M^3$, such that $f$ is harmonic and the 
surfaces $F$'s are minimal. Furthermore, the argument in [K] shows 
that the fibers of $f: \partial M^3 \rightarrow S^1$ are geodesic loops. 
The boundary $\partial M$ is \emph{sufficiently convex} in the sense of Meeks and Yau 
[MY] (in fact, it is flat).

According to Theorem 5.1 in [FHS], for a large class of 3-manifolds 
$M$ (so called, $P^2$-irreducible ones), any two-sided \emph{incompressible} 
oriented embedded surface $\alpha: \Sigma ' \subset M$, distinct from $S^2$, can 
be isotoped to 
1) a \emph{minimal embedded} surface $\Sigma \subset M$, or to 
2) a surface 
$\Sigma \subset M$, which is a boundary of a regular neighborhood of a 
\emph{one-sided minimal} embedded surface.
  
In the first case, the minimal surface minimizes the area in the homotopy class of $\alpha$.    
In the second case, the minimal surface realizes \emph{half} of the minimal area in the 
homotopy class of $\alpha$.\smallskip

Combining this result with the special  properties of the metrics $\{g_\mu\}$, we conclude 
that any embedded incompressible surface in the homology class of a fiber has a minimal or 
a "near-minimal" representative $\Sigma \subset M$, having only \emph{hyperbolic} 
tangencies with the $f$-fibers \emph{outside} of the disks $D_\mu(x_\alpha)$'s. In other 
words, all the elliptic tangent points must be located \emph{inside} the disks 
$\{D_\mu(x_\alpha)\}$. As $\mu \rightarrow 0$, they are attracted towards the singularities 
of $f$.
\smallskip

\section{notation list}
\begin{itemize}
\item $var_{\chi_-}(f)$ \quad  cf. formulas (3.1) and (4.1)

\item $Var_{\chi_-}(f)$ \quad  cf. formula (4.2)

\item $\rho(\Sigma, F)$ \quad cf. the $\rm 7^{th}$  paragraph in Section 6

\item $\rho^\circ(\Sigma, F)$ \quad cf. the $\rm 2^{nd}$ paragraph before Definition 6.1

\item $\rho_{\chi_-}(\Sigma, F)$ \quad cf. Definition 6.1 

\item $\rho_{\chi_-}([\Sigma], F_R)$ \quad cf. Definition 6.3

\item $\rho_{\chi_-}(f) := \rho_{\chi_-}([F], F_R)$,\; $F$ being a fiber

\item $\rho(\hat\Sigma, \tilde M)$ \quad cf. formula (6.1)

\item $h(\Sigma, F)$ \quad cf. formula (6.2)

\item $b(F, \Sigma)$ \; cf. Definition 6.7

\item $\mu^\circ(\Sigma, F)$ \; cf. formula (7.1)

\item $b_{\chi_-}(F_R, k[F])$ \quad cf. the paragraph before Corollary 6.9

\item $h_{\chi_-}(k[F], f)$ \quad cf. Definition 6.11 

\item $h_{\chi_-}(f) := h_{\chi_-}([F], f)$,\;  $F$ being a fiber

\item $\mu^\circ([\Sigma], F_R)$ \quad cf. the $\rm 1^{st}$ paragraph after Corollary 8.7 

\item $\|[\Sigma]\|_T$ \quad cf. the $\rm 3^{nd}$ paragraph in the Introduction

\item $\|[\Sigma]\|_{H^f}$ \quad cf. formula (8.2) and Definition 8.1

\item $\|\kappa\|_A$ \quad cf. formula (8.1)
 
\end{itemize}

\end{document}